\documentclass[notitlepage,leqno,10pt]{article}
\usepackage{latexsym}
\usepackage{amsmath}
\usepackage{amsfonts}
\usepackage{amssymb}
\usepackage{amscd}

\renewcommand{\d}{\delta }
\newcommand{\D }{\Delta }

\renewcommand{\l }{\lambda }
\renewcommand{\L }{\Lambda }
\newcommand{\m }{\mu }
\newcommand{\n }{\nabla }

\newcommand{\Sig }{\Sigma}

\newcommand{\ov}{\overline}
\newcommand{\intbar}{\mathop{\int\makebox(-13.5,0){\rule[4pt]{.7em}{0.3pt}}%
\kern-6pt}\nolimits}

\newcommand{\be}{\begin{equation}}
\newcommand{\ee}{\end{equation}}
\newcommand{\bes}{\begin{equation*}}
\newcommand{\ees}{\end{equation*}}
\newcommand{\ba}{\begin{eqnarray}}
\newcommand{\ea}{\end{eqnarray}}
\newcommand{\bas}{\begin{eqnarray*}}
\newcommand{\eas}{\end{eqnarray*}}
\newenvironment{pf}{\noindent{\sc Proof}.\enspace}{\rule{2mm}{2mm}\medskip}
\newenvironment{pfn}{\noindent{\sc Proof}}{\rule{2mm}{2mm}\medskip}

\newcommand{\R}{\mathbb{R}}

\newcommand{\Z}{\mathbb{Z}}

\newcommand{\N}{\mathbb{N}}

\author{ Cheikh Birahim NDIAYE$^a$ \;\& \;Mohameden OULD AHMEDOU$^b$ }

\date{}

\title{\bf Morse theory and the resonant $Q$-curvature problem}
\begin{document}

\newtheorem{lem}{Lemma}[section]
\newtheorem{pro}[lem]{Proposition}
\newtheorem{thm}[lem]{Theorem}
\newtheorem{rem}[lem]{Remark}
\newtheorem{cor}[lem]{Corollary}
\newtheorem{df}[lem]{Definition}

\maketitle

\begin{center}

{\small

\noindent  $^{a, b}$ Mathematisches Institut der Justus-Liebig-Universit\"at Giessen, \\Arndtstrasse 2, D-35392 Giessen, Germany.

}
{\small

\noindent

}
\end{center}

\footnotetext[1]{E-mail addresses: ndiaye@everest.mathematik.uni-tuebingen.de, Mohameden.Ahmedou@math.uni-giessen.de}

\

\

\begin{center}
{\bf Abstract}
\end{center}
In this paper, we study the prescribed $Q$-curvature problem on closed four-dimensional Riemannian manifolds when the total integral of the $Q$-curvature is a positive integer multiple of the one of the four-dimensional round sphere. This problem has a variational structure with a lack of compactness. Using some topological tools of the theory of { \it critical points at infinity}   combined with a refined blow-up analysis and some dynamical arguments, we identify the accumulations points of all noncompact flow lines of a pseudogradient flow, the so called critical points at infinity of the associated variational problem, and associate to them a natural Morse index. We then prove strong Morse type inequalities, extending the full Morse theory to this noncompact variational problem. Finally, we derive from our results Poincar\'e-Hopf index type criteria for existence, extending known results in the literature and deriving new ones.

\begin{center}


\bigskip

\centerline{\bf AMS subject classification: 53C21, 35C60, 58J60}

\end{center}


\section{Introduction and statement of the results}
One of the most important problem in conformal geometry is the problem of finding conformal metrics with a prescribed curvature quantity. An example of curvature quantity which has received a lot of attention in the last decades is the Branson's $Q$-curvature. It is a Riemannian scalar invariant introduced by Branson\cite{bran1} for closed four-dimensional Riemannian manifolds. Indeed, given $(M, g)$ a four-dimensional Riemannian manifold, with Ricci tensor $Ric_{g}$, scalar curvature $R_{g}$, and Lapalce-Beltrami operator $\D_g$, the $Q$-curvature of $(M, g)$ is defined by
\begin{equation}\label{eq:P}
Q_{g}=-\frac{1}{12}(\D_{g}R_{g}-R_{g}^{2}+3|Ric_{g}|^{2}).
\end{equation}
Under the conformal change of metric $g_u=e^{2u}g$ with $u$ a smooth function on $M$, the $Q$-curvature transforms in the following way
\begin{equation}\label{eq:P1}
P_gu+2Q_g=2Q_{g_{u}}e^{4u},
\end{equation}
where $P_g$ is the Paneitz operator introduced by Paneitz\cite{p1} in 1983 and is defined by the following formula
\ba\label{eq:P2}
P_{g}\varphi=\D_{g}^{2}\varphi+div_{g}\left((\frac{2}{3}R_{g}g-2Ric_{g})\n_g\varphi\right),
\ea
where\;$\varphi$\; is any smooth function on \;$M$ and $\n_g$ denotes the covariant derivative with respect to $g$. When, one changes conformally $g$ as before, namely by $g_u=e^{2u}g$ with $u$ a smooth function on $M$, $P_g$ obeys the following simple transformation law
\begin{equation}\label{eq:tranlaw}
P_{g_u}=e^{-4u}P_{g}.
\end{equation}
The equation \eqref{eq:P1} and the formula \eqref{eq:tranlaw} are analogous to classical ones which hold on closed Riemannian surfaces. Indeed, given a closed Riemannian surface $(\Sig, g)$ and $g_u=e^{2u}g$ a conformal change of $g$ with $u$ a smooth function on $M$, it is well know that
\ba\label{eq:g}
\Delta_{g_u}=e^{-2u}\Delta_{g},\;\;\;\;\;\;\;\;\;-\Delta_{g}u+K_{g}=K_{g_u}e^{2u},
\ea
where for a background metric $\tilde g$ on $\Sigma$ , $\Delta_{\tilde g}$\;and\;$K_{\tilde g}$ are respectively the Laplace-Beltrami operator and the Gauss curvature of ($\Sigma, \tilde g$). In addition to these,  we have an analogy with the classical Gauss-Bonnet formula
$$
\int_{\Sigma}K_{g}dV_{g}=2\pi\chi(\Sigma),
$$
where\;$\chi(\Sigma)$\; is the Euler characteristic of \;$\Sigma$. In fact, we have the Chern-Gauss-Bonnet formula
$$
\int_{M}(Q_{g}+\frac{|W_{g}|^{2}}{8})dV_{g}=4\pi^{2}\chi(M),
$$
where \;$W_{g}$\;denotes the Weyl tensor of \;($M,g$) and $\chi(M)$\; is the Euler characteristic of \;$M$. Hence, from the pointwise conformal invariance of \;$|W_{g}|^{2}dV_{g}$,\;it follows that \;$\int_{M}Q_{g}dV_g$\; is also conformally invariant and will be denoted by $\kappa_P$, namely
\begin{equation}\label{eq:kappa}
\kappa_{P}=\kappa_{P}[g]:=\int_{M}Q_{g}dV_g.
\end{equation}
\vspace{6pt}

\noindent
Like for the  Kazdan-Warner problem, given a smooth positive function $K$ defined on a closed four-dimensional Riemannian manifold $(M, g)$, one can ask under which condition on $K$, there exists a metric conformal to $g$ and having $K$ as $Q$-curvature. Thanks to $\eqref{eq:tranlaw}$, the problem is equivalent to finding a smooth solution of the equation
\ba\label{eq:qequation}
P_{g}u+2Q_{g}=2Ke^{4u}\;\;\;\;\;\;in\;\;M.
\ea
On the other hand, by the regularity result of Uhlenbeck-Viaclovsky\cite{uv},  smooth solutions of equation \eqref{eq:qequation} can be found as critical points of the following geometric functional:
\begin{equation}\label{eq:defj1}
J(u):= \, <P_g u,u> \, + \, 4 \int_M Q_g u dV_g \,  - \,\kappa_{P} \log\left(\int_{M}Ke^{4u} dV_g\right),\;\;\;\;u\in W^{2, 2}(M),
\end{equation}
where
$$
 <P_g u,v>=\int_M\D_gu\D_gvdV_g+\frac{2}{3}\int_{M}R_g\n_g u\cdot\n_gvdV_g-\int_M2Ric_g(\n_g u, \n_gv)dV_g,\;\;\;u,\;v\in W^{2, 2}(M),
$$
 and $W^{2, 2}(M)$ is the space of functions on $M$ which are of class $L^2$, together with their first and second derivatives.

 \vspace{4pt}

\noindent
The analytic features of equation $\eqref{eq:qequation}$ and of the associated Euler-Lagrange functional $J$ depend strongly on the conformal invariant $\kappa_P$. Indeed, depending on whether $\kappa_P$ is a positive integer multiple of $8\pi^2$ or not, the noncompactness of equation $\eqref{eq:qequation}$ and the way of finding critical points of $J$ changes drastically. As far as existence results are concerned, we have that problem\;$\eqref{eq:qequation}$\;has been solved in a celebrated work of Chang-Yang$\cite{cy}$ under the assumption that $ker P_g\simeq \R$,\;$P_{g}$\;is a non-negative operator and $\kappa_P<8\pi^{2}$. The same result has also been reproved and extended to arbitrary even dimensions (regarding high-order analogues of the Paneitz operator and the $Q$-curvature, see \cite{fg} and \cite{gjms}) by Brendle\cite{bren} using geometric flow methods. Later Djadli and Malchiodi \cite{dm} used a delicate min-max scheme to extend the above result of Chang-Yang\cite{cy} and Brendle\cite{bren}. More precisely they showed that problem $\eqref{eq:qequation}$ is solvable provided, $ker P_g \simeq \R$ and $\kappa_{P} \notin 8\pi^2 \N^{*}$. The latter result of Djadli-Malchiodi\cite{dm} has been extended to arbitrary dimensions by the first author, see \cite{nd1}, regarding high-order analogues of the Paneitz operator and $Q$-curvature (see \cite{fg} and \cite{gjms}) .
\vspace{4pt}

\noindent
The assumptions $ker P_g\simeq \R$ and $\kappa_P\notin 8\pi^2\N^{*}$ will be referred to as {\em nonresonant} case. This terminology is motivated by the fact that in that situation the set of solutions to some perturbations of equation $\eqref{eq:qequation}$ (including it) is compact, see \cite{dr}, \cite{mal}, \cite{mart1}, and \cite{nd1}.  Naturally, we call {\em resonant} case when $ker P_g\simeq \R$ and $\kappa_P\in 8\pi^2\N^{*}$. With these terminologies, we have that the works of Chang-Yang\cite{cy} and Djadli-Malchiodi\cite{dm} answer affirmatively the question raised  above in the {\em nonresonant} case. However, for the {\em resonant} case, the only existence results known to the best of our knowledge deal with the case $\kappa_P=8\pi^2$ and are of two types. The first type of existence results is based on violation of strong Morse inequalities at "infinity" and is due to Wei-Xu\cite{wx} and Malchiodi-Struwe \cite{MS} and is for $(S^4, g_{S^4})$ as background Riemannian manifold.  While the second type of existence results is based on a positive mass type assumption and is due to Li-Li-Liu\cite{lll}.

\vspace{4pt}

\noindent
In this work, we are interested in the {\em resonant} case, namely when $ker P_g\simeq \R$ and $\kappa_P=8\pi^2m$ with $m\in\N^*$. To present the main results of the paper, we need to set some notation and make some definitions. We define \;$\mathcal{F}_K: M^m\setminus F_m(M)\longrightarrow \R$\; as follows
\begin{equation}\label{eq:limitfs}
\mathcal{F}_K(a_1,\cdots, a_m):=2\sum_{i=1}^m\left(H(a_i, a_i)+\sum_{j\neq i}G(a_i, a_j)+\frac{1}{2}\log(K(a_i))\right),
\end{equation}
where $F(M^m)$ denotes the fat Diagonal of $M^m$, namely $F(M^m):=\{A:=(a_1, \cdots, a_m)\in M^m:\;\;\text{there exists} \;\;i\neq j\;\,\text{with}\;\, a_i=a_j\}$, $G$ is the Green's function of $P_g(\cdot)+\frac{2}{m}Q_g$ satisfying the  normalization $\int_M Q_g(x) G(\cdot, x)dV_g(x)=0$, and $H$ is its regular part, see section \ref{eq:notpre} for more information. Furthermore, we define
\begin{equation}\label{eq:critfk}
Crit(\mathcal{F}_K):=\{A\in M^m\setminus F_m(M), \;\;A\;\;\;\text{critical point of} \;\;\mathcal{F}_K\}.
\end{equation}
Moreover,  for $A=(a_1,\cdots, a_m)\in M^m\setminus F_m(M)$, we set
\begin{equation}\label{eq:partiallimit}
\mathcal{F}^A_i(x):=e^{4(H(a_i, x)+\sum_{j\neq i}G(a_j, x))+\frac{1}{4}\log(K(x))},
\end{equation}
and define
\begin{equation}\label{eq:defindexa}
\mathcal{L}_K(A):=-\sum_{i=1}^m L_g(\sqrt{\mathcal{F}^{A}_i})(a_i),
\end{equation}
where $L_g:=-\D_g+\frac{1}{6}R_g$ is the conformal Laplacian associated to $g$.  We set also
\begin{equation}\label{eq:critsett}
\mathcal{F}_{\infty}:=\{A\in Crit(\mathcal{F}_K)\:\;\mathcal{L}_K(A)<0\},
\end{equation}
\begin{equation}\label{eq:minf}
i_{\infty}(A):=5m-1-Morse(A, \mathcal{F}_K),
\end{equation}
and define
\begin{equation}\label{eq:mi}
m_i^m:=card\{A\in Crit(\mathcal{F}_K): i_{\infty}(A)=i\}, i=0, \cdots, 5m-1,
\end{equation}
where $Morse(\mathcal{F}_K, A)$ denotes the Morse index of $\mathcal{F}_K$ at  $A$. We would like to point out that for $m\geq 2$, $m_i^m=0$ for  $0\leq i\leq m-2$. For $m\geq 2$, we use the notation $B_{m-1}(M)$ to denote the set of formal barycenters of order $m-1$ of $M$, namely
\begin{equation}\label{eq:baryp}
B_{m-1}(M):=\{\sum_{i=1}^{m-1}\alpha_i\d_{a_i}, a_i\in M, \alpha_i\geq 0, i=1,\cdots, m-1,\;\,\sum_{i=1}^{m-1}\alpha_i=1\},
\end{equation}
and for more informations, see \cite{dm} and \cite{kk}. Furthermore, we define
\begin{equation}\label{eq:cpm}
c^{m-1}_p=dim \;H_p(B_{m-1}(M)), \;\;p=1, \cdots 5m-6,
\end{equation}
where $H_p(B_{m-1}(M)$ denotes the $p$-th homology group of $B_{m-1}(M)$ with $\Z_2$ coefficients. Finally, we say
\begin{equation}\label{eq:nondeg}
(ND) \;\;\;\text{holds if } \;\;\mathcal{F}_K\;\text{is a Morse function and for every}\;\;A\in Crit(\mathcal{F}_K)  \;\;\mathcal{L}_K(A)\neq 0.
\end{equation}
\vspace{4pt}

\noindent
Now, we are ready to state our main results starting with the {\em critical} case, namely when $m=1$.
\begin{thm}\label{eq:morsepoincare1}
Let $(M, g)$ be a closed four-dimensional Riemannian manifold such that $Ker P_g\simeq \R$ and $\kappa_P=8\pi^2$. Assuming that $K$ is a smooth positive function on $M$ such that $(ND)$ holds and the following system
\begin{equation}\label{eq:mp1}
\begin{cases}
m^1_0=1+k_0,\\
m^1_i=k_i+k_{i-1}, \;&i=1, \cdots, 4,\\
0=k_4\\
k_i\geq 0,\;\;& i=0, \cdots, 4
\end{cases}
\end{equation}
has no solutions, then $K$ is the $Q$-curvature of a metric conformal to $g$.
\end{thm}
\vspace{4pt}

\noindent
As a corollary of the above theorem, we derive the following Poincar\'e-Hopf index type criterion for existence.
\begin{cor}\label{eq:existence1}
Let $(M, g)$ be a closed four-dimensional Riemannian manifold such that $Ker P_g\simeq \R$ and $\kappa_P=8\pi^2$. Assuming that $K$ is a smooth positive function on $M$ such that $(ND)$ holds and
\begin{equation}\label{eq:ep1}
\sum_{A\in \mathcal{F}_{\infty}} (-1)^{i_{\infty}(A)}\neq 1,
\end{equation}
 then $K$ is the $Q$-curvature of a metric conformal to $g$.
\end{cor}
\vspace{4pt}

\noindent
We notice that Corollary \ref{eq:existence1} extends a previous Poincar\'e-Hopf index criterium for existence proved for the four-dimensional round sphere, to the case of arbitrary closed four-dimensional Riemannian manifold with Paneitz operator of trivial kernel and total integral of $Q$-curvature equal to $8\pi^2$, see \cite{wx}. We remark also that Theorem \ref{eq:morsepoincare1} extends a previous result of Malchiodi-Struwe\cite{MS}, to the case of arbitrary closed four-dimensional Riemannian manifold with Paneitz operator of trivial kernel and total integral of $Q$-curvature equal to $8\pi^2$. 
\vspace{4pt}

\noindent
As a by product of our analysis, we extend the above Poincar\'e Hopf criterium to include the case where the total sum equals 1 but a partial one is not, provided that there is a jump in the Morse indices of the function $\mathcal{F}_K$. Namely we have the following criterium:
\begin{thm}\label{t:C}
Let $(M, g)$ be a closed four-dimensional Riemannian manifold such that $Ker P_g\simeq \R$, $\kappa_P=8\pi^2$ and  $K$ be  a smooth positive function on $M$ satisfying the non degeneracy  condition $(ND)$. Assuming that  there exists a positive integer \;$1 \leq l \leq 4$\; such that 
\begin{equation*}\label{eq:ep1c}
\begin{split}
\sum_{A\in \mathcal{F}_{\infty},\;  i_{\infty}(A) \leq l -1 } &(-1)^{i_{\infty}(A)}\neq 1\\&\text{and}\\
 \forall A \in \mathcal{F}_{\infty},\;\; &\quad i_{\infty}(A) \neq l, 
\end{split}
\end{equation*}
then $K$ is the $Q$-curvature of a metric conformal to $g$.
\end{thm}
\vspace{4pt}

\noindent
Next, we state our result in the {\em supercritical} case, namely when $m\geq 2$. It read as follows:
\begin{thm}\label{eq:morsepoincare2}
Let $(M, g)$ be a closed four-dimensional Riemannian manifold such that $Ker P_g\simeq \R$ and $\kappa_P=8\pi^2m$, $m\in \N$ and $m\geq 2$. Assuming that $K$ is a smooth positive function on $M$ such that $(ND)$ holds and the following system
\begin{equation}\label{eq:mp3}
\begin{cases}
0=k_0,\\
m^m_1=k_1,\\
m_i^m=c^{m-1}_{i-1}+k_i+k_{i-1}, \;&i=2,\cdots, 5m-5,\\
m_{i}^m=k_{i}+k_{i-1},\;& i=5m-4,\cdots, 5m-1,\\
0=k_{5m-1},\\
k_i\geq 0,\;\;& i=0, \cdots, 5m-1,
\end{cases}
\end{equation}
has no solutions, then $K$ is the $Q$-curvature of a metric conformal to $g$.
\end{thm}
\vspace{4pt}

\noindent
As a corollary of Theorem \ref{eq:morsepoincare2}, we derive the following Poincar\'e-Hopf index type criterion for existence.
\begin{cor}\label{eq:existence2}
Let $(M, g)$ be a closed four-dimensional Riemannian manifold such that $Ker P_g\simeq \R$ and $\kappa_P=8\pi^2m$ with $m\in \N$ and $m\geq 2$. Assuming that $K$ is a smooth positive function on $M$ such that $(ND)$ holds and 
\begin{equation}\label{eq:ep2}
\sum_{A\in \mathcal{F}_{\infty}} (-1)^{i_{\infty}(A)}\neq \frac{1}{(m-1)!}\Pi_{i=1}^{m-1}(i-\chi(M)),
\end{equation}
 then $K$ is the $Q$-curvature of a metric conformal to $g$.
\end{cor}
\vspace{6pt}
Just as in the one mass case, we generalize the above criterium to the case where there is a jump in the indices of the function $\mathcal{F}_K$. Namely we prove the following:
\begin{thm}\label{eq:Cm}
Let $(M, g)$ be a closed four-dimensional Riemannian manifold such that $Ker P_g\simeq \R$, $\kappa_P=8\pi^2m$ with $m\in \N$ and $m\geq 2$, and  $K$ be  a smooth positive function on $M$ satisfying the non degeneracy condition $(ND)$.  Assuming that there exists a positive integer \;$1 \leq l \leq 5m-1$\;  and \;$A^l\in \mathcal{F}_{\infty}$ with \;$i_{\infty}(A^l) \leq l - 1$ such that
\begin{equation*}\label{eq:ep2c}
\begin{split}
\sum_{A\in \mathcal{F}_{\infty}, \; i_{\infty}(A) \leq l - 1} (-1)^{i_{\infty}(A)}&\neq \frac{1}{(m-1)!}\Pi_{j=1}^{m-1}(j-\chi(M))\\
&\text{and}\\
\forall A  \in \mathcal{F}_{\infty},\;\;  &\qquad i_{ \infty}(A)  \neq l,
\end{split}
\end{equation*}
 then $K$ is the $Q$-curvature of a metric conformal to $g$.
\end{thm}
\vspace{6pt}
\noindent
\begin{rem}\label{eq:remtheo}
We would like to point out that, Theorem \ref{eq:morsepoincare1}  implies directly Corollary \ref{eq:existence1} (see section \ref{eq:proofresults}), since \eqref{eq:mp1} not having a solution is equivalent to strong Morse inequalities  at "infinity" do not hold while \eqref{eq:ep1} is just a violation of Poincar\'e-Hopf index formula at "infinity". Similarly also and for the same reason, we have Theorem \ref{eq:morsepoincare2} implies directly Corollary \ref{eq:existence2}.
\end{rem}

\begin{rem}\label{eq:negativeimpact}
We would like to point out that in Theorem \ref{eq:morsepoincare1}- Theorem \ref{eq:Cm}, the presence of negative eigenvalues of the Paneitz operator $P_g$ plays no role. In fact, this is not surprising, since Theorem \ref{eq:morsepoincare1}- Theorem \ref{eq:Cm} are based on the topological contribution of the  critical points at infinity of $J$ which are due to the involved bubbling phenomena, and not the loss of coercivity of $J$, since even Palais-Smale sequences have bounded negative component, where the negative component of a $W^{2, 2}$-function is  defined by \eqref{eq:defnegativec}.
\end{rem}

\begin{rem}\label{eq:relililiu}
We remark that the estimates established in this work and the dynamical argument of section \ref{eq:morseinfinity} is used in a forthcoming paper to sharpen Proposition \ref{eq:escape} and to prove existence results in the spirit of the one of Li-Li-Liu\cite{lll}. Furthermore, it will be used to establish a compactness theorem for equation \eqref{eq:qequation} in the full resonant case when $(ND)$ holds and to give an interpretation of \eqref{eq:ep1} (resp. \eqref{eq:ep2}) as the non vanishing of the Leray-Schauder degree of \eqref{eq:qequation} when $(ND)$ holds and $m=1$ (resp. $m\geq 2$).
\end{rem}
\vspace{6pt}

\noindent
The strategy of the proof of our results goes along the ideas and techniques introduced by Bahri\cite{bah} (see also \cite{bc}, \cite{bd} and \cite{bcch}). A first difficulty in implementing such techniques, is the analysis of the sequences violating the Palais-Smale condition, which seems to be out of reach at the moment. To circumvent such a difficulty, we make use of a deformation Lemma due to Lucia\cite{lu}. It has the advantage of replacing the analysis of arbitrary Palais-Smale sequences to the one of some appropriate viscosity solutions, which are Palais-Smale ones, but with a prescribed rate of decay to zero for the gradient of $J$ along them. A second difficulty arises in the finite-dimensional reduction. Indeed, it turns out to be a difficult task to get the needed positive definiteness of second variation of the Euler-Lagrange functional $J$ when working with the usual neighborhood of potential critical points at infinity. A fact which jeopardizes the finite-dimensional Lyapunov-Schmidt reduction. To remedy to such a problem, we use a delicate blow-up analysis of the lack of compactness along the flow lines of Lucia's pseudogradient, to reduce the size of the usual neighborhood of potential critical points at infinity, and perform a Morse type reduction of $J$ in such a neighborhood. Next, using some refined estimates of the Euler-Lagrange functional $J$ and its gradient at infinity combined with the construction of a delicate pseudogradient at infinity, we identify all the critical points at infinity of our variational problem, and associate to them a kind of Morse index, which corresponds to the dimension of their unstable manifolds. We then extend the full Morse theory to this noncompact variational problem and derive the associated strong Morse inequalities. The existence criteria in Theorem \ref{eq:morsepoincare1}, Corollary \ref{eq:existence1},  Theorem \ref{eq:morsepoincare2} and Corollary \ref{eq:existence2} are just violation of the strong Morse inequalities quoted previously. While the existence criteria in Theorem \ref{t:C} and Theorem \ref{eq:Cm} are violation of a suitable  Poincar\'e-Hopf index formula.
\vspace{8pt}

\noindent
The structure of the paper is as follows. In section \ref{eq:notpre}, we collect some notation and preliminary results, like the Green's function  $G$ of $P_g(\cdot)+\frac{2}{m}Q_g$ satisfying the normalization $\int_MG(\cdot, y)Q_g(y)dV_g(y)=0$ and define a family of approximate solutions. In section \ref{eq:ninfty}, we introduce the notion of neighborhoods of potential critical points at infinity and derive a refinement of Lucia's deformation Lemma. In section \ref{eq:jexp}, we perform a first expansion of $J$ at infinity. In section \ref{eq:gjesp}, we present a useful expansion of the gradient of $J$ at infinity. Section \ref{eq:findim} is concerned with reducing variationally the problem to a finite dimensional one. In section \ref{eq:topsublevel}, we study the topology of some appropriate sublevels of $J$. In section \ref{eq:morseleminf}, we present a Morse Lemma at infinity for $J$ by constructing a suitable pseudogradient for the finite dimension reduction functional. Section \ref{eq:proofresults} is concerned with the proofs of Theorem \ref{eq:morsepoincare1}-Theorem \ref{eq:Cm}. Finally, in section \ref{eq:appendix}, we collect some technical Lemmas.
\vspace{4pt}

\noindent
\begin{center}
{\bf Acknowledgements}
\end{center}
The authors started working on this problem when C. B. Ndiaye was still  at the university of T\"ubingen.  Part of this paper was written when C. B. Ndiaye was visiting the department of mathematics of the university of Rutgers in USA. He is very grateful to all these institutions for support and good working conditions, and to DFG Sonderforschungsbereich TR 71 Freiburg-T\"ubingen. The authors are also grateful to Pr. Abbas Bahri for  useful discussions. The authors have been supported  by the DFG project "Fourth-order uniformization type theorems for $4$-dimensional Riemannian manifolds".
\section{Notation and preliminaries}\label{eq:notpre}
In this brief section, we fix our notation, and give some preliminaries. First of all, we recall that $(M, g)$  and $K$ are respectively the given underlying closed four-dimensional Riemannian manifold and the prescribed function  with the following properties
\begin{equation}\label{eq:trivialker}
ker P_g\simeq \R\;\;\text{ and} \;\;\kappa_P=8m\pi^2\;\; \text {for some integer }\;\;m\geq 1, \;\;\text{and}\;\;K\;\;\text{is a smooth positive function on}\;M.
\end{equation}

\noindent
We are going to discuss the asymptotics near the singularity of the Green's function $G$ of the operator \;$P_{g}(\cdot)+\frac{2}{m}Q_g$\; satisfying the normalization $\int_MG(\cdot, y)Q_g(y)dV_g(y)=0$ \; and make some related definitions.
\vspace{6pt}

\noindent
In the following, for a Riemmanian metric $\bar g$ on $M$, we will use the notation\;$B^{\bar g}_{p}(r)$\; to denote the geodesic ball with respect to $\bar g$ of radius \;$r$\;and center \;$p$. We also denote by \;$d_{\bar g}(x,y)$\; the geodesic distance with respect to $\bar g$ between two points \;$x$\;and \;$y$\; of \;$M$, $exp_x^{\bar g}$ the exponential map with respect to $\bar g$ at $x$. $inj_{\bar g}(M)$\;stands for the injectivity radius of \;$(M, \bar g)$, $dV_{\bar g}$\;denotes the Riemannian measure associated to the metric\;$\bar g$.  Furthermore, we recall that $\n_{\bar g}$, \;$\D_{\bar g}$, \;$R_{\bar g}$\; and \;$Ric_{\bar g}$\; will denote respectively the covariant derivative, the Laplace-Beltrami operator, the scalar curvature and Ricci curvature with respect to \;$\bar g$. For simplicity, we will use the notation $B_p(r)$ to denote $B^g_{p}(r)$, namely $B_p(r)=B^g_p(r)$. $M^2$\;stands for the cartesian product \;$M\times M$, while \;$Diag(M)$\; is the diagonal of \;$M^2$.
\vspace{4pt}

\noindent
$\N$\;denotes the set of non-negative integers, $\N^*$\;stands for the set of positive integers, $\R_+$ the set of positive real numbers, $\bar \R_+$ the set of non-negative real numbers, $(\R^4, g_{\R^4})$ the  $4$-dimensional Euclidean space, $\D_4$ the Euclidean Laplacian on $\R^4$ and $(S^4, g_{S^4})$ the $4$-dimensional round sphere. Moreover, for $r>0$, $B^{4}_0(r)$ denotes the $4$-dimensional Euclidean ball of center $0$ and radius $r$. $M^m$ denotes the cartesian product of $m$ copies of $M$.
\vspace{4pt}

\noindent
For $1\leq p\leq \infty$ and $k\in \N$, $\theta\in  ]0, 1[$, $L^p(M)$, $W^{k, p}(M)$, $C^k(M)$, and $C^{k, \theta} (M)$ stand respectively for the standard Lebesgue space, Sobolev space, $k$-continuously differentiable space and $k$-continuously differential space of H\"older exponent $\beta$, all with respect $g$ (if the definition needs a metric structure, and for precise definitions and properties, see \cite{aubin} or \cite{gt}). Given a function \;$u\in L^1(M)$,\;$\bar u_{Q}$\; denotes its average on \;$M$ with respect to the sign measure $Q_gdV_g$, that is
\begin{equation}\label{eq:qva}
\ov{u}_Q=\frac{1}{8m\pi^2}\int_{M} Q_g(x)u(x)dV_{g}(x).
\end{equation}
Given a generic Riemannian metric $\tilde g$ on $M$ and a function \;$F(x, y)$\; defined on a open subset of\;$M^2$ which is symmetric and\; with $F(\cdot, \cdot)\in C^2$ with respect to $\tilde g$, we define $\frac{\partial F(a, a)}{\partial a}:=\frac{\partial F(x, a)}{\partial x}|_{x=a}=\frac{\partial F(a, y)}{\partial y}|_{y=a}=$, and $\D_{\tilde g} F(a_1, a_2):=\D_{\tilde g, x}F(x, a_2)|_{x=a_1}=\D_{\tilde g, y}F(a_2, y)|_{y=a_1}.$
\vspace{6pt}

 \noindent
 For  $\epsilon>0$ and small,  $\l\in \R_+$, $\l\geq \frac{1}{\epsilon}$, and $a\in M$, $O_{\l, \epsilon}(1)$ stands for quantities bounded uniformly in \;$\l$, and $\epsilon$, and \;$O_{a, \epsilon}(1)$ stands for quantities bounded uniformly in $a$ and $\epsilon$. For  $l\in \N^*$, $O_{l}(1)$ stands for quantities bounded uniformly in \;$l$\; and \;$o_l(1)$ stands for quantities which tends to $0$ as $l\rightarrow +\infty$.  For  $\epsilon$ positive and small, \;$a\in M$\; and \;$\l\in \R_+$ large, $\l\geq \frac{1}{\epsilon}$,\;$O_{a, \l, \epsilon}(1)$\; stands for quantities bounded uniformly in \;$a$, \;$\l$, and $\epsilon$. For $\epsilon$ positive and small, $p\in \N^{*}$, $\bar \l:=(\l_1, \cdots, \l_p)\in (\R_+)^p$, $\l_i\geq \frac{1}{\epsilon}$  for $i=1, \cdots, p$, and $A:=(a_1, \cdots, a_p)\in M^p$ (where $ (\R_+)^p$ and $M^p$ denotes respectively the cartesian product of $p$ copies of $\R_+$ and $M$), $O_{A, \bar \l, \epsilon}(1)$ stands for quantities bounded uniformly in $A$, $\bar \l$, and $\epsilon$. Similarly for $\epsilon $ positive and small,  $p\in \N^{*}$, $\bar \l:=(\l_1, \cdots, \l_p)\in (\R_+)^p$, $\l_i\geq \frac{1}{\epsilon}$ for $i=1, \cdots, p$, $\bar \alpha:=(\alpha_1, \cdots, \alpha_p)\in \R^p$, $\alpha_i$ close to $1$ for $i=1, \cdots, p$, and $A:=(a_1, \cdots, a_p)\in M^p$ (where $ \R^p$  denotes the cartesian product of $p$ copies of $\R$, $O_{\bar\alpha, A, \bar \l, \epsilon}(1)$ will mean quantities bounded from above and below independent of $\bar \alpha$, $A$, $\bar \l$, and $\epsilon$. For $x\in \R$, we will use the notation $O(x)$ to mean $|x|O(1)$ where $O(1)$ will be specified in all the contexts where it is used. Large positive constants are usually denoted by $C$ and the value of\;$C$\;is allowed to vary from formula to formula and also within the same line. Similarly small positive constants are also denoted by $c$ and their value may varies from formula to formula and also within the same line.
\vspace{6pt}

\noindent
 Now, for $X$ a topological space, $(X, A)$ a pair of spaces and $q\in \N$, we denote respectively by $H_q(X)$, $H_{q}(X, A)$,  the $q$-th homology group of $X$, and the $q$-th homology group of $(X, A)$, all with  respect to $\Z_2$.  Furthermore, we use the the notation $b_q(X)$ and $b_q(X, A)$ to denote respectively the $q$-th betti number  of $X$ and $(X, A)$. We will use $\chi(X)$ and $\chi(X, A)$ to denote respectively the Euler characteristic of $X$ and $(X, A)$.
\vspace{6pt}

\noindent
We call $\bar m$ the number of negative eigenvalues (counted with multiplicity) of $P_g$. We point out that $\bar m$ can be zero, but it is always finite. If $\bar m\geq 1$, then we will denote by \;$V\subset W^{2, 2}(M)$ the direct sum of the eigenspaces corresponding to the negative eigenvalues of $P_g$. The dimension of $V$ is of course $\bar m$. On the other hand, we have the existence of an $L^2$-orthonormal basis of eigenfunctions $v_1,\cdots, v_{\bar m}$ of $V$ satisfying
\begin{equation}\label{eq:defeigen1}
P_g v_i=\mu_r v_r\;\;\;\forall\;\;r=1\cdots \bar m,
\end{equation}
\begin{equation}\label{eq:defeigen2}
\mu_1\leq \mu_2\leq \cdots\leq \mu_{\bar m}<0<\mu_{\bar m+1}\leq\cdots,
\end{equation}
where $\mu_r$'s are the eigenvalues of $P_g$ counted with multiplicity. We define also the pseudo-differential operator $P^+_g$ as follows
\begin{equation}\label{eq:operatorrev}
P_g^+u=P_gu-2\sum_{r=1}^{\bar m}\m_r(\int_M u v_rdV_g) v_r.
\end{equation}
Basically $P_g^+$ is obtained from $P_g$ by reversing the sign of the negative eigenvalues and we extend the latter definition to $\bar m=0$ for uniformity in the analysis and recall that in that case $P_g^+=P_g$.  Now, using $P_g^+$, we set for $t>0$
\begin{equation}\label{eq:jt}
J_t(u):=<P_g^+u, u>+2t\sum_{r=1}^{\bar m}\mu_r(u^{r})^2+4t\int_MQ_gudV_g-8\pi^2tm\log\int_MKe^{4u}dV_g, \,\;\;u\in W^{2, 2}(M),
\end{equation}
with
\begin{equation}\label{eq:defnegativecr}
u^{r}:=\int_M uv_rdV_g,\;\;\; r=1, \cdots, m, u\in W^{2, 2}(M).
\end{equation}
Now, using \eqref{eq:operatorrev} and \eqref{eq:jt}, we obtain
\begin{equation}\label{eq:jt1}
J_t(u):=<P_gu, u>+2(t-1)\sum_{r=1}^{\bar m}\mu_r(u^{r})^2+4t\int_MQ_gudV_g-8\pi^2tm\log\int_MKe^{4u}dV_g, \,\;\;u\in W^{2, 2}(M),
\end{equation}
and hence $J=J_1$. Furthermore, using \eqref{eq:defnegativecr}, we define
\begin{equation}\label{eq:defnegativec}
u^{-}=\sum_{r=1}^{\bar m}u^{r}v_r.
\end{equation}
\vspace{6pt}

\noindent
We will use the notation $<\cdot, \cdot>$ to denote the $L^2$ scalar product and $<\cdot, \cdot>_{W^{2, 2}}$ for the $W^{2, 2}$ scalar product. On the other hand, it is easy to see that
\begin{equation}\label{eq:productpr}
<u, v>_{P}:=<P^{+}_gu, v>, \;\,\;u, v\in\{w\in W^{2, 2}(M):\;\;\;u_Q=0\}
\end{equation}
defines a inner product on $\{u\in W^{2, 2}(M):\;\;\;u_Q=0\}$ which induces a norm equivalent to the standard $W^{2, 2}(M)$ (on $\{u\in W^{2, 2}(M):\;\;\;u_Q=0\}$) and denoted by
\begin{equation}\label{eq:normpr}
||u||:=\sqrt{<u, u>_{P}} \;\;\;u\in\{w\in W^{2, 2}(M):\;\;\;u_Q=0\}.
\end{equation}
As above, in the general case, namely $\bar m\geq 0$, for $\epsilon$ small and positive, $\bar \beta:=(\beta_1, \cdots, \beta_{\bar m})\in \R^{\bar m}$  with $\beta_i$ close to $0$, $i=1, \cdots, \bar m$) (where $\R^{\bar m}$ is the empty set when $\bar m=0$), $\bar \l:=(\l_1, \cdots, \l_p)\in (\R_+)^p$, $\l_i\geq \frac{1}{\epsilon}$  for $i=1, \cdots, p$, $\bar \alpha:=(\alpha_1, \cdots, \alpha_p)\in \R^p$, $\alpha_i$ close to $1$ for $i=1, \cdots, p$, and $A:=(a_1, \cdots, a_p)\in M^p$, $p\in \N^{*}$, $w\in W^{2, 2}$ with $||w||$ small, $O_{\bar\alpha, A, \bar \l, \bar \beta, \epsilon}(1)$ will stand quantities bounded independent of $\bar \alpha$, $A$, $\bar \l$, $\bar \beta$, and $\epsilon$, and $O_{\bar\alpha, A, \bar \l, \bar \beta, w, \epsilon}(1)$ will stand quantities bounded independent of $\bar \alpha$, $A$, $\bar \l$, $\bar \beta$,  $w$ and $\epsilon$.
\vspace{6pt}

\noindent
Given a point $b\in \R^4$ and $\lambda$ a positive real number, we define $\delta_{b, \lambda}$ to be the {\em standard bubble}, namely
\begin{equation}\label{eq:standarbubble}
\delta_{b, \lambda}(y):=\log\left(\frac{2\lambda}{1+\lambda^2|y-b|^2}\right),\;\;\;\;\;\;y\in \R^4.
\end{equation}
The functions $\delta_{b, \lambda}$ verify the following equation
\begin{equation}\label{eq:bubbleequation}
\D^2_4\delta_{b,\lambda}=6e^{4\delta_{b,\lambda}}\;\;\;\text{in}\;\;\;\R^4.
\end{equation}
Geometrically, equation $\eqref{eq:bubbleequation}$ means that the metric $g=e^{2\delta_{b, \lambda}} dx^2$ (after pull-back by stereographic projection) has constant $Q$-curvature equal to $3$, where $dx^2$ is the standard metric on $\R^4$.\\
Using the existence of conformal normal coordinates (see \cite{cao} or \cite{gun}), we have that, for $a \in M$\; there exists a function $u_a\in C^{\infty}(M)$ such that
\begin{equation}\label{eq:detga}
g_a = e^{2u_a} g\;\; \text{verifies}\;\;det g_a(x)=1\;\;\text{for}\;\;\; x\in B^{g_a}_a( \varrho_a).
\end{equation}
with $0<\varrho_a<\frac{inj_{g_a}(M)}{10}$. Moreover, we can take the families of functions \;$u_a$, $g_a$ and $\varrho_a$ such that
\begin{equation}\label{eq:varro0}
\text{the maps}\;\;\;a\longrightarrow u_a, \;g_a\;\;\text{are}\;\;C^1\;\;\;\text{and}\;\;\;\;\varrho_a\geq \varrho_0>0,
\end{equation}
for some small positive $\varrho_0$ satisfying $\varrho_0<\frac{inj_g(M)}{10}$, and
\begin{equation}\label{eq:proua}
\begin{split}
&||u_a||_{C^4(M)}=O_a(1),\;\;\frac{1}{\ov C^2} g\leq g_a\leq \ov C^2 g, \\\;&u_a(x)= O_a(d^2_{g_a}(a, x))=O_a(d_{g}^2(a, x)) \;\;\text{for}\;\;x\in\;\;B_a^{g_a}(\varrho_0)\supset B_a(\frac{\varrho_0}{2\ov C}),\;\;\text{and}\\&
u_a(a)=0,\;\;\;R_{g_a}(a)=0, \;\;\;\n_{g_a} R_{g_a}(a)=0,
\end{split}
\end{equation}
for some large positive constant $\ov C$ independent of $a$. For the meaning of $O_a(1)$ in \eqref{eq:proua}, see section \ref{eq:notpre}. Furthermore, using the scalar curvature equation satisfied by $e^{-u_a}$, namely $-\D_{g_a} (e^{-u_a})+\frac{1}{6}R_{g_a}e^{-u_a}=\frac{1}{6}R_g(a)e^{-3u_a}$\; in $M$, and \eqref{eq:detga}-\eqref{eq:proua}, it is easy to see that the following holds
\begin{equation}\label{eq:lapemua}
\D_{g_a}(e^{-u_{a}})(a)=-\frac{1}{6}R_g(a).
\end{equation}
For $a\in M$, and $r>0$, we set
\begin{equation}\label{eq:expballaga}
exp_a^{a}:=exp_{a}^{g_a}\;\;\;\;\text{and}\;\;\;B_a^{a}(r):=B_a^{g_a}(r).
\end{equation}
On the other hand, using the properties of $g_a$ (see \eqref{eq:detga}-\eqref{eq:proua})), it is easy to check that for every $u\in C^2(B_a^{a}(\eta))$ with $0<\eta<\varrho_0$ there holds
\begin{equation}\label{eq:gradlapga}
\n_{g_a}u(a)=\n_gu(a)=\n_4\hat u(0),\;\;\;\;\text{and}\;\;\;\D_{g_a}u(a)=\D_4\hat u(0),
\end{equation}
where
\begin{equation}\label{eq:hatu}
\hat u(y)=u(exp_a^{a}(y)),\;\;y\in B_0^4(\eta).
\end{equation}
Now, for $0<\varrho<\min\{\frac{inj_g(M)}{4}, \frac{\varrho_0}{4}\}$ where $\varrho_0$ is as in \eqref{eq:varro0}, we define a cut-off function $\chi_{\varrho} : \bar\R_+ \rightarrow \bar\R_+$ satisfying the following properties:
\begin{equation}\label{eq:cutoff}
\begin{cases}
\chi_{\varrho}(t)  = t \;\;&\text{ for } \;\;t \in [0,\varrho],\\
\chi_{\varrho}(t) = 2 \varrho \;\;&\text{ for } \;\; t \geq 2 \varrho, \\
 \chi_{\varrho}(t) \in [\varrho, 2 \varrho] \;\;\;&\text{ for } \;\; t\in [\varrho, 2 \varrho].
\end{cases}
\end{equation}
Using the cut-off function $\chi_{\varrho}$, we define for $a\in M$ and $\lambda\in \R_+$  the function $\hat{\delta}_{a, \lambda}$ as follows
\begin{equation}\label{eq:hatdelta}
\hat{\delta}_{a, \lambda}(x):=\log\left(\frac{2\lambda}{1+\lambda^2\chi_{\varrho}^2(d_{g_a}(x, a))}\right).
\end{equation}
For every $a\in M$ and $\lambda\in \R_+$, we define $\varphi_{a, \lambda}$ to be the unique the solution of the following projection problem
\begin{equation}\label{eq:projbubble}
\begin{cases}
P_g \varphi_{a,\l} \, + \, \frac{2}{m} Q_g\, = 16\pi^2\,\frac{ e^{4 (\hat{\d}_{a,\l} \,  +  \, u_a)}}{\int_M e^{4 (\hat{\d}_{a,\l} \,  +  \, u_a)}dV_g}\;\; \mbox{ in } \;\;M, \\
\int_M Q_g(x)\varphi_{a,\l}(x) \,dV_g(x) = 0.
\end{cases}
\end{equation}
So differentiating with respect to $\l$ and $a$ (respectively) the relation $\int_MQ_g(x)\varphi_{a, \l}(x)dV_g(x)=0$, we get (respectively)
\begin{equation}\label{eq:normlambda}
\int_M Q_g(x)\frac{\partial \varphi_{a, \l}(x)}{\partial \l}dV_g(x)=0,
\end{equation}
and
\begin{equation}\label{eq:norma}
\int_M Q_g(x)\frac{\partial \varphi_{a, \l}(x)}{\partial a}dV_g(x)=0.
\end{equation}
Now, we recall that $G$ is the unique solution of the following PDE
\begin{equation}\label{eq:defG4}
\begin{cases}
P_g G(a, \cdot)+\frac{2}{m}Q_g(\cdot)=16\pi^2\delta_a(\cdot),\\
\int_M Q_g(x) G(a, x)dV_g(x)=0.
\end{cases}
\end{equation}
Using  \eqref{eq:defG4}, it is easy to see that the following integral representation formula holds
\begin{equation}\label{eq:G4integral}
u(x)-\ov{u}_Q=\frac{1}{16\pi^2}\int_M G(x, y)P_gu(y), \;\;\;u\in C^4(M), \;x\in M,
\end{equation}
where $u_Q$ is defined as in section \ref{eq:notpre}. It is a well know fact that $G$ has a logarithmic singularity. In fact $G$ decomposes as follows
\begin{equation}\label{eq:decompG4}
 G(a,x)=\log \left(\frac{1}{\chi_{\varrho}^2(d_{{g}_a}(a, x))}\right)+H(a, x).
\end{equation}
where $H$ is the regular par of $G$, see for example \cite{zw}. Furthermore, it is also a well-know fact that (see  still \cite{zw})
\begin{equation}\label{eq:regH4}
G\in C^{\infty}(M^2\setminus Diag(M)),\;\;\;H\in W^{4, p}(M^2)\;\;\forall 1\leq p<\infty, \,\;\text{hence} \;H\in C^{3, \beta}(M^2)\;\;\;\forall \beta\in (0, 1).
\end{equation}
Now, using \eqref{eq:limitfs} and \eqref{eq:partiallimit} combined with the symmetry of $H$, it is easy to see that
\begin{equation}\label{eq:relationderivative}
\frac{\partial \mathcal{F}(a_1, \cdots, a_m)}{\partial a_i}=\frac{\n_g\mathcal{F}^{A}_i(a_i)}{\mathcal{F}^{A}_i(a_i)}, \;\;\;i=1, \cdots, m.
\end{equation}
Next, we set
\begin{equation}\label{eq:deflA}
l_K(A):=\sum_{i=1}^m\left(\frac{\D_{g_{a_i}} \mathcal{F}^{A}_i(a_i)}{\sqrt{\mathcal{F}^{A}_i(a_i)}}-\frac{2}{3}R_g(a_i)\sqrt{\mathcal{F}^{A}_i(a_i)}\right),
\end{equation}
and use the properties of the metrics $g_{a_i}$, the transformation rule of the conformal Laplacian under conformal change of metrics and direct calculations, to get
\begin{equation}\label{eq:auxiindexa1}
l_K(A)=2\mathcal{L}_K(A), \;\;\forall A\in Crit(\mathcal{F}_K).
\end{equation}
Finally, for  $m\geq 2$ and $\bar m\geq 1$, we define
\begin{equation}\label{eq:defbarynega}
A_{m-1, \bar m}:=\widetilde{B_{m-1}(M)\times B^{\bar m}_1},
\end{equation}
 where the equivalence relation \;$\widetilde{}$\; means that $B_{m-1}(M)\times \partial B^{\bar m}_1$ is identified with $\partial B_1^{\bar m}$, and for more informations, see \cite{dm}.
On the other hand, covering the set $A_{m-1, \bar m}$ with the two sets  $A=B_{m-1}(M)\times B_{\frac{3}{4}}^{\bar m}$ and	$B=B_{m-1}(M)\times (B_{1}^{\bar m}\setminus B_{\frac{1}{4}}^{\bar m})$, and using the exactness of the Mayer-Vietoris sequence and the K\"unneth theorem, one can easily see that the following holds (see for example \cite{dm}, \cite{dem}, and \cite{maldeg})
\begin{equation}\label{eq:relationbary}
H_{q}(A_{m-1, \bar m})\simeq 0,\,\;\text{for}\;\; 1\leq q\leq\bar m\;\;\text{and}\;\;H_q(A_{m-1, \bar m})\simeq  H_{q-\bar m}(B_{m-1}(M))\;\;\text{for}\;\;q\geq \bar m+1.
\end{equation}
\section{Neighborhood of potential critical points at infinity}\label{eq:ninfty}
In this section, we introduce the neighborhood of potential critical points at infinity and derive a deformation Lemma which takes into account the possible noncompactness phenomena. We start with the following deformation Lemma which is due to Lucia\cite{lu} and inspired from the works of Bahri \cite{bah2}, \cite{bah1} on contact forms. We anticipate that its proof is exactly the same as in \cite{lu} and for this reason we will not present it here.
\begin{lem}\label{eq:deformlem} (Lucia deformation Lemma)\\
 Assuming that $a, b\in \R$ such that $a<b$ and  there is no critical values of $J$ in $[a, b]$, then there are two possibilities\\
1) Either  $$J^a\;\text{is a deformation retract of}\;\; J^b.$$
2) Or there exists a sequence $t_l\rightarrow 1$ as $l\rightarrow +\infty$ and a sequence of critical point $u_l$ of $J_{t_l}$ verifying $a\leq J(u_l)\leq b$ for all $l\in \N^*$.
\end{lem}
\vspace{6pt}
\begin{df}
Following A. Bahri we call  "critical point at infinity" an accumulation point of a sequence of the type given by the second alternative of the above Lemma and which imped the deformation of the sublevel $J^b$ onto the sublevel $J^a$, for some real numbers  $a, b \in \R$ with $a<b$.
\end{df}
\noindent
Clearly, Lemma \ref{eq:deformlem} implies that to understand the variational analysis of $J$, one has to understand the asymptotic behaviour of blowing-up sequence of the type given by alternative 2) of Lemma \ref{eq:deformlem}. In order to do that, we first fix $\L$ to  a large positive constant and $\gamma \in (0, \frac{1}{2})$. Next, for $\epsilon$ and $\eta$ small positive real numbers, we denote by $V(m, \epsilon, \eta)$ the $(m,\epsilon, \eta)$-{\em neighborhood of potential critical points at infinity}, namely
\begin{equation}\label{eq:ninfinity}
\begin{split}
V(m, \epsilon, \eta):=\{u\in W^{2, 2}(M):a_1, \cdots, a_m\in M, \;\;\l_1,\cdots, \l_m>0, \;\;||u-\ov{u}_Q-\sum_{i=1}^m\varphi_{a_i,\lambda_i}||+\\ ||\n^{W^{2,2}} J(u)||=O\left( \sum_{i=1}^m\frac{1}{\l_i^{1-\gamma}}\right)\;\;\;\lambda_i\geq \frac{1}{\epsilon},\;\;\frac{2}{\Lambda}\leq \frac{\l_i}{\l_j}\leq \frac{\Lambda}{2},\;\;\text{and}\;\;d_g(a_i, a_j)\geq 4\ov C\eta\;\;\text{for}\;i\neq j\},
\end{split}
\end{equation}
where $\ov C$ is as in \eqref{eq:proua},  $\n^{W^{2, 2}}J$ is the gradient of $J$ with respect to the $W^{2, 2}$-topology, $O(1):=O_{A, \bar \l, u, \epsilon}(1)$ meaning bounded uniformly in  $\bar\l:=(\l_1, \cdots, \l_n)$, $A:=(a_1, \cdots, m)$, $u$, $\epsilon$.
\vspace{6pt}

\noindent
Now, having introduced the sets $V(m, \epsilon, \eta)$, we are ready to characterize blowing-up sequence of the type given by alternative 2) of Lemma \ref{eq:deformlem}. Indeed, we have:
\begin{pro}\label{eq:escape}
Let $\epsilon$ and $\eta$ be small positive real numbers with $0<2\eta<\varrho$ where $\varrho$ is as in \eqref{eq:cutoff}. Assuming that $u_l$ is a sequence of blowing up critical point  of $J_{t_l}$ with $\ov{(u_l)}_Q=0, l\in \N$ and $t_l\rightarrow 1$ as $l\rightarrow +\infty$ , then  there exists $l_{\epsilon, \eta}$ a large positive integer such that for every $l\geq l_{\epsilon, \eta}$, we have $u_l\in V(m, \epsilon, \eta)$.
\end{pro}
\begin{pf}
Since $u_l\in Crit (J_{t_l})$, then using \eqref{eq:jt1}, and setting
\begin{equation}\label{eq:dfvl}
v_l=u_l-\frac{1}{4}\log\left(\int_M Ke^{4u}dV_g\right)+\frac{1}{4}\log \left(8\pi^2m\right),
\end{equation}
we derive that $v_l$ solves the following fourth order elliptic equation
\begin{equation}\label{eq:equvl}
P_g v_l+2(t_l-1)\sum_{r=1}^{\bar m}\mu_r v_l^{r}v_r+2t_lQ_g=2 t_lKe^{4v_l}.
\end{equation}
Now, testing equation \eqref{eq:equvl} with $v_l^{-}$ (for its definition see \eqref{eq:defnegativecr}), integrating \eqref{eq:equvl} and using the fact that $P_g$ annihilates constant and its eigenfunctions has zero mean value with respect to $dV_g$, we get
\begin{equation}\label{eq:quantization}
\int_M K e^{4v_l}dV_g=\int_MQ_gdV_g=8\pi^2m,
\end{equation}
and that $v_l^{-}$ verifies the following estimate
\begin{equation}\label{eq:estvlmo1}
||v_l^{-}||=O_l(1).
\end{equation}
Thus, setting $Q_l:=\sum_{r=1}^{\bar m}\mu_r v_l^{r}v_r$,  and recalling that the $v_r$'s verify $\int_M v_rdV_g=0$, we infer that
\begin{equation}\label{eq:defql}
\int_M Q_ldV_g=0, \;\;\;l\in N,
\end{equation}
and  up to a subsequence
\begin{equation}\label{eq:conql}
Q_l\longrightarrow Q_0\;\;\text{in}\;\;C^2(M), \;\;l\longrightarrow +\infty..
\end{equation}
Now, combining \eqref{eq:equvl} and \eqref{eq:defql}, we obtain
\begin{equation}\label{eq:equvl1}
P_g v_l+2(t_l-1)Q_l+2t_lQ_g=2 t_lKe^{4v_l}.
\end{equation}
On the other hand, using the definition of $v_l$ (see \eqref{eq:dfvl}) and the fact that $\ov{(u_l)}_{Q}=0$, we obtain
\begin{equation}\label{eq:relavlul}
v_l-\ov{(v_l)}_{Q}=u_l.
\end{equation}
Now, using \eqref{eq:quantization}, \eqref{eq:relavlul}, and the assumption $u_l$ blows-up, we deduce that $v_l$ blows-up. Thus, since  $\int_M Ke^{4v_l}dV_g=8m\pi^2$, $m\in \N^*$, then from the work of Druet-Robert\cite{dr} and Malchiodi\cite{mal} (see \cite{mart1}, \cite{nd1}), we have that  (up to a subsequence) there exists $m$ sequence of real numbers $\l_{i, l}$, $i=1, \cdots, m$, $l\in \N$ and $m$ sequence of converging points $x_{i, l}$ $i=1, \cdots, m$, $l\in \N$ such that
\begin{equation}\label{eq:infopl}
\begin{split}
d_g(x_{i, l}, x_{j, l})\geq 4\ov C\eta, \;\;\;i\neq j=1, \cdots, m,\; l\in \N\\
\frac{2}{\Lambda}\leq \frac{\l_{i, l}}{\l_{j, l}}\leq \frac{\Lambda}{2},\;\;\;i, j=1, \cdots, m, \;l\in \N\\
\l_{i, l}\geq \frac{1}{\epsilon},\;\;\;i=1, \cdots, m, \;l\in \N.
\end{split}
\end{equation}
Now, to continue, we set
\begin{equation}\label{eq:defwl}
w_l=v_l-\ov{(v_l)}_Q-\sum_{i=1}^m\varphi_{x_{i, l}, \l_{i, l}}, \;\;\;l\in \N.
\end{equation}
Thus, combining \eqref{eq:projbubble} and \eqref{eq:defwl}, we infer that
\begin{equation}\label{eq:meanqwl}
\ov{(w_l)}_Q=0,\;\;l\in \N.
\end{equation}
Next, using \eqref{eq:projbubble} and \eqref{eq:equvl1}, we derive that $w_l$ solve the following fourth order elliptic equation
\begin{equation}\label{eq:equieqwl}
\begin{split}
P_g w_l+2(t_l-1)Q_l+2(t_l-1)Q_g=&6e^{4(v_l+\frac{1}{4}\log( \frac{t_lK}{3}))}-\sum_{i=1}^me^{4(\hat\d_{x_{i, l}, \l_{i, l}}+u_{x_{i, l}})}\\&+\sum_{i=1}^m\left(6-\frac{16\pi^2}{\int_M e^{4(\hat\d_{x_{i, l}, \l_{i, l}}+u_{x_{i, l}})}dV_g}\right)e^{4(\hat\d_{x_{i, l}, \l_{i, l}}+u_{x_{i, l}})}, \;\;l\in \N.
\end{split}
\end{equation}
On the other hand, since $\gamma\in (0, 1)$, then there exists $p_1>1$ and close to $1$ and $r_1>0$ and close to $0$ such that the following holds
\begin{equation}\label{eq:p1r}
1-\gamma+16p_1-16+r_1\in (0, 1)\;\;\;\text{and}\;\;\;\gamma-16p_1+16\in (0, 1).
\end{equation}
Thus, we have
\begin{equation}\label{eq:nu}
\nu:=1-\gamma+16p_1-16+r_1\in (0, 1).
\end{equation}
Now, using the Weinstein-Zhang estimate with \eqref{eq:nu},  and the properties of $u_{x_{i, l}}$, we  infer that
\begin{equation}\label{eq:wza1}
v_l(\cdot)+\frac{1}{4}\log\frac{t_lK(x_{i, l})}{3}=\hat \d_{x_{i, l}, \l_{i, l}}(\cdot)+u_{x_{i, l}}+O\left( (d_g(\cdot, x_{i, l}))^{\nu}\right)\;\;\;B_{x_{i, l}}(\eta).
\end{equation}
Next, we are going to estimate $\ov{(v_l)}_Q$. In order to do that, we first use a classical argument (see for example \cite{mal}, page 16, proof of formula (88)), and obtain
\begin{equation}\label{eq:estvlq}
e^{4\ov{(v_l)}_Q}\leq C\sum_{i=1}^m\frac{1}{\l_{i, l}^{4-\tilde \epsilon}},
\end{equation}
where $\tilde \epsilon$ is small and positive, and $C$ is a large positive constant , both independent of $l$. On the other hand, using  a classical argument based on the Green's representation formula \eqref{eq:G4integral} and suitable application of Jensens's inequality and Fubini theorem (see for example \cite{mal}, proof of Proposition 3.1), we derive that for every $q\geq 1$, the following estimate holds
\begin{equation}\label{eq:estout}
||e^{4v_l}||_{L^q(M\setminus\cup_{1=1}^m B_{x_{i, l}}(\eta))}\leq C e^{4\ov{(v_l)}_Q},
\end{equation}
where $C$ is a large positive constant independent of $l$.  Thus, combining \eqref{eq:estvlq} and \eqref{eq:estout}, we obtain that for every $q\geq 1$, the following estimate holds
\begin{equation}\label{eq:estoutvlq}
||e^{4v_l}||_{L^q(M\setminus\cup_{1=1}^m B_{x_{i, l}}(\eta))}\leq  C\sum_{i=1}^m\frac{1}{\l_{i, l}^{4-\tilde \epsilon}},
\end{equation}
where $C$ is a large positive constant independent of $l$. To continue, we set
\begin{equation}\label{eq:el}
E_l=e^{4(v_l+\frac{1}{4}\log( \frac{t_lK}{3}))}-\sum_{i=1}^me^{4(\hat\d_{x_{i, l}, \l_{i, l}}+u_{x_{i, l}})}
\end{equation}
and
\begin{equation}\label{eq:fl}
F_l=\sum_{i=1}^m\left(6-\frac{16\pi^2}{\int_M e^{4(\hat\d_{x_{i, l}, \l_{i, l}}+u_{x_{i, l}})}dV_g}\right)e^{4(\hat\d_{x_{i, l}, \l_{i, l}}+u_{x_{i, l}})}.
\end{equation}
So , combining \eqref{eq:equieqwl}, \eqref{eq:el} and \eqref{eq:fl}, we get
\begin{equation}\label{eq:equieqwl1}
P_gw_l+2(t_l-1)Q_l+2(t_l-1)Q_g=6E_l+F_l.
\end{equation}
Now, let us estimate $F_l$ and after $E_l$. For this, we first use  the explicit expression of $\hat\d_{x_{i, l}, \l_{i, l}}$ and the properties of $u_{x_{i, l}}$,  and have that for every $p\geq 1$, there holds
\begin{equation}\label{eq:linfel1}
||F_l||_{L^p(M)}\leq C\sum_{i=1}^m\frac{1}{\l_{i, l}^{\frac{4}{p}}},
\end{equation}
where $C$ is a large positive constant independent of $l$. Now, we are going to estimate $E_l$. In order to do that, we first set
\begin{equation}\label{eq:e1il}
E^{1, i}_l:=\left(e^{4(v_l+\frac{1}{4}\log(\frac{t_lK(x_{i, l}}{3}))}-e^{4(\hat\d_{x_{i, l}, \l_{i, l}}+u_{x_{i, l}})}\right)\chi_{B_{x_{i, l}}(\eta)},\;\;\;i=1, \cdots, m,
\end{equation}
\begin{equation}\label{eq:e2il}
E^{2, i}_l:=\left(e^{4(v_l+\frac{1}{4}\log(\frac{t_lK}{3}))}-e^{4(v_l+\frac{1}{4}\log(\frac{t_lK(x_{i, l}}{3}))}\right)\chi_{B_{x_{i, l}}(\eta)},\;\;\;i=1, \cdots, m.
\end{equation}
Furthermore, we define
\begin{equation}\label{eq:e3l}
E^{3}_l:=e^{4(v_l+\frac{1}{4}\log(\frac{t_lK}{3}))}\chi_{(M\setminus \cup_{i=1}^m B_{x_{i, l}}(\eta))}-\sum_{i=1}^me^{4(\hat\d_{x_{i, l}, \l_{i, l}}+u_{x_{i, l}})}\chi_{(M\setminus B_{x_{i, l}}(\eta))}.
\end{equation}
Moreover, we set
\begin{equation}\label{eq:e1e2l}
E^1_l=\sum_{1=1}^mE_l^{1, i}\;\;\;\text{and}\;\;E^2_l=\sum_{i=1}^mE^{2, i}_l.
\end{equation}
On the other hand, using \eqref{eq:el}, \eqref{eq:e1il}-\eqref{eq:e1e2l}, it is easy to see that
\begin{equation}\label{eq:equiel}
E_l=E^1_l+E^2_l+E^3_l.
\end{equation}
Now, using again the explicit expression of $\hat\d_{x_{i, l}, \l_{i, l}}$ and the properties of $u_{x_{i, l}}$ combined with \eqref{eq:estoutvlq}, we infer that  for every $q\geq 1$, there holds
\begin{equation}\label{eq:este3l}
||E^3_l||_{L^q(M)}\leq C\sum_{i=1}^m\frac{1}{\l_{i, l}^3},
\end{equation}
where $C$ is a large positive constant independent of $l$. Next, using again the explicit expression of $\hat\d_{x_{i, l}, \l_{i, l}}$, the properties of $u_{x_{i, l}}$, normal coordinate at $x_{i, l}$ with respect to $g_{x_{i, l}}$ and \eqref{eq:wza1}, we obtain that for every $i=1, \cdots, m$ and for every $p\in (1, \frac{4}{4-\nu})$, the following estimate holds
\begin{equation}\label{eq:e1il}
||E^{1, i}_l||_{L^p(M)}^p=\int_{B_{x_{i, l}}(\eta)}|E^{1, i}_ldV_g\leq C\l_{i, l}^{4p-4-p\nu}\int_{B_0(\l_{i, l}\eta)}\left(\frac{2}{1+|y|^2}\right)^{4p}|y|^{p\nu}dy\leq C\l_{i, l}^{4p-4-p\nu},
\end{equation}
where $C$ is a large positive constant independent of $l$. Hence, combining \eqref{eq:e1e2l} and \eqref{eq:e1il}, we derive that for every $p\in (1, \frac{4}{4-\nu})$, there holds
\begin{equation}\label{eq:e1l}
||E^{1}_l||_{L^p(M)}\leq C\sum_{i=1}^m\frac{1}{\l_{i, l}^{\frac{4-4p+p\nu}{p}}},
\end{equation}
where $C$ is a large positive constant independent of $l$. Now, let us estimate $E_l^2$. For this, we use the mean value Theorem, the explicit expression of $\hat\d_{x_{i, l}, \l_{i, l}}$, the properties of $u_{x_{i, l}}$, normal coordinate at $x_{i, l}$ with respect to $g_{x_{i, l}}$ and \eqref{eq:wza1}, we obtain that for every $i=1, \cdots, m$ and for every $p\in (1, \frac{4}{3})$, the following estimate holds
\begin{equation}\label{eq:e2il}
||E^{2, i}_l||_{L^p(M)}^p=\int_{B_{x_{i, l}}(\eta)}|E^{2, i}_ldV_g\leq C\l_{i, l}^{3p-4}\int_{B_0(\l_{i, l}\eta)}\left(\frac{2}{1+|y|^2}\right)^{4p}|y|^{p}dy\leq C\l_{i, l}^{3p-1},
\end{equation}
where $C$ is a large positive constant independent of $l$. Thus, combining \eqref{eq:e1e2l} and \eqref{eq:e2il}, we derive that for every $p\in (1, \frac{4}{3})$, there holds
\begin{equation}\label{eq:e2l}
||E^{2}_l||_{L^p(M)}\leq C\sum_{i=1}^m\frac{1}{\l_{i, l}^{\frac{4-3p}{p}}},
\end{equation}
where $C$ is a large positive constant independent of $l$. On the other hand, setting
\begin{equation}\label{eq:defp}
p:=\frac{4}{4-r_1},
\end{equation}
and using the fact that $\nu=1-(\gamma-16p_1+16)+r$, $r_1\in (0, 1)$ and $\gamma-16p_1+16\in (0, 1)$, we have that
\begin{equation}\label{eq:goodr}
p\in (1, \frac{4}{3})\;\;\;\;\text{and}\;\;\;\;p\in (1, \frac{4}{4-\nu}).
\end{equation}
Furthermore, using again the fact that $\nu=1-(\gamma-16p_1+16)+r_1$ combined with \eqref{eq:defp}, we obtain
\begin{equation}\label{eq:realexp}
\frac{4-4p+p\nu}{p}=1-\gamma+16p_1-16 \;\;\;\;\text{and}\;\;\;\;\frac{4p-3}{p}=1+\frac{3}{4}r_1.
\end{equation}
Thus, combining \eqref{eq:e1l}, \eqref{eq:e2l}, \eqref{eq:goodr} and \eqref{eq:realexp},  we get
\begin{equation}\label{eq:e1lgp}
||E^{1}_l||_{L^p(M)}\leq C\sum_{i=1}^m\frac{1}{\l_{i, l}^{1-\gamma+16p_1-16}},
\end{equation}
and
\begin{equation}\label{eq:e2lgp}
||E^{2}_l||_{L^p(M)}\leq C\sum_{i=1}^m\frac{1}{\l_{i, l}^{1+\frac{3r_1}{4}}},
\end{equation}
where $C$ is a large positive constant independent of $l$. Hence, combining \eqref{eq:equiel}, \eqref{eq:este3l}, \eqref{eq:e1lgp}, and \eqref{eq:e2lgp}, we obtain
\begin{equation}\label{eq:estel}
||E_l||_{L^p(M)}\leq C\sum_{i=1}^m\frac{1}{\l_{i, l}^{1-\gamma+16p_1-16}},
\end{equation}
where $C$ is a large positive constant independent of $l$. Next, we set
\begin{equation}\label{eq:tfl}
\tilde F_l=6E_l+F_l,
\end{equation}
and   use \eqref{eq:equieqwl} to have
\begin{equation}\label{eq:eqtfl}
P_g w_l+2(t_l-1)Q_l+2(t_l-1)Q_g=\tilde F_l.
\end{equation}
On the other hand, using \eqref{eq:linfel1}  and \eqref{eq:defp}, we obtain
\begin{equation}\label{eq:linfel2}
||F_l||_{L^{p}(M)}\leq C\sum_{i=1}^m\frac{1}{\l_{i, l}^{4-r_1}}.
\end{equation}
 Thus, combining \eqref{eq:linfel2}, \eqref{eq:estel}, \eqref{eq:tfl},  and recalling that $r_1\in (0, 1)$, we get
\begin{equation}\label{eq:esttfl}
||\tilde F_l||_{L^p(M)}\leq C\sum_{i=1}^m\frac{1}{\l_{i, l}^{1-\gamma+16p_1-16}}.
\end{equation}
Now, integration \eqref{eq:eqtfl} and using \eqref{eq:defql}, the assumption $\int_M Q_gdV_g=8\pi^2m$, and \eqref{eq:esttfl}, we derive that the following estimate holds
\begin{equation}\label{eq:esttl}
|t_l-1|\leq C\sum_{i=1}^m\frac{1}{\l_{i, l}^{1-\gamma+16p_1-16}},
\end{equation}
where $C$ is large positive constant independent of $l$. Now, recalling that $\ov{(w_l)}_Q=0$, $l\in \N $ (see \eqref{eq:meanqwl})) and using \eqref{eq:eqtfl}, \eqref{eq:esttfl}, and \eqref{eq:esttl}, we have thanks to standard elliptic regularity theory  that the following estimate holds
\begin{equation}\label{eq:w4pwl}
||w_l||_{W^{4, p}(M)}\leq \sum_{i=1}^m\frac{1}{\l_{i, l}^{1-\gamma+16p_1-16}},
\end{equation}
where $C$ is a large positive constant independent of $l$. Thus, recalling that $p>1$ and using Sobolev embedding theorem, we infer that
\begin{equation}\label{eq:w22wl}
||w_l||_{W^{2, 2}(M)}\leq C\sum_{i=1}^m\frac{1}{\l_{i, l}^{1-\gamma+16p_1-16}},
\end{equation}
where $C$ is a large positive constant independent of $l$. On the other hand, since $\ov{(w_l)}_Q=0$, then using \eqref{eq:normpr} and recalling the discussions before it, we have that \eqref{eq:w22wl} implies that
\begin{equation}\label{eq:w22wl1}
||w_l||\leq C\sum_{i=1}^m\frac{1}{\l_{i, l}^{1-\gamma+16p_1-16}},
\end{equation}
where $C$ is a large positive constant independent of $l$. Next, recalling that $p_1>1$ and using \eqref{eq:defwl} and \eqref{eq:w22wl1}, we infer the following weaker estimate
\begin{equation}\label{eq:estvlw22a}
||v_l-\ov{(v_l)}_Q-\sum_{i=1}^m\varphi_{x_{i, l}, \l_{i, l}}||\leq C\sum_{i=1}^m\frac{1}{\l_{i, l}^{1-\gamma}},
\end{equation}
where $C$ is a large positive constant independent of $l$. Thus, recalling that $\ov {(u_l)}_Q=0$ and using \eqref{eq:dfvl}, we infer that \eqref{eq:estvlw22a} implies
\begin{equation}\label{eq:estvlw22}
||u_l-\ov{(u_l)}_Q-\sum_{i=1}^m\varphi_{x_{i, l}, \l_{i, l}}||\leq C\sum_{i=1}^m\frac{1}{\l_{i, l}^{1-\gamma}},
\end{equation}
where $C$ is a large positive constant independent of $l$.
\vspace{2pt}

\noindent
Now, we are going to estimate $||\n^{W^{2, 2}}J(u)||_{W^{2, 2}(M)}$ and for this we will need the stronger estimate \eqref{eq:w22wl} on $w_l$. First of all, since $u_l\in Crit(J_{t_l})$, using the definition of $J_{t_l}$ (see \eqref{eq:jt} with $t$ replaced by $t_l$),  and \eqref{eq:defwl}, we obtain
\begin{equation}\label{eq:derbu2}
\begin{split}
&<\n^{W^{2, 2}} J(\sum_{i=1}^m\varphi_{x_{i, l}, \l_{i, l}}), \varphi>_{W^{2, 2}(M)}=-2<P_gw_l, \varphi>+4(1-t_l)\int_MQ_g\varphi dV_g\\&+2(1-t_l)\sum_{i=1}^{\bar m}\mu_iu^{i}\varphi^{i}+4\times 32m\pi^2(1-t_l)\frac{\int_M Ke^{4\sum_{i=1}^m\varphi_{x_{i, l}, \l_{i, l}}}\varphi dV_g}{\int_M Ke^{4\sum_{i=1}^m\varphi_{x_{i, l}, \l_{i, l}}}dV_g}\\+&4\times 32mt_l\pi^2\left(\frac{\int_M Ke^{4\sum_{i=1}^m\varphi_{x_{i, l}, \l_{i, l}}}\varphi dV_g}{\int_M Ke^{4\sum_{i=1}^m\varphi_{x_{i, l}, \l_{i, l}}}dV_g}-\frac{\int_M Ke^{4(v_l-\ov{(v_l)}_Q)}\varphi dV_g}{\int_M Ke^{4(v_l-\ov{(v_l)}_Q)}dV_g}\right).
\end{split}
\end{equation}
On the other hand, using Lemma \ref{eq:intbubbleest}, we infer that the following estimates hold
\begin{equation}\label{eq:lp1}
||Ke^{4\sum_{i=1}^m\varphi_{x_{i, l}, \l_{i, l}}}||_{L^{p_1}(M)}\leq C\sum_{i=1}^m\l_{i, l}^{8p_1-4},
\end{equation}
and
\begin{equation}\label{eq:l1}
||Ke^{4\sum_{i=1}^m\varphi_{x_{i, l}, \l_{i, l}}}||_{L^{1}(M)}\geq C^{-1}\sum_{i=1}^m\l_{i, l}^{4},
\end{equation}
where $C$ is a large positive constant independent of $l$. Thus, combining \eqref{eq:infopl}, \eqref{eq:w22wl}, \eqref{eq:lp1} and \eqref{eq:l1}, we get
\begin{equation}\label{eq:diff0}
||w_l||_{W^{2, 2}(M)}\frac{||Ke^{4\sum_{i=1}^m\varphi_{x_{i, l}, \l_{i, l}}}||_{L^{p_1}(M)}}{||Ke^{4\sum_{i=1}^m\varphi_{x_{i, l}, \l_{i, l}}}||_{L^{1}(M)}}\leq C\sum_{i=1}^m\frac{1}{\l_{i, l}^{1-\gamma+8p_1-8}},
\end{equation}
where $C$ is a large positive constant independent of $l$. Using again \eqref{eq:defwl}, Lagrange theorem, H\"older inequality, Sobolev embedding theorem, and Moser-Trudinger inequality, we get
\begin{equation}\label{eq:diff5}
\int_M Ke^{4(v_l-\ov{(v_l)}_Q)}dV_g=\int_M Ke^{4\sum_{i=1}^m\varphi_{x_{i, l}, \l_{i, l}}}dV_g\left(1+O\left(\frac{||Ke^{4\sum_{i=1}^m\varphi_{x_{i, l}, \l_{i, l}}}||_{L^{p_1}(M)}||||w_l||_{W^{2, 2}(M)}}{||Ke^{4\sum_{i=1}^m\varphi_{x_{i, l}, \l_{i, l}}}||_{L^1(M)}}\right)\right).
\end{equation}
Now, we are going to estimate $\frac{\int_M Ke^{4(v_l-\ov{(v_l)}_Q)}\varphi dV_g}{\int_M Ke^{4(v_l-\ov{(v_l)}_Q)}dV_g}-\frac{\int_M Ke^{4\sum_{i=1}^m\varphi_{x_{i, l}, \l_{i, l}}}\varphi dV_g}{\int_M Ke^{4\sum_{i=1}^m\varphi_{x_{i, l}, \l_{i, l}}}dV_g}$. For this, we use \eqref{eq:diff0}, \eqref{eq:diff5}, H\"older inequality, Moser-Trudinger inequality,  Sobolev embedding theorem, to get
\begin{equation}\label{eq:diff11}
\begin{split}
&\frac{\int_M Ke^{4(v_l-\ov{(v_l)}_Q)}\varphi dV_g}{\int_M Ke^{4(v_l-\ov{(v_l)}_Q)}dV_g}-\frac{\int_M Ke^{4\sum_{i=1}^m\varphi_{x_{i, l}, \l_{i, l}}}\varphi dV_g}{\int_M Ke^{4\sum_{i=1}^m\varphi_{x_{i, l}, \l_{i, l}}}dV_g}=\\&O\left(\frac{||Ke^{4\sum_{i=1}^m\varphi_{x_{i, l}, \l_{i, l}}}||_{L^{p_1}(M)}||\varphi ||_{W^{2, 2}(M)}||w_l||_{W^{2, 2}(M)}}{||Ke^{4\sum_{i=1}^m\varphi_{x_{i, l}, \l_{i, l}}}||_{L^1(M)}}\right)\\&+O\left(\frac{||Ke^{4\sum_{i=1}^m\varphi_{x_{i, l}, \l_{i, l}}}||_{L^{p_1}(M)}^2||\varphi ||_{W^{2, 2}(M)}||w_l||_{W^{2, 2}(M)}}{||Ke^{4\sum_{i=1}^m\varphi_{x_{i, l}, \l_{i, l}}}||_{L^1(M)}^2}\right).
\end{split}
\end{equation}
Next, we are going to estimate $\int_MKe^{4\sum_{i=1}^m\varphi_{x_{i, l}, \l_{i, l}}}\varphi dV_g$. For this, we use again  H\"older inequality, Sobolev embedding Theorem and Moser-Trudinger inequality to infer that the following estimate holds
\begin{equation}\label{eq:diff12}
\left|\int_MKe^{4\sum_{i=1}^m\varphi_{x_{i, l}, \l_{i, l}}}\varphi dV_g\right|\leq C||Ke^{4\sum_{i=1}^m\varphi_{x_{i, l}, \l_{i, l}}}||_{L^{p_1}(M)}||\varphi ||_{W^{2, 2}(M)},
\end{equation}
where $C$ is a large positive constant independent of $l$. Now, combining \eqref{eq:estvlmo1}, \eqref{eq:derbu2}, \eqref{eq:diff11} and \eqref{eq:diff12}, we obtain
\begin{equation}\label{eq:diff13}
\begin{split}
\left|<\n^{W^{2, 2}} J(\sum_{i=1}^m\varphi_{x_{i, l}, \l_{i, l}}), \varphi>_{W^{2, 2}(M)}\right|\leq  &C\left(||w_l||_{W^{2, 2}(M)}||\varphi||_{W^{2, 2}(M)}+|1-t_l|||\varphi ||_{L^2(M)}\right)\\&+C\left(|1-t_l|\frac{||Ke^{4\sum_{i=1}^m\varphi_{x_{i, l}, \l_{i, l}}}||_{L^{p_1}(M)}||\varphi ||_{W^{2, 2}(M)}}{||\int_M Ke^{4\sum_{i=1}^m\varphi_{x_{i, l}, \l_{i, l}}}||_{L^1(M)}}\right)\\&+C\left(\frac{||Ke^{4\sum_{i=1}^m\varphi_{x_{i, l}, \l_{i, l}}}||_{L^{p_1}(M)}||\varphi ||_{W^{2, 2}(M)}||w_l||_{W^{2, 2}(M)}}{||Ke^{4\sum_{i=1}^m\varphi_{x_{i, l}, \l_{i, l}}}||_{L^1(M)}}\right)\\&+C\left(\frac{||Ke^{4\sum_{i=1}^m\varphi_{x_{i, l}, \l_{i, l}}}||_{L^{p_1}(M)}^2||\varphi ||_{W^{2, 2}(M)}||w_l||_{W^{2, 2}(M)}}{||Ke^{4\sum_{i=1}^m\varphi_{x_{i, l}, \l_{i, l}}}||_{L^1(M)}^2}\right).
\end{split}
\end{equation}
Next, combining \eqref{eq:infopl}, \eqref{eq:esttl}\eqref{eq:w22wl}, \eqref{eq:lp1}, \eqref{eq:l1}, and \eqref{eq:diff13}, we arrive to
\begin{equation}\label{eq:diff14}
||\n^{W^{2, 2}} J(\sum_{i=1}^m\varphi_{x_{i, l}, \l_{i, l}})||_{W^{2, 2}(M)}\leq C\sum_{i=1}^m\frac{1}{\l_{i, l}^{1-\gamma}},
\end{equation}
where $C$ is a large positive constant independent of $l$. Now, we are going to estimate  $||\n^{W^{2, 2}} J(v_l)||_{W^{2, 2}(M)}$.  For this, we first use the fact that $\n^{W^{2, 2}}J(v_l)=\n^{W^{2, 2}}J(v_l-\ov{(v_l)}_Q)$ and write
\begin{equation}\label{eq:diff15}
\begin{split}
<\n^{W^{2, 2}} J(v_l), \varphi>_{W^{2, 2}(M})=&<\n^{W^{2, 2}}J(v_l-\ov{(v_l)}_Q)-\n^{W^{2, 2}} J(\sum_{i=1}^m\varphi_{x_{i, l}, \l_{i, l}}), \varphi>_{W^{2, 2}(M)}\\&+<\n^{W^{2, 2}} J(\sum_{i=1}^m\varphi_{x_{i, l}, \l_{i, l}}), \varphi>_{W^{2, 2}(M)}.
\end{split}
\end{equation}
Now, we have \eqref{eq:diff15} implies that
\begin{equation}\label{eq:diff16}
\begin{split}
<\n^{W^{2, 2}} J(v_l), \varphi>_{W^{2, 2}(M})=&2<P_g w_l,  \varphi>+4\times 32m\pi^2\frac{\int_M Ke^{4(v_l-\ov{(v_l)}_Q)}\varphi dV_g}{\int_M Ke^{4(v_l-\ov{(v_l)}_Q)}dV_g}-\frac{\int_M Ke^{4\sum_{i=1}^m\varphi_{x_{i, l}, \l_{i, l}}}\varphi dV_g}{\int_M Ke^{4\sum_{i=1}^m\varphi_{x_{i, l}, \l_{i, l}}}dV_g}\\&+<\n^{W^{2, 2}} J(\sum_{i=1}^m\varphi_{x_{i, l}, \l_{i, l}}), \varphi>_{W^{2, 2}(M)}.
\end{split}
\end{equation}
On the other hand, it is easy to see that \eqref{eq:diff16} implies
\begin{equation}\label{eq:diff17}
\begin{split}
\left|<\n^{W^{2, 2}} J(v_l), \varphi>_{W^{2, 2}(M})\right|\leq &C||w_l|_{W^{2, 2}(M)}|||\varphi||_{W^{2, 2}(M)}\\&+C\left|\frac{\int_M Ke^{4(v_l-\ov{(v_l)}_Q)}\varphi dV_g}{\int_M Ke^{4(v_l-\ov{(v_l)}_Q)}dV_g}-\frac{\int_M Ke^{4\sum_{i=1}^m\varphi_{x_{i, l}, \l_{i, l}}}\varphi dV_g}{\int_M Ke^{4\sum_{i=1}^m\varphi_{x_{i, l}, \l_{i, l}}}dV_g}\right|\\&+C||\n^{W^{2, 2}} J(v_l)||_{W^{2, 2}(M)}||\varphi||_{W^{2, 2}(M)}.
\end{split}
\end{equation}
Thus, combining \eqref{eq:infopl}, \eqref{eq:w22wl}, \eqref{eq:lp1}, \eqref{eq:l1}, \eqref{eq:diff11},  \eqref{eq:diff14} and \eqref{eq:diff17}, we obtain
\begin{equation}\label{eq:diff18a}
||\n^{W^{2, 2}} J(v_l)||_{W^{2, 2}(M)}\leq C\sum_{i=1}^m\frac{1}{\l_{i, l}^{1-\gamma}},
\end{equation}
where $C$ is large positive constant independent of $l$. Thus, using \eqref{eq:dfvl} and the fact that $J$ is invariant by translations by constants, we infer that \eqref{eq:diff18a} implies
\begin{equation}\label{eq:diff18}
||\n^{W^{2, 2}} J(u_l)||_{W^{2, 2}(M)}\leq C\sum_{i=1}^m\frac{1}{\l_{i, l}^{1-\gamma}}.
\end{equation}
Hence the Proposition follows from \eqref{eq:infopl}, \eqref{eq:estvlw22}, \eqref{eq:w22wl} and \eqref{eq:diff18}.
\end{pf}
\vspace{6pt}

\noindent
Now, using Lemma \ref{eq:deformlem} and Proposition \ref{eq:escape}, we derive the following refined deformation  Lemma.
\begin{lem}\label{eq:deformlemr} (Refined deformation Lemma)\\
Assuming that $\epsilon$ and $\eta$ are small positive real numbers with $0<2\eta<\varrho$, then for $a, b\in \R$ such that $a<b$, we have that if  there is no critical values of $J$ in $[a, b]$, then there are two possibilities\\
1) Either  $$J^a\;\text{is a deformation retract of}\;\; J^b.$$
2) Or there exists a sequence $t_l\rightarrow 1$ as $l\rightarrow +\infty$ and a sequence of critical point $u_l$ of $J_{t_l}$ verifying $a\leq J(u_l)\leq b$ for all $l\in \N^*$ and $l_{\epsilon, \eta}$ a large positive integer such that  $u_l\in V(m, \epsilon, \eta)$ for all $l\geq l_{\epsilon, \eta}$.
\end{lem}
\begin{pf}
It follows directly from Lemma \ref{eq:deformlem}, Proposition \ref{eq:escape} and the fact that $J_t$ is invariant by translation by constants for every $t>0$.
\end{pf}

\section{An useful expansion of the functional $J$ at infinity}\label{eq:jexp}
In this section, we perform an expansion of the Euler-Lagrange functional $J$ in a subset of a neighborhood of potential critical points at infinity. For this, we first fix $R$ to be a very large positive number. Next, following the ideas of Bahri-Coron \cite{bc1}, and using Lemma \ref{eq:intbubbleest} and Lemma \ref{eq:outbubbleest} (see also \cite{bah} and \cite{cl2}), we have that for $\eta$ a small positive real number with $0<2\eta<\varrho$, there exists $\epsilon_0=\epsilon_0(\eta)>0$ such that
\begin{equation}\label{eq:mini}
\forall\;0<\epsilon\leq \epsilon_0,\;\;\forall u\in V(m, \epsilon, \eta), \text{the minimization problem }\min_{B_{\epsilon, \eta}}||u-\ov{u}_Q-\sum_{i=1}^m\alpha_i\varphi_{a_i, \l_i}-\sum_{r=^1}^{\bar m}\beta_r(v_r-\ov{(v_r)}_Q)||
\end{equation}
has a unique solution, up to permutations, where $B_{\epsilon, \eta}$ is defined as follows
\begin{equation}
\begin{split}
{B_{\epsilon, \eta}:=\{(\bar\alpha, A, \bar \l, \bar \beta)\in \R^m\times M^m\times (0, +\infty)^m\times \R^{\bar m}}:|\alpha_i-1|\leq \epsilon, \l_i\geq \frac{1}{\epsilon}, i=1, \cdots, m, \\d_g(a_i, a_j)\geq 4\ov C\eta, i\neq j, |\beta_r|\leq R, r=1, \cdots, \bar m\}.
\end{split}
\end{equation}
 Moreover, using the solution of \eqref{eq:mini}, we have that every $u\in V(m, \epsilon, \eta)$ can be written as
\begin{equation}\label{eq:para}
u-\ov{u}_Q=\sum_{i=1}^m\alpha_i\varphi_{a_i, \l_i}+\sum_{r=1}^{\bar m}\beta_r(v_r-\ov{(v_r}_Q)+w,
\end{equation}
where $w$ verifies the following orthogonality conditions
\begin{equation}\label{eq:ortho}
\begin{split}
<Q_g, w>=<\varphi_{a_i, \l_i}, w>_{P}=<\frac{\partial\varphi_{a_i, \l_i}}{\partial \l_i}, w>_{P}=<\frac{\partial\varphi_{a_i, \l_i}}{\partial a_i}, w>_{P}=<v_r, w>=0, i=1, \cdots, m,\\r=1, \cdots, \bar m
\end{split}
\end{equation}
and the estimate
\begin{equation}\label{eq:estwmin}
||w||=O\left(\sum_{i=1}^m\frac{1}{\l_i^{1-\gamma}}\right),
\end{equation}
where here $O\left(1\right):=O_{\bar \alpha, A, \bar \l, \bar\beta , w, \epsilon}\left(1\right)$,  and for the meaning of , see section \ref{eq:notpre}. Furthermore, the concentration points $a_i$,  the masses $\alpha_i$, the concentrating parameters $\l_i$ and the negativity parameter $\beta_r$ in \eqref{eq:para} verify also
\begin{equation}\label{eq:afpara}
\begin{split}
d_g(a_i, a_j)\geq 4\ov C\eta,\;i\neq j=1, \cdots, m, \frac{1}{\L}\leq\frac{\l_i}{\l_j}\leq\L\;\;i, j=1, \cdots, m, \;\;\l_i\geq\frac{1}{\epsilon},\;\;\text{and}\\\;\;\;\sum_{r=1}^{\bar m}|\beta_r|+\sum_{i=1}^m|\alpha_i-1|\sqrt{\log \l_i}=O\left(\sum_{i=1}^m\frac{1}{\l_i^{1-\gamma}}\right)
\end{split}
\end{equation}
with still $O\left(1\right)$ as in \eqref{eq:estwmin}.
\vspace{6pt}

\noindent
Now as pointed out at the beginning of the section, we are going to give the following Lemma which gives a first expansion of $J$ in a useful subset of $V(m, \epsilon, \eta)$, precisely for elements in $V(m, \epsilon, \eta)$ having $w=0$  in the representation \eqref{eq:para}. Indeed, we have:
\begin{lem}\label{eq:energyest} (Energy estimate in $V(m, \epsilon, \eta)\cap\{w=0\}$)\\
Assuming that $\eta$ is a small positive real number with $0<2\eta<\varrho$ where $\varrho$ is as in \eqref{eq:cutoff}, and $0<\epsilon\leq \epsilon_0$ where $\epsilon_0$ is as in \eqref{eq:mini}, then  for $a_i\in M$ concentration points,  $\alpha_i$ masses, $\l_i$ concentration parameters ($i=1,\cdots,m$), and $\beta_r$ negativity parameters ($r=1, \cdots, \bar m$) satisfying \eqref{eq:afpara}, we have
\begin{equation*}
\begin{split}
J(\sum_{i=1}^m\alpha_i\varphi_{a_i, \l_i}+\sum_{r=^1}^{\bar m}\beta_r(v_r-\ov{(v_r)}_Q))&=-\frac{40}{3}m\pi^2-8m\pi^2\log(\frac{m\pi^2}{6})-8\pi^2\mathcal{F}_K(a_1, \dots, a_m)\\&+\sum_{i=1}^m(\alpha_i-1)^2\left(32\pi^2\log\l_i+16\pi^2H(a_i, a_i)+\frac{472\pi^2}{9}\right)+\sum_{r=1}^{\bar m}\mu_r\beta_r^2\\&+16\pi^2\sum_{i=1}^m(\alpha_i-1)\left[\sum_{r=1}^{\bar m}2\beta_r(v_r-\ov{(v_r)}_Q)(a_i)-\sum_{ j=1, j\neq i}^m (\alpha_j-1)G(a_i, a_j)\right]\\&-2\pi^2\sum_{i=1}^m\frac{1}{\l_i^2}\left(\frac{\D_{g_{a_i}}\mathcal{F}^{A}_i(a_i)}{\mathcal{F}^{A}_i(a_i)}-\frac{2}{3}R_g(a_i)\right)\\&+2\pi^2\sum_{i=1}^m\frac{\tilde \tau_i}{\l_i^2}\left(\frac{\D_{g_{a_i}}\mathcal{F}^{A}_i(a_i)}{\mathcal{F}^{A}_i(a_i)}-\frac{2}{3}R_g(a_i)\right)\\&+8\pi^2\sum_{i=1}^m\log(1-\tilde\tau_i)+O\left(\sum_{i=1}^m|\alpha_i-1|^3+\sum_{r=1}^{\bar m}|\beta_r|^3+\sum_{i=1}^m\frac{1}{\l^3_i}\right),
\end{split}
\end{equation*}
where \;$O\left(1\right)$ means here \;$O_{\bar\alpha, A, \bar \l, \bar \beta, \epsilon}\left(1\right)$ \;with \;$\bar\alpha=(\alpha_1, \cdots, \alpha_m)$, $A:=(a_1, \cdots, a_m)$ $\bar \l:=(\l_1, \cdots, \l_m)$, $\bar \beta:=(\beta_1, \cdots, \beta_{\bar m})$ and for $i=1, \cdots, m$, $$\tilde \tau_i:=1-\frac{m\gamma_i}{\Gamma},\;\;\;\; \Gamma:=\sum_{i=1}^m\gamma_i, \;\;\;\gamma_i:=\frac{\pi^2}{2(2\alpha_i-1)(2\alpha_i-1)}\l_i^{8\alpha_i-4}\mathcal{F}^{A}_i(a_i)\mathcal{G}_i(a_i),$$ with $$\mathcal{G}_i(a_i):=e^{4((\alpha_i-1)H(a_i, a_i)+\sum_{j=1, j\neq i}^m(\alpha_j-1)G(a_j, a_i))}e^{\sum_{j=1, j\neq i}^m\frac{\alpha_j}{\l_j^2}\D_{g_{a_j}}G(a_j, a_i)}e^{\frac{\alpha_i}{\l_i^2}\D_{g_{a_i}}H(a_i, a_i)}e^{\sum_{r=1}^{\bar m}\beta_rv_r(a_i)},$$
and for the meaning of $O_{\bar\alpha, A, \bar \l, \bar \beta, \epsilon}\left(1\right)$, see section \ref{eq:notpre}.
\end{lem}
\begin{pf}
By definition of $J$ (see \eqref{eq:defj1}) and the fact that it is invariant by translations by constants combined with the linearity of $P_g$ and of scalar products, and the normalization $\int_MQ_g(x)\varphi_{a_j, \l_j}(x)dV_g(x)=0$ (see \eqref{eq:projbubble} with $a$ replaced by $a_j$ and $\l$ replaced by $\l_j$) for $j=1,\cdots, m$, we have
\begin{equation}\label{eq:5imp1}
J(\sum_{i=1}^m\alpha_i\varphi_{a_i, \l_i})=E_1+E_0-8m\pi^2 E_2,
\end{equation}
where
\begin{equation}\label{eq:e1}
E_1:=\sum_{i=1}^m\alpha_i^2<P_g\varphi_{a_i, \l_i}),\varphi_{a_i,\l_i})>+\sum_{i=1}^m\sum_{j=1, j\neq i}^m\alpha_i\alpha_j<P_g\varphi_{a_i, \l_i}, \varphi_{a_j,\l_j}>,
\end{equation}
\begin{equation}\label{eq:e0}
E_0:=\sum_{j=1}^{\bar m}\mu_r\beta_r^2+2\sum_{i=1}^m\sum_{r=1}^{\bar m}\alpha_i\beta_r<P_g \varphi_{a_i, \l_i}, v_r>+4\int_M Q_g(x)(\sum_{r=1}^{\bar m}\beta_rv_r)dV_g(x)
\end{equation}
and
\begin{equation}\label{eq:e2e3}
E_3:=\int_M Ke^{4(\sum_{j=1}^m\alpha_j\varphi_{a_j,\l_j}(x)+4\sum_{r=1}^{\bar m}\beta_r v_r(x))}dV_g(x),\;\;\;\;E_2=\log E_3.
\end{equation}
Now let us estimate $E_1$, $E_0$ and $E_2$. We start with $E_1$. For this, we use \eqref{eq:e1}, the first estimate of Lemma \ref{eq:selfintest},  the first estimate of Lemma \ref{eq:interactest}, the relations $\alpha_i^2=1+2(\alpha_i-1)+(\alpha_i-1)^2$ ($i=1, \cdots, m$),  $\alpha_i\alpha_j=1+(\alpha_i-1)+(\alpha_j-1)+(\alpha_i-1)(\alpha_j-1)$ ($i, j=1, \cdots, m$) and $|\alpha_i-1|=O\left(\frac{1}{\l_i^{1-\gamma}\sqrt{\log\l_i}}\right)$ ($i=1, \cdots, m$) to obtain
\begin{equation}\label{eq:5imp10}
\begin{split}
E_1&=32\pi^2 \sum_{i=1}^m\log \l_i+64\pi^2\sum_{i=1}^m(\alpha_i-1)\log\l_i+32\pi^2\sum_{i=1}^m(\alpha_i-1)^2\log \l_i-\frac{40\pi^2}{3}m-\frac{80\pi^2}{3}\sum_{i=1}^m(\alpha_i-1)\\&+16\pi^2\left(\sum_{i=1}^m H(a_i, a_i)+\sum_{i=1}^m\sum_{j=1, j\neq i}^mG_(a_j, a_i)\right)-\frac{40\pi^2}{3}\sum_{i=1}^m(\alpha_i-1)^2+16\pi^2\sum_{i=1}^m(\alpha_i-1)^2H(a_i, a_i)\\&+32\pi^2\left(\sum_{i=1}^m(\alpha_i-1) H(a_i, a_i)+\sum_{i=1}^m\sum_{j=1, j\neq i}^m(\alpha_j-1) G(a_j, a_i)\right)+16\pi^2\sum_{i=1}^m\sum_{j=1, j\neq i}^m(\alpha_i-1)(\alpha_j-1)G(a_j, a_i)\\+&8\pi^2\left(\sum_{i=1}^m \frac{1}{\l_i^2} \D_{g_{a_i}} H(a_i, a_i)+\sum_{i=1}^m\sum_{j=1, j\neq i}^m\frac{1}{\l_j^2} \D_{g_{a_j}} G(a_j, a_i))\right)\\&+O\left(\sum_{i=1}^m\frac{1}{\l^3_i}+\sum_{i=1}^m\frac{|\alpha_i-1|}{\l_i^2}\right).
\end{split}
\end{equation}
Next, let us estimate $E_2$. For this, we start by estimating $E_3$. In order to estimate $E_3$, we first use \eqref{eq:e2e3}, the fact that $\sum_{r=1}^{\bar m}|\beta_r|+\sum_{i=1}^m|\alpha_i-1|\sqrt{\log \l_i}=O\left(\sum_{i=1}^m\frac{1}{\l_i^{1-\gamma}}\right)$,  the fact that $d_g(a_i, a_j)\geq 4\ov C\eta$ for $i, j=1, \cdots, m$ and $i\neq j$, the properties of the metrics $g_{a_i}$ (see \eqref{eq:detga}-\eqref{eq:proua}) $i=1, \cdots, m$, the first estimate of Lemma \ref{eq:outbubbleest}, \eqref{eq:decompG4}, and \eqref{eq:regH4}, to get
\begin{equation} \label{eq:5imp13}
E_3=\sum_{i=1}^mE_3^{i}+O\left(1\right),
\end{equation}
where
\begin{equation}\label{eq:5imp12}
E_3^{i}:=\int_{B^{a_i}_{a_i}(\eta)} K(x)e^{(4\alpha_i\varphi_{a_i, \l_i}(x)+4(\sum_{j=1, j\neq i}^m\alpha_j\varphi_{a_j, \l_j}(x)+\sum_{r=1}^{\bar m}\beta_r v_r(x)))}dV_g(x),\;\;\;i=1,\cdots, m,
\end{equation}
Now, let us estimate each $E_3^{i}$ ($i=1, \cdots, m$). In order to do that, we use again the fact that $d_g(a_i, a_j)\geq 4\ov C\eta$ for $i\neq j=1, \cdots, m$, the properties of the metrics $g_{a_j}$ (see \eqref{eq:detga}-\eqref{eq:proua} with $a$ replaced by $a_j$), $j=1, \cdots, m$, the first estimate of Lemma \ref{eq:intbubbleest}, the first estimate of Lemma \ref{eq:outbubbleest}, the explicit expression of $\hat \d_{a_i, \l_i}$ (see \eqref{eq:hatdelta} with $a$ replaced by $a_j$ and $\l$  replaced by $\l_j$), the definition of $\chi_{\varrho}$ (see \eqref{eq:cutoff}), the fact that $0<2\eta<\varrho$, and \eqref{eq:5imp12}, to obtain
\begin{equation}\label{eq:5imp14ab}
E_3^{i}=\tilde{E}_3^{i}\left(1+O\left(\sum_{j=1}^m\frac{1}{\l^3_j}\right)\right),
\end{equation}
where
\begin{equation}\label{eq:5imp17o}
\tilde {E}_3^{i}=\int_{B_{a_i}^{a_i}(\eta)}\frac{\l_i^{8\alpha_i}}{(1+\l_i^2d_{g_{a_i}}^2(a_i, x))^{4\alpha_i}}\mathcal{F}_i^{A}(x)\tilde {\mathcal{G}_i}(x)e^{\frac{\alpha_i}{\l_i^2}\D_{g_{a_i}}H(a_i, x)}dV_{g_{a_i}}(x),
\end{equation}
with
\begin{equation}\label{eq:5imp16}
\mathcal{F}_i^{A}(x)=e^{4(H(a_i, x)+\sum_{j=1, j\neq i}^mG(a_j, x))+\frac{1}{4}\log(K(x))}, \;\;A=(a_1,\cdots,a_m),
\end{equation}
\begin{equation}\label{eq:5imp15}
\tilde{\mathcal{G}_i}(x)=\bar{\mathcal{G}_i}(x)e^{-4u_{a_i}(x)},
\end{equation}
and
\begin{equation}\label{eq:5imp14}
\bar{\mathcal{G}_i}(x)=e^{4((\alpha_i-1)H(a_i, x)+\sum_{j=1, j\neq j}^m(\alpha_j-1)G(a_j, x))}\times e^{\sum_{j=1, j\neq i}^m\frac{\alpha_j}{\l_j^2}\D_{g_{a_j}}G(a_j, x)}e^{4\sum_{r=1}^{\bar m}\beta_r v_r(x)}.
\end{equation}
Next, we define
\begin{equation}\label{eq:5imp14new}
\mathcal{G}_i(x):=\bar{\mathcal{G}_i}(x)e^{\frac{\alpha_i}{\l_i^2}\D_{g_{a_i}}H(a_i, x)}.
\end{equation}
Furthermore, setting
\begin{equation}\label{eq:5imp17}
\tilde{E}^{i, 1}_3:=\int_{B_{a_i}^{a_i}(\eta)}\frac{\l_i^{8\alpha_i}}{(1+\l_i^2d_{g_{a_i}}^2(a_i, x))^{4\alpha_i}}\mathcal{F}_i^{A}(x)\tilde {\mathcal{G}_i}(x)dV_{g_{a_i}}(x),
\end{equation}
and
\begin{equation}\label{eq:5imp17bu2}
\tilde{E}_3^{i, 2}:=-\int_{B_{a_i}^{a_i}(\eta)}\frac{\l_i^{8\alpha_i}}{(1+\l_i^2d_{g_{a_i}}^2(a_i, x))^{4\alpha_i}}\mathcal{F}_i^{A}(x)\tilde {\mathcal{G}_i}(x)\D_{g_{a_i}}H(a_i, x)dV_{g_{a_i}}(x),
\end{equation}
we have that \eqref{eq:5imp17o} gives
\begin{equation}\label{eq:7imp17bu3}
\tilde{E}^{i}_3=\tilde{E}^{i, 1}_3\left(1+O\left(\frac{1}{\l^3_i}\right)\right)-\frac{\alpha_i}{\l_i^2}\tilde{E}^{i, 2}_3.
\end{equation}
Now, let us estimate first \;$\tilde{E}_3^{i, 1}$\; and after \;$\tilde{E}_3^{i, 2}$. For this, we use Taylor expansion at $a_i$, the change of  perform the change of variable \;$exp_{a_i}^{a_i}(\frac{\zeta}{\l_i})=x$ in  \eqref{eq:5imp17}, the properties of $g_{a_i}$ (see \eqref{eq:detga} with $a$ replaced by $a_i$), \eqref{eq:5imp15}, $|\alpha_i-1|=O\left(\frac{1}{\l_{i}^{1-\gamma}\sqrt{\log \l_i}}\right)$, and direct calculations to have
\begin{equation}\label{eq:5imp26}
\begin{split}
\tilde {E}_3^{i, 1}&=\l_i^{8\alpha_i-4}\mathcal{F}_i^{A}(a_i)\bar{\mathcal{G}_i}(a_i)\left(\frac{\pi^2}{2(4\alpha_i-1)(2\alpha_i-1)}+O\left(\frac{1}{\l^4_i}\right)\right)\\&+\frac{1}{8}\l_i^{8\alpha_i-6}\D_{g_{a_i}}(\mathcal{F}_i^{A}\tilde {\mathcal{G}_i})(a_i)\left(\frac{\pi^2}{(4\alpha_i-1)(2\alpha_i-1)(4\alpha_i-3)}+O\left(\frac{1}{\l_i}\right)\right)\\&+O\left(\l_i\right).
\end{split}
\end{equation}
Now, using \eqref{eq:lapemua}, \eqref{eq:decompG4}, \eqref{eq:regH4}, the fact that $\sum_{r=1}^{\bar m}|\beta_r|+\sum_{i=1}^m|\alpha_i-1|\sqrt{\log\l_i}=O\left(\sum_{i=1}^m\frac{1}{\l_i^{1-\gamma}}\right)$, \eqref{eq:5imp15}, and \eqref{eq:5imp14}, it is easy to see that
\begin{equation}\label{eq:5imp27a}
\begin{split}
\bar{\mathcal{G}_i}(a_i)=\tilde{\mathcal{G}_i}(a_i)=1+O\left(\sum_{j=1}^m|\alpha_j-1|+\sum_{r=1}^{\bar m}|\beta_r|+\sum_{j=1}^m\frac{1}{\l_j^2}\right),
\end{split}
\end{equation}
and
\begin{equation}\label{eq:5imp28}
\D_{g_{a_i}}(\mathcal{F}_i^{A}\tilde {\mathcal{G}_i})(a_i)=\D_{g_{a_i}}\mathcal{F}_i^{A}(a_i)-\frac{2}{3}R_g(a_i)\mathcal{F}_i^{A}(a_i)+O\left(\sum_{j=1}^m|\alpha_j-1|+\sum_{r=1}^{\bar m}|\beta_r|+\sum_{j=1}^m\frac{1}{\l_j^2}\right).
\end{equation}
Combining \eqref{eq:5imp26} and \eqref{eq:5imp28}, we arrive to
\begin{equation}\label{eq:5imp29}
\begin{split}
\tilde {E}_3^{i, 1}=&\l_i^{8\alpha_i-4}\mathcal{F}_i^{A}(a_i)\bar{\mathcal{G}_i}(a_i)\left(\frac{\pi^2}{2(4\alpha_i-1)(2\alpha_i-1)}+O\left(\frac{1}{\l^4_i}\right)\right)\\&+\frac{1}{8}\l_i^{8\alpha_i-6}\left(\D_{g_{a_i}}\mathcal{F}_i^{A}(a_i)-\frac{2}{3}R_g(a_i)\mathcal{F}_i^{A}(a_i)+O\left(\sum_{j=1}^m|\alpha_j-1|+\sum_{r=1}^{\bar m}|\beta_r|+\sum_{j=1}^m\frac{1}{\l_j^2}\right)\right)\times\\&\left(\frac{\pi^2}{(4\alpha_i-1)(2\alpha_i-1)(4\alpha_i-3)}+O\left(\frac{1}{\l_i}\right)\right)\\&+O\left(\l_i\right).
\end{split}
\end{equation}
Next, setting
\begin{equation}\label{eq:5imp30}
 \gamma_i^1=\frac{\pi^2}{2(4\alpha_i-1)(2\alpha_i-1)}\l_i^{8\alpha_i-4}\mathcal{F}_i^{A}(a_i)\bar{\mathcal{G}_i}(a_i),
\end{equation}
and using\eqref{eq:decompG4}, \eqref{eq:regH4}, \eqref{eq:5imp14}, \eqref{eq:5imp16} and the fact that $\sum_{r=1}^{\bar m}|\beta_r|+\sum_{j=1}^m|\alpha_j-1|\sqrt{\log\l_j}=O\left(\sum_{j=1}^m\frac{1}{\l^{1-\gamma}_j}\right)$, it is easy  to see that
\begin{equation}\label{eq:5imp31}
C_0^{-1}\l_i^4\leq  \gamma_i^1\leq C_0\l_i^4,
\end{equation}
for some large positive constant $C_0$ independent of $\epsilon$, $\l_j$, $a_j$, $\alpha_j$, $\beta_r$, $j=1, \cdots, m$ and $r=1, \cdots, \bar m$.. Thus, setting
\begin{equation}\label{eq:tildegammai}
\gamma_i:=\tilde\gamma_i^1e^{\frac{\alpha_i}{\l_i^2}\D_{g_{a_i}}H(a_i, a_i)}=\frac{\pi^2}{2(4\alpha_i-1)(2\alpha_i-1)}\l_i^{8\alpha_i-4}\mathcal{F}_i^{A}(a_i)\mathcal{G}_i(a_i),
\end{equation}
and using  the fact that $|\alpha_i-1|=O\left(\frac{1}{\l^{1-\gamma}_i\sqrt{\log\l_i}}\right)$, \eqref{eq:regH4}, \eqref{eq:5imp31}, we obtain
\begin{equation}\label{eq:5imp31a}
C_1^{-1}\l_i^4\leq \gamma_i\leq C_1\l_i^4,\;\;\;\text{and}\;\,\;\gamma_i=\gamma_i^1\left(1+O\left(\frac{1}{\l^2_i}\right)\right),
\end{equation}
for some large positive constant $C_1$ independent of independent of $\epsilon$, $\l_k$, $a_k$ and $\alpha_k$, $\beta_j$, $k=1, \cdots, m$ and $j=1, \cdots, \bar m$. On the other hand, using again the fact that $|\alpha_i-1|=O\left(\frac{1}{\l^{1-\gamma}_i\sqrt{\log\l_i}}\right)$ combined \eqref{eq:decompG4}, \eqref{eq:regH4}, \eqref{eq:partiallimit}, \eqref{eq:5imp27a} \eqref{eq:5imp30}, \eqref{eq:5imp31}, we get
\begin{equation}\label{eq:5imp34}
\begin{split}
\tilde {E}_3^{i, 1}=\gamma_i^1\left(1+\frac{1}{4\l_i^2}\left(\frac{\D_{g_{a_i}}\mathcal{F}_i^{A}(a_i)}{\mathcal{F}^{A}_i(a_i)}-\frac{2}{3}R_g(a_i)\right)+O\left(\sum_{j=1}^m\frac{|\alpha_j-1|}{\l_j^2}+(\sum_{r=1}^{\bar m}|\beta_r|)(\sum_{j=1}^m\frac{1}{\l_j^2})+\sum_{j=1}^m\frac{1}{\l_j^3}\right)\right).
\end{split}
\end{equation}
Now, we are going to estimate $\tilde{E}_3^{i, 2}$. For this, use the same arguments as above and get
\begin{equation}\label{eq:5imp23c}
\begin{split}
\tilde {E}_3^{i, 2}&=\gamma_i^1\left(-\D_{g_{a_i}}H(a_i, a_i)+O\left(\frac{1}{\l_i}\right)\right).
\end{split}
\end{equation}
On the other hand, using \eqref{eq:regH4} and the fact that  $|\alpha_i-1|=O\left(\frac{1}{\l_i^{1-\gamma}\sqrt{\log \l_i}}\right)$, we infer that
\begin{equation}\label{eq:deficit}
e^{\pm\frac{\alpha_i}{\l_i^2} \D_{g_{a_i}} H(a_i, a_i)}=1\pm\frac{\alpha_i}{\l_i^2} \D_{g_{a_i}} H(a_i, a_i)+O\left(\frac{1}{\l^3_i}\right).
\end{equation}
So combining \eqref{eq:decompG4}, \eqref{eq:regH4}, \eqref{eq:5imp13}, \eqref{eq:5imp14ab},  \eqref{eq:7imp17bu3}, \eqref{eq:tildegammai}, \eqref{eq:5imp34}, \eqref{eq:5imp23c}, \eqref{eq:deficit}, the fact that $|\alpha_i-1|=O\left(\frac{1}{\l^{1-\gamma}_i\sqrt{\log \l_i}}\right)$, and  factorizing by $e^{\frac{\alpha_i}{\l_i^2} \D_{g_{a_i}} H(a_i, a_i)}$, we obtain
\begin{equation}\label{eq:5imp36}
\begin{split}
E_3=\sum_{i=1}^m\gamma_i\left(1+\frac{1}{4\l_i^2}\left(\frac{\D_{g_{a_i}}\mathcal{F}_i^{A}(a_i)}{\mathcal{F}^{A}_i(a_i)}-\frac{2}{3}R_g(a_i)\right)+O\left(\sum_{i=1}^m\frac{|\alpha_i-1|}{\l_i^2}+(\sum_{r=1}^{\bar m}|\beta_r|)(\sum_{i=1}^m\frac{1}{\l_i^2})+\sum_{i=1}^m\frac{1}{\l^3_i}\right)\right)\\+O\left(1\right).
\end{split}
\end{equation}
Thus, setting
\begin{equation}\label{eq:5im38}
\Gamma=\sum_{i=1}^m \gamma_i,
\end{equation}
and
\begin{equation}\label{eq:5imp42}
\tilde\tau_i=1-\frac{m\gamma_i}{ \Gamma},\;\;\;\;i=1, \cdots, m,
\end{equation}
we have that \eqref{eq:5imp36} becomes
\begin{equation}\label{eq:5imp43}
\begin{split}
E_3=& \Gamma\left(1+\frac{1}{4m}\sum_{i=1}^m\frac{1}{\l_i^2}\left(\frac{\D_{g_{a_i}}\mathcal{F}_i^{A}(a_i)}{\mathcal{F}^{A}_i(a_i)}-\frac{2}{3}R_g(a_i)\right)\right)\\&-\Gamma\left(\frac{1}{4m}\sum_{i=1}^m\tilde\tau_i\frac{1}{\l_i^2}\left(\frac{\D_{g_{a_i}}\mathcal{F}_i^{A}(a_i)}{\mathcal{F}^{A}_i(a_i)}-\frac{2}{3}R_g(a_i)\right)+O\left(\sum_{i=1}^m\frac{|\alpha_i-1|}{\l_i^2}+(\sum_{r=1}^{\bar m}|\beta_r|)(\sum_{i=1}^m\frac{1}{\l_i^2})+\sum_{i=1}^m\frac{1}{\l^3_i}\right)\right).
\end{split}
\end{equation}
Now, we are ready to estimate $E_2$. In order to do that, we use \eqref{eq:decompG4}, \eqref{eq:regH4}, \eqref{eq:e2e3}, \eqref{eq:tildegammai}, \eqref{eq:5imp31a}, \eqref{eq:5im38}, \eqref{eq:5imp42}, and \eqref{eq:5imp43} to infer that
\begin{equation}\label{eq:5imp49}
\begin{split}
E_2=&\log(\frac{m\pi^2}{6})-\frac{1}{m}\sum_{i=1}^m\frac{10}{3}(\alpha_i-1)-\frac{74}{9}(\alpha_i-1)^2+\frac{4}{m}\sum_{i=1}^m\log\l_i+\frac{8}{m}\sum_{i=1}^m(\alpha_i-1)\log\l_i\\&+\frac{1}{m}\sum_{i=1}^m\left(\log(\mathcal{F}^{A}_i(a_i))+\log(\mathcal{G}_i(a_i))\right)+\frac{1}{4m}\sum_{i=1}^m\frac{1}{\l_i^2}\left(\frac{\D_{g_{a_i}}\mathcal{F}_i^{A}(a_i)}{\mathcal{F}^{A}_i(a_i)}-\frac{2}{3}R_g(a_i)\right)\\&-\frac{1}{4m}\sum_{i=1}^m\tilde\tau_i\frac{1}{\l_i^2}\left(\frac{\D_{g_{a_i}}\mathcal{F}_i^{A}(a_i)}{\mathcal{F}^{A}_i(a_i)}-\frac{2}{3}R_g(a_i)\right) -\frac{1}{m}\sum_{i=1}^m\log(1-\tilde\tau_i)\\&+O\left(\sum_{i=1}^m\frac{|\alpha_i-1|}{\l_i^2}+(\sum_{r=1}^{\bar m}|\beta_r|)(\sum_{i=1}^m\frac{1}{\l_i^2})+\sum_{i=1}^m|\alpha_i-1|^3+\sum_{i=1}^m\frac{1}{\l^3_i}\right).
\end{split}
\end{equation}
Next, we are going to estimate $E_0$. For this, we use \eqref{eq:qva} and \eqref{eq:projbubble} combined with \eqref{eq:e0} to obtain
\begin{equation}\label{eq:este01}
E_0=\sum_{r=1}^{\bar m}\mu_r\beta_r^2-32\pi^2\sum_{i=1}^m\sum_{r=1}^{\bar m}(\alpha_i-1)\beta_r\ov{(v_r)}_Q+32\pi^2\sum_{i=1}^m\sum_{r=1}^{\bar m}\frac{\alpha_i\beta_r}{\int_{M}e^{4\hat \d_{a_i, \l_i}}dV_{g{a_i}}}\int_{M}e^{4\hat \d_{a_i, \l_i}}v_rdV_{g_{a_i}}.
\end{equation}
On the other hand, using the explicit expression of $\hat \d_{a_i, \l_i}$, the smoothness of $v_r$, normal coordinate at $a_i$ with respect to $g_{a_i}$, and  Taylor expansion at $a_i$, we get
\begin{equation}\label{eq:massneg}
\frac{1}{\int_{M}e^{4\hat \d_{a_i, \l_i}}dv_{g{a_i}}}\int_{M}e^{4\hat \d_{a_i, \l_i}}v_rdV_{g_{a_i}}= v_j(a_i)+O(\frac{1}{\l_i^2})
\end{equation}
Thus, using \eqref{eq:este01} and \eqref{eq:massneg},  and recalling for $i=1, \cdots, m$ that $|\alpha_i-1|=O\left(\frac{1}{\l_i^{1-\gamma}\sqrt{\log \l_i}}\right)$, we get
\begin{equation}\label{eq:este02}
E_0=\sum_{r=1}^{\bar m}\mu_r\beta_jr^2-32\pi^2\sum_{i=1}^m\sum_{r=1}^{\bar m}(\alpha_i-1)\beta_r\ov{(v_r)}_Q+32\pi^2\sum_{i=1}^m\sum_{r=1}^{\bar m}\alpha_i\beta_rv_r(a_i)+O\left(\sum_{r=1}^{\bar m}|\beta_r|)(\sum_{i=1}^m\frac{1}{\l_i^2})\right).
\end{equation}
Hence, combining \eqref{eq:5imp1}, \eqref{eq:5imp10}, \eqref{eq:5imp14}, \eqref{eq:5imp14new}, \eqref{eq:5imp16},  \eqref{eq:5imp49}, \eqref{eq:este02} and using Young's inequality, we arrive to
\begin{equation*}
\begin{split}
&J(\sum_{i=1}^m\alpha_i\varphi_{a_i, \l_i}+\sum_{r=1}^{\bar m}\beta_r(v_r-\ov{(v_r)}_Q))=-\frac{40\pi^2}{3}m-8m\pi^2\log(\frac{m\pi^2}{6})-8\pi^2\mathcal{F}_K(a_1, \dots, a_m)\\&+\sum_{i=1}^m(\alpha_i-1)^2\left(32\pi^2\log\l_i+\frac{472\pi^2}{9}+16\pi^2 H(a_i, a_i)\right)+\sum_{r=1}^m\mu_r\beta_r^2\\&+16\pi^2\sum_{i=1}^m(\alpha_i-1)\left[\sum_{r=1}^{\bar m}2\beta_r(v_j-\ov{(v_r)}_Q)(a_i)-\sum_{j=1, j\neq i}^m (\alpha_j-1)G(a_j, a_i)\right]\\&-2\pi^2\sum_{i=1}^m\frac{1}{\l_i^2}\left(\frac{\D_{g_{a_i}}\mathcal{F}^{A}_i(a_i)}{\mathcal{F}^{A}_i(a_i)}-\frac{2}{3}R_g(a_i)\right)+2\pi^2\sum_{i=1}^m\tilde\tau_i\frac{1}{\l_i^2}\left(\frac{\D_{g_{a_i}}\mathcal{F}^{A}_i(a_i)}{\mathcal{F}^{A}_i(a_i)}-\frac{2}{3}R_g(a_i)\right)\\&+ 8\pi^2\sum_{i=1}^m\log(1-\tilde\tau_i)+O\left(\sum_{i=1}^m|\alpha_i-1|^3+\sum_{r=1}^{\bar m}|\beta_r|^3+\sum_{i=1}^m\frac{1}{\l^3_i}\right),
\end{split}
\end{equation*}
thereby ending the proof of the Lemma.
\end{pf}

\section{An useful expansion of the gradient of $J$ at infinity}\label{eq:gjesp}
In this section, we perform an expansion of the $L^2$-gradient of $J$ on the same set as in section \ref{eq:jexp} and in the direction of the masses parameters $\alpha_i$, the concentrations points $a_i$, concentrations parameters $\l_i$ and the negativity parameters $\beta_r$, $i=1, \cdots, m$, $r=1, \cdots, \bar m$. We start with the gradient estimate with respect to the $\l_j$'s.  Precisely, we have:
\begin{lem}\label{eq:gradientlambdaest} (Gradient estimate with respect to $\bar\l$)\\
Assuming that $\eta$ is a small positive real number with $0<2\eta<\varrho$ where $\varrho$ is as in \eqref{eq:cutoff}, and $\epsilon\leq \epsilon_0$ where $\epsilon_0$ is as in \eqref{eq:mini}, then  for $a_i\in M$ concentration points,  $\alpha_i$ masses , $\l_i$ concentration parameters ($i=1,\cdots,m$) and $\beta_r$ negativity parameters ($r=1, \cdots, \bar m$) satisfying \eqref{eq:afpara}, we have that for every $r=1,\cdots, m$, there holds
\begin{equation*}
\begin{split}
&<\n J(\sum_{i=1}^m\alpha_i\varphi_{a_i, \l_i}+\sum_{r=1}^{\bar m}\beta_r(v_r-\ov{(v_r)}_Q), \l_j\frac{\partial \varphi_{a_j, \l_j}}{\partial \l_j}>= 32\pi^2\alpha_j\tau_j-\frac{4\pi^2}{\l_j^2}\left(\frac{\D_{g_{a_j}} \mathcal{F}^{A}_j(a_j)}{\mathcal{F}^{A}_j(a_j)}-\frac{2}{3}R_g(a_j)\right)\\&-\frac{16\pi^2}{\l_j^2}\tau_j \D_{g_{a_j}}H(a_j, a_j)-\frac{16\pi^2}{\l_j^2}\sum_{i=1, i\neq j}^m\tau_i \D_{g_{a_j}}G(a_j, a_i)+ \frac{4\pi^2}{\l_j^2}\tau_j\left(\frac{\D_{g_{a_j}} \mathcal{F}^{A}_j(a_j)}{\mathcal{F}^{A}_j(a_j)}-\frac{2}{3}R_g(a_j)\right)\\&+O\left(\sum_{i=1}^m|\alpha_i-1|^2+\sum_{r=1}^{\bar m}|\beta_r|^3+\sum_{i=1}^m\frac{1}{\l_i^3}\right),
\end{split}
\end{equation*}
where $A:=(a_1, \cdots, a_m)$, $O\left(1\right)$ is as in Lemma \ref{eq:energyest}, and for $i=1, \cdots, m$, $$\tau_i:=1-\frac{m\gamma_i}{D}, \;\;\;\;\;D:=\int_M K(x)e^{4(\sum_{i=1}^m\alpha_i\varphi_{a_i, \l_i}(x)+\sum_{r=1}^{\bar m}\beta_r v_r(x))}dV_g(x),$$ with $\gamma_i$ as in Lemma \ref{eq:energyest}.
\end{lem}
\begin{pf}
First of all, using the definition of $J$ (see \eqref{eq:defj1}), the fact that $J$ is invariant by translations by constants, and the normalization $\int_MQ_g(x)\frac{\partial \varphi_{a_j, \l_j}(x)}{\partial\l_j}dV_g(x)=0$ (see \eqref{eq:normlambda} with $a$ replaced by $a_j$ and $\l$ replaced by $\l_j$), we obtain
\begin{equation}\label{eq:6imp3}
<\n J(\sum_{i=1}^m\alpha_i\varphi_{a_i, \l_i}+\sum_{r=1}^{\bar m}\beta_r(v_r-\ov{(v_r)}_Q))), \l_j\frac{\partial \varphi_{a_j, \l_j}}{\partial\l_j}>=2(E_{1, \l_j}+E_{0, \l_j}-16m\pi^2E_{2, \l_j}),
\end{equation}
where
\begin{equation}\label{eq:6imp1}
E_{1, \l_j}:=<P_g(\sum_{i=1}^m\alpha_i\varphi_{a_i, \l_i}), \l_j\frac{\partial \varphi_{a_j, \l_j}}{\partial\l_j}>,
\end{equation}
\begin{equation}\label{eq:e0lj}
E_{0, \l_j}=<P_g(\sum_{r=1}^{\bar m}\beta_r v_r),  \l_j\frac{\partial \varphi_{a_j, \l_j}}{\partial\l_j}>,
\end{equation}
and
\begin{equation}\label{eq:6imp2}
E_{2, \l_j}:=\frac{\int_M K(x)e^{4\sum_{i=1}^m\alpha_i\varphi_{a_i, \l_i}(x)+4\sum_{r=1}^{\bar m}\beta_r v_r(x)}\l_j\frac{\partial \varphi_{a_j, \l_j}(x)}{\partial\l_j}(x)dV_g(x)}{\int_M K(x)e^{4\sum_{i=1}^m\alpha_i\varphi_{a_i, \l_i}(z)+4\sum_{r=1}^{\bar m}\beta_r v_r(z)}dV_g(z)}.
\end{equation}
Now, we are going to estimate \;$E_{1, \l_j}$,\; $E_{0, \l_j}$,\;and \;$E_{2, \l_j}$. We start with \;$E_{1, \l_j}$. For this, we use \eqref{eq:6imp1}, the linearity of $P_g$ and of the inner product $<\cdot, \cdot>$, the fact that $d_g(a_i, a_j)\geq 4\ov C\eta$ for $i\neq j$, the second estimate of Lemma \ref{eq:selfintest}, the second estimate of Lemma \ref{eq:interactest}, the fact that $|\alpha_i-1|=O\left(\frac{1}{\l_i^{1-\gamma}\sqrt{\log\l_i}}\right)$ $i=1, \cdots, m$, and the $\l_i$'s are comparable, to get
\begin{equation}\label{eq:6imp5}
E_{1, \l_j}=16\alpha_j\pi^2+\alpha_j\frac{8\pi^2}{\l_j^2}H_2(a_j, a_j)+\sum_{i=1, i\neq j}^m \alpha_i\frac{8\pi^2}{\l_j^2}G_2(a_j, a_i)+O\left(\frac{1}{\l^3_j}\right).
\end{equation}
Now, we turn to the estimate of $E_{2, \l_j}$. For this, we first set
\begin{equation}\label{eq:6imp6}
\begin{split}
D:=\int_M K(x)e^{4\sum_{i=1}^m\alpha_i\varphi_{a_i, \l_i}(z)+4\sum_{r=1}^{\bar m}\beta_r v_r(z)}dV_g(z),\\\;\;\;E_{3, \l_j}:=\int_M K(x)e^{4\sum_{i=1}^m\alpha_i\varphi_{a_i, \l_i}(x)+4\sum_{r=1}^{\bar m}\beta_r v_r(x)}\l_j\frac{\partial \varphi_{a_j, \l_j}(x)}{\partial\l_j}(x)dV_g(x),
\end{split}
\end{equation}
and use \eqref{eq:6imp2} to have
\begin{equation}\label{eq:6imp7}
E_{2, \l_j}=\frac{E_{3, \l_j}}{D}.
\end{equation}
Next, let us estimate $E_{3, \l_j}$.  In order to do that, we first use the fact that  $d_{g}(a_i, a_j)\geq 4\ov C\eta$ for $i\neq j=1, \cdots, $, the properties of the metrics $g_{a_i}$ (see \eqref{eq:detga}-\eqref{eq:proua} with $a$ replaced by $a_i$) $i=1, \cdots, m$, the first  and second estimate of Lemma \ref{eq:outbubbleest}, \eqref{eq:decompG4}, \eqref{eq:regH4}, the fact that $\sum_{r=1}^{\bar m}|\beta_r|+\sum_{i=1}^m|\alpha_i-1|\sqrt{\log \l_i}=O\left(\sum_{i=1}^m\frac{1}{\l_i^{1-\gamma}}\right)$ and \eqref{eq:6imp6} to get
\begin{equation}\label{eq:6imp8}
E_{3, \l_j}=\sum_{k=1}^m\int_{B_{a_k}^{a_k}(\eta)} K(x)e^{4\sum_{i=1}^m\alpha_i\varphi_{a_i, \l_i}(x)+4\sum_{r=1}^{\bar m}\beta_r v_r(x)}\l_j\frac{\partial \varphi_{a_j, \l_j}(x)}{\partial\l_j}(x)dV_g(x)+O\left(1\right).
\end{equation}
Now, to continue we set
\begin{equation}\label{eq:6imp9}
E^k_{3, \l_j}=\int_{B_{a_k}^{a_k}(\eta)}K(x)e^{4\sum_{i=1}^m\alpha_i\varphi_{a_i, \l_i}(x)+4\sum_{r=1}^{\bar m}\beta_r v_r(x)}\l_j\frac{\partial \varphi_{a_j, \l_j}(x)}{\partial\l_j}(x)dV_g(x), \,\,\;k=1, \cdots, m,
\end{equation}
and use \eqref{eq:6imp8} to obtain
\begin{equation}\label{eq:6imp10}
E_{3, \l_j}=\sum_{k=1}^mE^k_{3, \l_j}+O\left(1\right).
\end{equation}
Furthermore, we write \eqref{eq:6imp10} into the following equivalent form
\begin{equation}\label{eq:6im10}
E_{3, \l_j}=E_{3, \l_j}^j+\sum_{k=1, k\neq j}^mE^k_{3, \l_j}+O\left(1\right).
\end{equation}
Now, we are going to estimate $E_{3, \l_j}^k$ for $k\neq j$, $k=1, \cdots, m$ and after $E_{3, \l_j}^j$. To do that, we use \eqref{eq:6imp9} and write it in the following equivalent form
\begin{equation}\label{eq:6imp11}
E_{3, \l_j}^k=\int_{B_{a_k}^{a_k}(\eta)}K(x)e^{4\alpha_k\varphi_{a_k, \l_k}}e^{4\sum_{i=1, i\neq k}^m\alpha_i\varphi_{a_i, \l_i}(x)+4\sum_{r=1}^{\bar m}\beta_r v_r(x)}\l_j\frac{\partial \varphi_{a_j, \l_j}(x)}{\partial\l_j}(x)dV_g(x).
\end{equation}
Next, using the fact that $k\neq j$, the properties of the metrics $g_{a_i}$ (see \eqref{eq:detga}-\eqref{eq:proua} with $a$ replaced by $a_i$ $i=1, \cdots, m$, the fact that $d_{g}(a_l, a_s)\geq 4\ov C\eta$ or $l, s=1,\cdots, m$, $l\neq s$,  \eqref{eq:decompG4}, \eqref{eq:regH4}, the first estimate of Lemma \ref{eq:intbubbleest}, the first estimate of Lemma \ref{eq:outbubbleest}, the second estimate of Lemma \ref{eq:outbubbleest}, the fact that $|\alpha_i-1|=O\left(\frac{1}{\l_i^{1-\gamma}\sqrt{\log \l_i}}\right)$ $i=1, \cdots, m$, and \eqref{eq:6imp11}, we obtain
\begin{equation}\label{eq:6imp12}
\begin{split}
E_{3, \l_j}^k=\int_{B_{a_k}^{a_k}(\eta)}K(x)\frac{\l_k^{8\alpha_k}}{(1+\l_k^2 d_{g_{a_k}}^2(a_k, x))^{4\alpha_k}}e^{4(\alpha_k H(a_k, x)+\frac{\alpha_k}{4\l_k^2}\D_{g_{a_k}} H(a_k, x))}e^{4\sum_{i=1, i\neq k}^m(\alpha_i G(a_i, x)+\frac{\alpha_i}{4\l_i^2}\D_{g_{a_i}}G(a_i, x))}\\\times e^{4\sum_{r=1}^{\bar m}\beta_r v_r(x)}\left(-\frac{1}{2\l_j^2}\D_{g_{a_j}}G(a_j, x)+O\left(\frac{1}{\l^3_j}\right)\right)dV_g(x)
\end{split}
\end{equation}
Now, arguing as in section \ref{eq:jexp}, we obtain
\begin{equation}\label{eq:6imp23}
E_{3, \l_j}^k=\gamma_k\left(-\frac{1}{2\l_j^2}\D_{g_{a_j}}G(a_j, a_k)+O\left(\frac{1}{\l^3_j}\right)\right), \;\;k=1,\cdots, m, k\neq j.
\end{equation}
Next, let us estimate $E_{3, \l_j}^j$. For this, we use \eqref{eq:6imp9} and write
\begin{equation}\label{eq:6imp24}
E_{3, \l_j}^j=\int_{B_{a_j}^{a_j}(\eta)}K(x)e^{4\alpha_j\varphi_{a_j, \l_j}(x)}e^{4\sum_{i=1, i\neq j}^m\alpha_i\varphi_{a_i, \l_i}(x)+4\sum_{r=1}^{\bar m}\beta_r v_r(x)}\l_j\frac{\partial \varphi_{a_j, \l_j}(x)}{\partial\l_j}(x)dV_g(x).
\end{equation}
Using the properties of the metrics $g_{a_i}$ $i=1, \cdots, m$ (see \eqref{eq:detga}-\eqref{eq:proua} with $a$ replaced by $a_i$),  the fact that $d_g(a_i, a_j)\geq 4\ov C\eta$ for $i\neq j$ $i=1, \cdots, m$, \eqref{eq:regH4}, the first  and second estimates of Lemma \ref{eq:intbubbleest}, and  the first estimate of Lemma \ref{eq:outbubbleest} combined with \eqref{eq:6imp24}, we get
\begin{equation}\label{eq:6imp25}
\begin{split}
E_{3, \l_j}^j=\int_{B_{a_j}^{a_j}(\eta)}K(x)\frac{\l_j^{8\alpha_j}}{(1+\l_j^2d^2_{g_{a_j}}(a_j, x))^{4\alpha_j}}e^{4(\alpha_j H(a_j, x)+\frac{\alpha_j}{4\l_j^2}\D_{g_{a_j}}H(a_j, x))}e^{4\sum_{i=1, i\neq j}^m(\alpha_i G(a_i, x)+\frac{\alpha_i}{4\l_i^2}\D_{g_{a_i}} G(a_i, x))}\\\times e^{4\sum_{r=1}^{\bar m}\beta_r v_r(x)}\left(\frac{2}{1+\l_j^2d_{g_{a_j}}^2(a_j, x)}-\frac{1}{2\l_j^2}\D_{g_{a_j}}H(a_j, x)+O\left(\frac{1}{\l^3_j}\right)\right)dV_g(x).
\end{split}
\end{equation}
Now, arguing again as in section \ref{eq:jexp}, we get
 \begin{equation}\label{eq:6imp44o1}
\begin{split}
E_{3, \l_j}^j=&\gamma_j^1\left(2-\frac{1}{\alpha_j}+\frac{\alpha_j}{\l_j^2}\D_{g_{a_j}}H(a_j, a_j)-\frac{1}{2\l_j^2}\D_{g_{a_j}}H(a_j, a_j)+\frac{1}{8\l_j^2}\left(\frac{\D_{g_{a_j}}\mathcal{F}_j^{A}(a_j)}{\mathcal{F}_j^{A}(a_j)}-\frac{2}{3}R_g(a_i)\right)\right)\\&+O\left(\sum_{i=1}^m|\alpha_i-1|\l_i^2+(\sum_{r=1}^{\bar m}|\beta_r|)(\sum_{i=1}^m\l_i^2)+\sum_{i=1}^m\l_i\right).
\end{split}
\end{equation}
Combining  \eqref{eq:5imp31}, \eqref{eq:5imp31a}, \eqref{eq:deficit}, the fact that the $\l_i$'s are comparable, \eqref{eq:6imp10}, \eqref{eq:6imp23}, and \eqref{eq:6imp44o1}, we obtain
\begin{equation}\label{eq:6imp45a}
\begin{split}
E_{3, \l_j}=&\gamma_j^1\left(1-\frac{1}{\alpha_j}+e^{+\frac{\alpha_j}{\l_j^2}\D_{g_{a_j}}H(a_j, a_j)}-\frac{1}{2\l_j^2}\D_{g_{a_j}} H(a_j, a_j)+\frac{1}{8\l_j^2}\left(\frac{\D_{g_{a_j}}\mathcal{F}_j^{A}(a_j)}{\mathcal{F}_j^{A}(a_j)}-\frac{2}{3}R_g(a_i)\right)\right)\\&-\sum_{i=1, i\neq j}^m\gamma_i\frac{1}{2\l_j^2}\D_{g_{a_j}} G(a_j, a_i)+O\left(\sum_{i=1}^m|\alpha_i-1|\l_i^2+(\sum_{r=1}^{\bar m}|\beta_r|)(\sum_{i=1}^m\l_i^2)+\sum_{i=1}^m\l_i\right).
\end{split}
\end{equation}
On the other hand, factorizing by $e^{+\frac{\alpha_j}{\l_j^2}\D_{g_{a_j}}H(a_j, a_j)}$ the first equation of the right hand side of \eqref{eq:6imp45a} and using \eqref{eq:decompG4}, \eqref{eq:regH4}, \eqref{eq:5imp31a}, and \eqref{eq:deficit}, we get
\begin{equation}\label{eq:6imp45b}
\begin{split}
E_{3, \l_j}=&\gamma_j\left((1-\frac{1}{\alpha_j})e^{-\frac{\alpha_j}{\l_j^2}\D_{g_{a_j}}H(a_j, a_j)}+1-\frac{1}{2\l_j^2}\D_{g_{a_j}} H(a_j, a_j)+\frac{1}{8\l_j^2}\left(\frac{\D_{g_{a_j}}\mathcal{F}_j^{A}(a_j)}{\mathcal{F}_j^{A}(a_j)}-\frac{2}{3}R_g(a_i)\right)\right)\\&-\sum_{i=1, i\neq j}^m\gamma_i\frac{1}{2\l_j^2}\D_{g_{a_j}}G(a_j, a_i)+O\left(\sum_{i=1}^m|\alpha_i-1|\l_i^2+(\sum_{r=1}^{\bar m}|\beta_r|)(\sum_{i=1}^m\l_i^2)+\sum_{i=1}^m\l_i\right).
\end{split}
\end{equation}
Furthermore, using \eqref{eq:regH4}, \eqref{eq:deficit}, and the assumption $|\alpha_j-1|=O\left(\frac{1}{\l_j^{1-\gamma}\sqrt{\log \l_j}}\right)$, it is easy to see that
\begin{equation}\label{eq:6imp45c}
(1-\frac{1}{\alpha_j})e^{-\frac{\alpha_j}{\l_j^2}\D_{g_{a_j}}H(a_j, a_j)}+1=2-\frac{1}{\alpha_j}+O\left(\frac{|\alpha_j-1|}{\l_j^2}+\frac{1}{\l^3_j}\right).
\end{equation}
Thus, combining \eqref{eq:5imp31a}, \eqref{eq:6imp45b} and \eqref{eq:6imp45c}, we derive
\begin{equation}\label{eq:6imp45}
\begin{split}
E_{3, \l_j}=&\gamma_j\left(2-\frac{1}{\alpha_j}-\frac{1}{2\l_j^2}\D_{g_{a_j}}H(a_j, a_j)+\frac{1}{8\l_j^2}\left(\frac{\D_{g_{a_j}}\mathcal{F}_j^{A}(a_j)}{\mathcal{F}_j^{A}(a_j)}-\frac{2}{3}R_g(a_i)\right)\right)\\&-\sum_{i=1, i\neq j}^m\gamma_i\frac{1}{2\l_j^2}\D_{g_{a_j}}G(a_j, a_i)+O\left(\sum_{i=1}^m|\alpha_i-1|\l_i^2+(\sum_{r=1}^{\bar m}|\beta_r|)(\sum_{i=1}^m\l_i^2)+\sum_{i=1}^m\l_i\right).
\end{split}
\end{equation}
So, using \eqref{eq:5imp31a}, \eqref{eq:6imp7} and \eqref{eq:6imp45}, we infer that
\begin{equation}\label{eq:6imp46}
\begin{split}
E_{2, \l_j}=&\frac{\gamma_j}{D}\left(2-\frac{1}{\alpha_j}-\frac{1}{2\l_j^2}\D_{g_{a_j}}H(a_j, a_j)+\frac{1}{8\l_j^2}\left(\frac{\D_{g_{a_j}}\mathcal{F}_j^{A}(a_j)}{\mathcal{F}_j^{A}(a_j)}-\frac{2}{3}R_g(a_i)\right)\right)\\&+\frac{\gamma_j}{D}O\left(\sum_{i=1}^m\frac{|\alpha_i-1|}{\l_i^2}+(\sum_{r=1}^m|\beta_r|)(\sum_{i=1}^m\frac{1}{\l_i^2})+\sum_{i=1}^m\frac{1}{\l^3_i}\right)\\&-\sum_{i=1, i\neq j}^m\frac{\gamma_i}{D}\left(\frac{1}{2\l_j^2}\D_{g_{a_j}}G(a_j, a_i)\right).
\end{split}
\end{equation}
Next, we are going to estimate $E_{0, \l_j}$. For this, we use the fact that the $v_r$'s ($r=1, \cdots, \bar m$) are eigenfunctions of $P_g$ and Lemma \ref{eq:intbubbleest}, to infer that
\begin{equation}\label{eq:e0ljestfin}
E_{0, \l_j}=O\left(\sum_{r=1}^{\bar m}|\beta_r|)(\sum_{i=1}^m\frac{1}{\l_i^2})\right).
\end{equation}
Now, using the relation $2-\frac{1}{\alpha_j}-\alpha_j=-\frac{(\alpha_j-1)^2}{\alpha_j}$, $|\alpha_j-1|=O\left(\frac{1}{\l_j^{1-\gamma}\sqrt{\log \l_j}}\right)$, \eqref{eq:6imp5}, \eqref{eq:6imp46}, \eqref{eq:e0ljestfin},  and setting
\begin{equation}\label{eq:6imp53}
\tau_i:=1-\frac{m\gamma_i}{D}, \,\;\;i=1, \cdots, m,
\end{equation}
we obtain
\begin{equation}\label{eq:6imp54}
\begin{split}
E_{1, \l_j}+E_{0, \l_j}-16m\pi^2 E_{2, \l_j}=&16\pi^2\alpha_j\tau_j-\frac{8\pi^2}{\l_j^2}\tau_j \D_{g_{a_j}}H(a_j, a_j)-\frac{2\pi^2}{\l_j^2}\left(\frac{\D_{g_{a_j}}\mathcal{F}_j^{A}(a_j)}{\mathcal{F}_j^{A}(a_j)}-\frac{2}{3}R_g(a_j)\right)\\&-\frac{8\pi^2}{\l_j^2}\sum_{i=1, i\neq j}^m \tau_i \D_{g_{a_j}} G(a_j, a_i)+\tau_j\frac{2\pi^2}{\l_j^2}\left(\frac{\D_{g_{a_j}}\mathcal{F}_j^{A}(a_j)}{\mathcal{F}_j^{A}(a_j)}-\frac{2}{3}R_g(a_j)\right)\\&+O\left(\sum_{i=1}^m|\alpha_i-1|^2+\sum_{i=1}^m\frac{|\alpha_i-1|}{\l_i^2}+(\sum_{r=1}^m|\beta_r|)(\sum_{i=1}^m\frac{1}{\l_i^2})+\sum_{i=1}^m\frac{1}{\l^3_i}\right).
\end{split}
\end{equation}
Hence, combining \eqref{eq:6imp3} , \eqref{eq:6imp54}, and using Young's inequality, we obtain
\begin{equation*}
\begin{split}
&<\n J(\sum_{i=1}^m\alpha_i\varphi_{a_i, \l_i}+\sum_{r=1}^{\bar m}\beta_r (v_r-\ov{(v_r)}_Q)), \l_j\frac{\partial \varphi_{a_j, \l_j}}{\partial\l_j}>=32\pi^2\alpha_j\tau_j-\frac{16\pi^2}{\l_j^2}\tau_j \D_{g_{a_j}}H(a_j, a_j)\\&-\frac{4\pi^2}{\l_j^2}\left(\frac{\D_{g_{a_j}}\mathcal{F}_j^{A}(a_j)}{\mathcal{F}_j^{A}(a_j)}-\frac{2}{3}R_g(a_j)\right)-\frac{16\pi^2}{\l_j^2}\sum_{i=1, i\neq j}^m \tau_i \D_{g_{a_j}}G(a_j, a_i)+\tau_j\frac{4\pi^2}{\l_j^2}\left(\frac{\D_{g_{a_j}}\mathcal{F}_j^{A}(a_j)}{\mathcal{F}_j^{A}(a_j)}-\frac{2}{3}R_g(a_j)\right)\\&+O\left(\sum_{i=1}^m|\alpha_i-1|^2+\sum_{r=1}^{\bar m}|\beta_r|^3+\sum_{i=1}^m\frac{1}{\l^3_i}\right).
\end{split}
\end{equation*}
thereby ending the proof of the Lemma.
\end{pf}
\vspace{6pt}

\noindent
Lemma \ref{eq:gradientlambdaest} implies the following corollary.
\begin{cor}\label{eq:cgradientlambdaest}(Corollary gradient estimate with respect to $\bar\l$)\\
Assuming that $\eta$ is a small positive real number with $0<2\eta<\varrho$ where $\varrho$ is as in \eqref{eq:cutoff},  and $0<\epsilon\leq \epsilon_0$ where $\epsilon_0$ is as in \eqref{eq:mini}, then for $a_i\in M$ concentration points,  $\alpha_i$ masses, $\l_i$ concentration parameters ($i=1,\cdots,m$), and $\beta_r$ negativity parameters  ($r=1, \cdots, m$)  satisfying \eqref{eq:afpara}, we have
\begin{equation*}
\begin{split}
&<\n J(\sum_{i=1}^m\alpha_i\varphi_{a_i, \l_i}+\sum_{r=1}^{\bar m}\beta_r(v_r-\ov{(v_r)}_Q)), \sum_{i=1}^m\frac{\l_i}{\alpha_i}\frac{\partial \varphi_{a_i, \l_i}}{\partial \l_i}>= 4\pi^2\sum_{i=1}^m\frac{1}{\l_i^2}\left(\frac{\D_{g_{a_i}} \mathcal{F}^{A}_i(a_i)}{\mathcal{F}^{A}_i(a_i)}-\frac{2}{3}R_g(a_i)\right)\\&+O\left(\sum_{i=1}^m|\alpha_i-1|^2+\sum_{r=1}^{\bar m}|\beta_r|^3+\sum^{m}_{i=1}\tau^3_i+\sum_{i=1}^m\frac{1}{\l_i^3}\right),
\end{split}
\end{equation*}
where $A:=(a_1, \cdots, a_m)$, $O\left(1\right)$ as as in Lemma \ref{eq:energyest}, and for $i=1, \cdots, m$, $\tau_i$ is as in Lemma \ref{eq:gradientlambdaest}.
\end{cor}
\vspace{6pt}
\begin{pf}
Using Lemma \ref{eq:gradientlambdaest}, the fact that $|\alpha_j-1|=O\left(\frac{1}{\l_j^{1-\gamma}\sqrt{\log \l_j}}\right)$, $j=1, \cdots, m$, and Young's inequality, we get
\begin{equation}\label{eq:coreql1}
\begin{split}
<\n J(\sum_{i=1}^m\alpha_i\varphi_{a_i, \l_i}+\sum_{r=1}^{\bar m}\beta_r(v_r-\ov{(v_r)}_Q)), \sum_{i=1}^m\frac{\l_i}{\alpha_i}\frac{\partial \varphi_{a_i, \l_i}}{\partial \l_i}>=-4\pi^2\sum_{i=1}^m\frac{1}{\l_i^2}\left(\frac{\D_{g_{a_i}} \mathcal{F}^{A}_i(a_i)}{\mathcal{F}^{A}_i(a_i)}-\frac{2}{3}R_g(a_i)\right)\\+32\pi^2\sum_{i=1}^m\tau_i+O\left(\sum_{i=1}^m|\alpha_i-1|^2+\sum_{r=1}^{\bar m}|\beta_r |^3+\sum^{m}_{i=1}\tau^2_3+\sum_{i=1}^m\frac{1}{\l_i^3}\right).
\end{split}
\end{equation}
Now, let us estimate $\sum_{j=1}^m\tau_j$. For this, we start by estimation $D$ in terms of $\Gamma$. Indeed, using \eqref{eq:5imp43} and \eqref{eq:6imp6}, we get
\begin{equation}\label{eq:coreql2}
D=\Gamma\left(1+\sum_{i=1}^m\frac{\gamma_i}{4\Gamma\l_i^2}\left(\frac{\D_{g_{a_i}} \mathcal{F}^{A}_i(a_i)}{\mathcal{F}^{A}_i(a_i)}-\frac{2}{3}R_g(a_i)\right)+O\left(\sum_{i=1}^m\frac{|\alpha_i-1|}{\l_i^2}+(\sum_{r=1}^{\bar m}|\beta_r|)(\sum_{i=1}^m\frac{1}{\l_i^2})+\sum_{i=1}^m\frac{1}{\l_i^3}\right)\right).
\end{equation}
Next, it is easy to see that \eqref{eq:coreql2}, implies
\begin{equation}\label{eq:coreql3}
D=\Gamma\left(1+O\left(\sum_{i=1}^m\frac{1}{\l_i^2}\right)\right).
\end{equation}
Using \eqref{eq:5imp42}, \eqref{eq:6imp53}, \eqref{eq:coreql2}, and \eqref{eq:coreql3}, we infer that
\begin{equation}\label{eq:coreql6}
\tau_j=\tilde \tau_j+\frac{\gamma_j}{\Gamma}\sum_{i=1}^m\frac{m\gamma_i}{4D\l_i^2}\left(\frac{\D_{g_{a_i}} \mathcal{F}^{A}_i(a_i)}{\mathcal{F}^{A}_i(a_i)}-\frac{2}{3}R_g(a_i)\right)+O\left(\sum_{i=1}^m\frac{|\alpha_i-1|}{\l_i^2}+(\sum_{r=1}^{\bar m}|\beta_r|)(\sum_{i=1}^m\frac{1}{\l_i^2})+\sum_{i=1}^m\frac{1}{\l_i^3}\right).
\end{equation}
Now, we emphasize the following estimate which will be useful later and follows directly from \eqref{eq:coreql6}
\begin{equation}\label{eq:coreql7}
\tau_j=\tilde\tau_j+O\left(\sum_{i=1}\frac{1}{\l_i^2}\right)
\end{equation}
Coming back to our goal of estimating $\sum_{j=1}^m\tau_j$, we use  \eqref{eq:5im38}, \eqref{eq:5imp42},  \eqref{eq:6imp53}, and \eqref{eq:coreql6}, to derive that
\begin{equation}\label{eq:coreql9}
\begin{split}
\sum_{j=1}^m\tau_j=&\sum_{i=1}^m\frac{1}{4\l_i^2}\left(\frac{\D_{g_{a_i}} \mathcal{F}^{A}_i(a_i)}{\mathcal{F}^{A}_i(a_i)}-\frac{2}{3}R_g(a_i)\right)-\sum_{i=1}^m\frac{\tau_i}{4\l_i^2}\left(\frac{\D_{g_{a_i}} \mathcal{F}^{A}_i(a_i)}{\mathcal{F}^{A}_i(a_i)}-\frac{2}{3}R_g(a_i)\right)\\&+O\left(\sum_{i=1}^m\frac{|\alpha_i-1|}{\l_i^2}+(\sum_{r=1}^{\bar m}|\beta_r|)(\sum_{i=1}^m\frac{1}{\l_i^2})+\sum_{i=1}^m\frac{1}{\l_i^3}\right).
\end{split}
\end{equation}
Next, combining \eqref{eq:coreql1} , \eqref{eq:coreql9},  and using again Young's inequality, we arrive to
\begin{equation}\label{eq:coreql10}
\begin{split}
&<\n J(\sum_{i=1}^m\alpha_i\varphi_{a_i, \l_i}+\sum_{r=1}^{\bar m}\beta_r(v_r-\ov{(v_r)}_Q)), \sum_{i=1}^m\frac{\l_i}{\alpha_i}\frac{\partial \varphi_{a_i, \l_i}}{\partial \l_i}>=4\pi^2\sum_{i=1}^m\frac{1}{\l_i^2}\left(\frac{\D_{g_{a_i}} \mathcal{F}^{A}_i(a_i)}{\mathcal{F}^{A}_i(a_i)}-\frac{2}{3}R_g(a_i)\right)\\&+O\left(\sum_{i=1}^m|\alpha_i-1|^2+\sum_{r=1}^{\bar m}|\beta_r|^3+\sum^{m}_{i=1}\tau^3_i+\sum_{i=1}^m\frac{1}{\l_i^3}\right),
\end{split}
\end{equation}
thereby ending the proof of the Corollary.
\end{pf}
\vspace{6pt}

\noindent
Now, we turn to the gradient estimate with respect to the $\alpha_i$'s. Precisely, we obtain:
\begin{lem}\label{eq:gradientalpha} (Gradient estimate with respect to $\bar \alpha$)\\
Assuming that $\eta$ is a small positive real number with $0<2\eta<\varrho$ where $\varrho$ is as in \eqref{eq:cutoff}, and $0<\epsilon\leq \epsilon_0$ where $\epsilon_0$ is as in \eqref{eq:mini}, then  for $a_i\in M$ concentration points,  $\alpha_i$ masses, $\l_i$ concentration parameters ($i=1,\cdots,m$), and $\beta_r$ negativity parameters ($r=1, \cdots, \bar m$) satisfying \eqref{eq:afpara}, we have that for every $j=1, \cdots, m$, there holds
\begin{equation*}
\begin{split}
&<\n J(\sum_{i=1}^m\alpha_i\varphi_{a_i, \l_i}+\sum_{r=1}^{\bar m}\beta_r(v_r-\ov{(v_r)}_Q), \varphi_{a_j, \l_j}>=\\&\left(2\log\l_j+H(a_j, a_j)-\frac{5}{6}\right)\frac{1}{\alpha_j}<\n J(\sum_{i=1}^m\alpha_i\varphi_{a_i, \l_i}+\sum_{r=1}^{\bar m}\beta_r(v_r-\ov{(v_r)}_Q), \l_j\frac{\partial \varphi_{a_j, \l_j}}{\partial \l_j}>\\&+\sum_{i=1, i\neq j}^m G(a_j, a_i)<\n J(\sum_{i=1}^m\alpha_i\varphi_{a_i, \l_i}+\sum_{r=1}^{\bar m}\beta_r(v_r-\ov{(v_r)}_Q), \l_i\frac{\partial \varphi_{a_i, \l_i}}{\partial \l_i}>\\&+64\pi^2(\alpha_j-1)\log\l_j+O\left(\log\l_j\left(\sum_{i=1}^m\frac{|\alpha_i-1|}{\log \l_i}+(\sum_{r=1}^{\bar m}|\beta_r|)(\sum_{i=1}^m\frac{1}{\log \l_i})+\sum_{i=1}^m\frac{1}{\l^2_i}\right)\right).
\end{split}
\end{equation*}
\end{lem}
\begin{pf}
First of all, using the same arguments as in Lemma \ref{eq:gradientlambdaest}, we obtain
\begin{equation}\label{eq:8imp3}
<\n J(\sum_{i=1}^m\alpha_i\varphi_{a_i, \l_i}+\sum_{r=1}^{\bar m}\beta_r(v_r-\ov{(v_r)}_Q)), \varphi_{a_j, \l_j}>=2(E_{1, \alpha_j}+E_{0, \alpha_j}-16m\pi^2E_{2, \alpha_j}).
\end{equation}
where
\begin{equation}\label{eq:8imp1}
E_{1, \alpha_j}:=<P_g(\sum_{i=1}^m\alpha_i\varphi_{a_i, \l_i}), \varphi_{a_j, \l_j}>,
\end{equation}
\begin{equation}\label{eq:e0alphaj}
E_{0, \alpha_j}:=<P_g(\sum_{r=1}^{\bar m}\beta_r v_r), \varphi_{a_j, \l_j}>
\end{equation}
and
\begin{equation}\label{eq:8imp2}
E_{2, \alpha_j}:=\frac{\int_M K(x)e^{4\sum_{i=1}^m\alpha_i\varphi_{a_i, \l_i}(x)+4\sum_{r=1}^{\bar m}\beta_r v_r(x)}\varphi_{a_j, \l_j}(x)dV_g(x)}{\int_M K(z)e^{4\sum_{i=1}^m\alpha_i\varphi_{a_i, \l_i}(z)+\sum_{r=1}^{\bar m}\beta_r v_r(z)}dV_g(z)},
\end{equation}
Now, we are going to estimate $E_{1, \alpha_j}$, $E_{0, \alpha_j}$, and $E_{2, \alpha_j}$. We start with $E_{1, \alpha_j}$. For this, we argue again as in Lemma \ref{eq:gradientlambdaest} and obtain
\begin{equation}\label{eq:8imp5a}
\begin{split}
E_{1, \alpha_j}=&\alpha_j\left(32\pi^2\log \l_j-\frac{40}{3}\pi^2+16\pi^2 H_(a_j, a_j)\right)\\&+
16\pi^2\sum_{i=1, i\neq j}^m\alpha_i G(a_j, a_i)+O\left(\sum_{i=1}^m\frac{1}{\l^2_i}\right).
\end{split}
\end{equation}
Next, we turn to the estimate of $E_{2, \alpha_j}$. For this, we first set
\begin{equation}\label{eq:8imp6}
E_{3, \alpha_j}:=\int_M K(x)e^{4\sum_{i=1}^m\alpha_i\varphi_{a_i, \l_i}(x)+4\sum_{r=1}^{\bar m}\beta_r v_r(x)}\varphi_{a_j, \l_j}(x)dV_g(x),
\end{equation}
and use \eqref{eq:8imp2} to have
\begin{equation}\label{eq:8imp7}
E_{2, \alpha_j}=\frac{E_{3, \alpha_j}}{D},
\end{equation}
where $D$ is as in \eqref{eq:6imp6}. Now, arguing again as in Lemma \ref{eq:gradientlambdaest}, we obtain
\begin{equation}\label{eq:8imp10a}
E_{3, \alpha_j}=E_{3, \alpha_j}^j+\sum_{k=1, k\neq j}^mE^k_{3, \alpha_j}+O\left(1\right),
\end{equation}
where
\begin{equation}\label{eq:8imp9}
E^k_{3, \alpha_j}=\int_{B_{a_k}^{a_k}(\eta)}K(x)e^{4\sum_{i=1}^m\alpha_i\varphi_{a_i, \l_i}(x)+4\sum_{r=1}^{\bar m}\beta_r v_r(x)}\varphi_{a_j, \l_j}(x)dV_g(x), \;\;\;k=1, \cdots, m.
\end{equation}
Next, we are going to estimate \;$E_{3, \alpha_j}^k$ for $k\neq j$, $k=1, \cdots, m$, and after $E_{3, \alpha_j}^j$. To do that, we argue again as in Lemma \ref{eq:gradientlambdaest} and get
\begin{equation}\label{eq:8imp23}
E_{3, \alpha_j}^k=\gamma_k\left(G(a_j, a_k)+O\left(\frac{1}{\l^2_k}\right)\right), \;\;k=1,\cdots,m, k\neq j.
\end{equation}
Now, let us estimate $E_{3, \alpha_j}^j$. For this, we argue again as in Lemma \ref{eq:gradientlambdaest} and get
\begin{equation}\label{eq:8imp60a}
\begin{split}
E_{3, \alpha_j}^j&= \gamma_j\left(2\log\l_j-s_j+\frac{\log\l_j}{2\l_j^2}\left(\frac{\D_{g_{a_j}}\mathcal{F}_j^{A}(a_j)}{\mathcal{F}^{A}_j(a_j)}-\frac{2}{3}R_g(a_j)\right)+H(a_j, a_j)\right)\\&+\gamma_j O\left(\log\l_j\left(\sum_{i=1}^m\frac{|\alpha_i-1|}{\l_i^2}+(\sum_{r=1}^{\bar m}|\beta_r|)(\sum_{i=1}^m\frac{1}{\l_i^2})+\sum_{i=1}^m\frac{1}{\l^3_i}\right)\right),
\end{split}
\end{equation}
where
\begin{equation}\label{eq:8imp44}
s_j=2(4\alpha_j-1)(2\alpha_j-1)(\frac{1}{(4\alpha_j-2)^2}-\frac{1}{(4\alpha_j-1)^2}),
\end{equation}
Combining \eqref{eq:5imp31a}, \eqref{eq:8imp10a}, \eqref{eq:8imp23} and \eqref{eq:8imp60a}, we get
\begin{equation}\label{eq:8imp60b}
\begin{split}
E_{3, \alpha_j}&= \gamma_j\left(2\log\l_j-s_j+\frac{\log\l_j}{2\l_j^2}\left(\frac{\D_{g_{a_j}}\mathcal{F}_j^{A}(a_j)}{\mathcal{F}^{A}_j(a_j)}-\frac{2}{3}R_g(a_j)\right)+H(a_j, a_j)\right)\\&+O\left(\log\l_j\left(\sum_{i=1}^m\frac{|\alpha_i-1|}{\l_i^2}+(\sum_{r=1}^{\bar m}|\beta_r|)(\sum_{i=1}^m\frac{1}{\l_i^2})+\sum_{i=1}^m\frac{1}{\l^3_i}\right)\right)\\&+\sum_{i=1, i\neq j}^m \gamma_i\left(G(a_j, a_i)+O\left(\frac{1}{\l^2_i}\right)\right).
\end{split}
\end{equation}
Next, we are going to estimate $E_{0, \alpha_j}$. For this, we use \eqref{eq:e0alphaj} and the fact that the $v_r$'s are eigenfunctions of $P_g$ to infer that
\begin{equation}\label{eq:e0alphajest}
E_{0, \alpha_j}=O\left(\sum_{r=1}^{\bar m}|\beta_r|\right).
\end{equation}
Now, using \eqref{eq:8imp3}, \eqref{eq:8imp5a}, \eqref{eq:8imp7}, \eqref{eq:8imp44}, \eqref{eq:8imp60b}, \eqref{eq:e0alphajest}, and the fact that $|\alpha_j-1|=O\left(\frac{1}{\l_j^{1-\gamma}\sqrt{\log\l_j}}\right)$, we get
\begin{equation}\label{eq:finp4}
\begin{split}
&<\n J(\sum_{i=1}^m\alpha_i\varphi_{a_i, \l_i}+\sum_{r=1}^{\bar m}\beta_r(v_r-\ov{(v_r)}_Q), \varphi_{a_j, \l_j}>=64\pi^2\tau_j\log\l_j+32\pi^2\tau_j H(a_j, a_j)\\&+32\pi^2\sum_{i=1, i\neq j}^m \tau_i G(a_j, a_i)-\frac{80\pi^2}{3}\tau_j\pi^2+64\pi^2(\alpha_j-1)\log\l_j\\&+O\left(\log\l_j\left(\sum_{i=1}^m\frac{|\alpha_i-1|}{\log \l_i}+(\sum_{r=1}^{\bar m}|\beta_r|)(\sum_{i=1}^m\frac{1}{\log \l_i})+\sum_{i=1}^m\frac{1}{\l^2_i}\right)\right).
\end{split}
\end{equation}
Thus combining Lemma \ref{eq:gradientlambdaest} , \eqref{eq:finp4}, the fact that $|\alpha_i-1|=O\left(\frac{1}{\l_i^{1-\gamma}\sqrt{\log\l_i}}\right)$  $i=1, \cdots, m$ and the fact that the $\l_i$' are comparable, we get
\begin{equation*}
\begin{split}
&<\n J(\sum_{i=1}^m\alpha_i\varphi_{a_i, \l_i}+\sum_{r=1}^{\bar m}\beta_r(v_r-\ov{(v_r)}_Q), \varphi_{a_j, \l_j}>=\\&\left(2\log\l_j+H(a_j, a_j)-\frac{5}{6}\right)\frac{1}{\alpha_j}<\n J(\sum_{i=1}^m\alpha_i\varphi_{a_i, \l_i}+\sum_{r=1}^{\bar m}\beta_r(v_r-\ov{(v_r}_Q), \l_j\frac{\partial \varphi_{a_j, \l_j}}{\partial \l_j}>\\&+\sum_{i=1, i\neq j}^m G(a_j, a_i)<\n J(\sum_{i=1}^m\alpha_i\varphi_{a_i, \l_i}+\sum_{r=1}^{\bar m}\beta_r(v_r-\ov{(v_r)}_Q), \l_i\frac{\partial \varphi_{a_i, \l_i}}{\partial \l_i}>\\&+64\pi^2(\alpha_j-1)\log\l_j+O\left(\log\l_j\left(\sum_{i=1}^m\frac{|\alpha_i-1|}{\log \l_i}+(\sum_{r=1}^{\bar m}|\beta_r|)(\sum_{i=1}^m\frac{1}{\log \l_i})+\sum_{i=1}^m\frac{1}{\l^2_i}\right)\right).
\end{split}
\end{equation*}
thereby ending the proof of the Lemma.
\end{pf}
\vspace{6pt}

\noindent
Next, we are going to derive a gradient estimate for $J$ in the direction of the $a_i$'s. Indeed, we will show:
\begin{lem}\label{eq:gradientaest}(Gradient estimate with respect to $A$)\\
Assuming that $\eta$ is a small positive real number with $0<2\eta<\varrho$ where $\varrho$ is as in \eqref{eq:cutoff}, and $0<\epsilon\leq \epsilon_0$ where $\epsilon_0$ is as in \eqref{eq:mini}, then  for $a_i\in M$ concentration points,  $\alpha_i$ masses, $\l_i$ concentration parameters ($i=1,\cdots,m$), ad $\beta_r$ negativity parameters ($r=1, \cdots, \bar m$) satisfying \eqref{eq:afpara}, we have that for every $j=1, \cdots, m$, there holds
\begin{equation*}
\begin{split}
<\n J(\sum_{i=1}^m\alpha_i\varphi_{a_i, \l_i}+\sum_{r=1}^{\bar m}\beta_r(v_r-\ov{(v_r)}_Q), \frac{1}{\l_j}\frac{\partial \varphi_{a_j, \l_j}}{\partial a_j}>&=-\frac{8\pi^2}{\l_j}\frac{\n_g\mathcal{F}_j^{A}(a_j)}{\mathcal{F}_j^{A}(a_j)}\\&+O\left(\sum_{i=1}^m|\alpha_i-1|^2+\sum_{i=1}^m\frac{1}{\l_i^2}+\sum_{r=1}^{\bar m}|\beta_r|^2+\sum_{i=1}^m\tau_i^2\right),
\end{split}
\end{equation*}
where $A:=(a_1, \cdots, a_m)$, $O(1)$ is as in Lemma \ref{eq:energyest} and for  $i=1, \cdots, m$, \;$\tau_i$ is as in Lemma \ref{eq:gradientlambdaest}.
\end{lem}
\begin{pf}
Arguing again as in Lemma \ref{eq:gradientlambdaest}, we obtain
\begin{equation}\label{eq:7imp3}
<\n J(\sum_{i=1}^m\alpha_i\varphi_{a_i, \l_i}+\sum_{r=1}^{\bar m}\beta_r (v_r-\ov{(v_r)}_Q)), \frac{1}{\l_j}\frac{\partial \varphi_{a_j, \l_j}}{\partial a_j}>=2(E_{1, a_j}+E_{0, a_j}-16m\pi^2E_{2, a_j}).
\end{equation}
where
\begin{equation}\label{eq:7imp1}
E_{1, a_j}:=<P_g(\sum_{i=1}^m\alpha_i\varphi_{a_i, \l_i}), \frac{1}{\l_j}\frac{\partial \varphi_{a_j, \l_j}}{\partial a_j}>,
\end{equation}
\begin{equation}\label{eq:e0aj}
E_{0, a_j}:=<P_g(\sum_{r=1}^{\bar m}\beta_r v_r), \frac{1}{\l_j}\frac{\partial \varphi_{a_j, \l_j}}{\partial a_j}>,
\end{equation}
and
\begin{equation}\label{eq:7imp2}
E_{2, a_j}:=\frac{\int_M K(x)e^{4\sum_{i=1}^m\alpha_i\varphi_{a_i, \l_i}(x)+4\sum_{r=1}^{\bar m}\beta_r v_r(x)}\frac{1}{\l_j}\frac{\partial \varphi_{a_j, \l_j}(x)}{\partial a_j}dV_g(x)}{\int_M K(z)e^{4\sum_{i=1}^m\alpha_i\varphi_{a_i, \l_i}(z)+4\sum_{r=1}^{\bar m}\beta_r v_r(z)}dV_g(z)}.
\end{equation}
Now, we are going to estimate $E_{1, a_j}$, $E_{0, a_j}$ and $E_{2, a_j}$. We start with $E_{1, a_j}$. For this, we use again the same arguments as in Lemma \ref{eq:gradientlambdaest} and get
\begin{equation}\label{eq:7imp5a}
E_{1, a_j}=\frac{16\pi^2}{\l_j}\frac{\partial H(a_j, a_j)}{\partial a_j}+\sum_{i=1, i\neq j}^m \frac{16\pi^2}{\l_j}\frac{\partial G(a_j, a_i)}{\partial a_j}+O\left(\sum_{i=1}^m\frac{|\alpha_i-1|}{\l_i}+\sum_{i=1}^m\frac{1}{\l^3_i}\right)
\end{equation}
Next, we turn to the estimate of $E_{2, a_j}$. For this, we first set
\begin{equation}\label{eq:7imp6}
E_{3, a_j}:=\int_M K(x)e^{4\sum_{i=1}^m\alpha_i\varphi_{a_i, \l_i}(x)+4\sum_{r=1}^{\bar m}\beta_r v_r(x)}\frac{1}{\l_j}\frac{\partial \varphi_{a_j, \l_j}(x)}{\partial a_j}dV_g(x),
\end{equation}
and use \eqref{eq:6imp6} to have
\begin{equation}\label{eq:7imp7}
E_{2, a_j}=\frac{E_{3, a_j}}{D},
\end{equation}
where $D$ is as in \eqref{eq:6imp6}. Now, let us estimate $E_{3, a_j}$.  To do that, we argue again as in Lemma \ref{eq:gradientlambdaest}, and obtain
\begin{equation}\label{eq:7im10a}
E_{3, a_j}=E_{3, a_j}^j+\sum_{k=1, k\neq j}^mE^k_{3, a_j}+O\left(1\right),
\end{equation}
where
\begin{equation}\label{eq:7imp9}
E^k_{3, a_j}=\int_{B_{a_k}^{a_k}(\eta)} K(x)e^{4\sum_{i=1}^m\alpha_i\varphi_{a_i, \l_i}(x)+4\sum_{r=1}^{\bar m}\beta_r v_r(x)}\frac{1}{\l_j}\frac{\partial \varphi_{a_j, \l_j}(x)}{\partial a_j}(x)dV_g(x), \;\;\;k=1, \cdots, m,
\end{equation}
Next, we are going to estimate $E_{3, a_j}^k$ for $k\neq j$, $k=1, \cdots, m$ and after $E_{3, a_j}^j$. To do that, we use again the same strategy as in Lemma \ref{eq:gradientlambdaest} and get
\begin{equation}\label{eq:7imp23}
E_{3, a_j}^k=\gamma_k\left(\frac{1}{\l_j}\frac{\partial G(a_j, a_k)}{\partial a_j}+O\left(\sum_{i=1}^m\frac{1}{\l^3_i}\right)\right), \;\;k=1,\cdots, m, k\neq j.
\end{equation}
Now, let us estimate $E_{3, a_j}^j$. For this, we use  again the same arguments as in Lemma \ref{eq:gradientlambdaest} and obtain
\begin{equation}\label{eq:7imp44}
\begin{split}
E_{3, a_j}^j=\gamma_j\left(\frac{1}{\l_j}\frac{\partial H(a_j, a_j)}{\partial a_j}+\frac{1}{4\l_j}\frac{\n_g\mathcal{F}_j^{A}(a_j)}{\mathcal{F}_j^{A}(a_j)}+O\left(\sum_{i=1}^m\frac{|\alpha_i-1|}{\l_i}+(\sum_{r=1}^{\bar m}|\beta_r|)(\sum_{i=1}^m\frac{1}{\l_i})+\sum_{i=1}^m\frac{1}{\l_i^3}\right)\right).
\end{split}
\end{equation}
Next, appealing again to \eqref{eq:5imp31a} (with $i$ replaced by $j$), and using \eqref{eq:7im10a}, \eqref{eq:7imp23}, and \eqref{eq:7imp44}, we obtain
\begin{equation}\label{eq:7imp45}
\begin{split}
E_{3, a_j}=&\gamma_j\left(\frac{1}{\l_j}\frac{\partial H(a_j, a_j)}{\partial a_j}+\frac{1}{4\l_j}\frac{\n_g\mathcal{F}_j^{A}(a_j)}{\mathcal{F}_j^{A}(a_j)}+O\left(\sum_{i=1}^m\frac{|\alpha_i-1|}{\l_i}+(\sum_{r=1}^{\bar m}|\beta_r|)(\sum_{i=1}^m\frac{1}{\l_i})+\sum_{i=1}^m\frac{1}{\l_i^3}\right)\right)\\&+\sum_{i=1, i\neq j}^m\gamma_i\left(\frac{1}{\l_j}\frac{\partial G(a_j, a_i)}{\partial a_j}+O\left(\sum_{i=1}^m\frac{1}{\l^3_i}\right)\right).
\end{split}
\end{equation}
Now, we are going to estimate $E_{0, a_j}$. For this, we use Lemma \ref{eq:intbubbleest}, \eqref{eq:e0aj} and the fact that the $v_r$'s are eigenfunctions of $P_g$ to obtain
\begin{equation}\label{eq:e0ajest1}
E_{0, a_j}=O\left((\sum_{r=1}^{\bar m}|\beta_r|)(\sum_{i=1}^m\frac{1}{\l_i})\right).
\end{equation}
Hence, combining \eqref{eq:7imp3}, \eqref{eq:7imp5a}, \eqref{eq:7imp7}, \eqref{eq:7imp45}, \eqref{eq:e0ajest1}, and Cauchy inequality we obtain
\begin{equation*}
\begin{split}
<\n J(\sum_{i=1}^m\alpha_i\varphi_{a_i, \l_i}+\sum_{r=1}^{\bar m}\beta_r(v_r-\ov{(v_r)}_Q), \frac{1}{\l_j}\frac{\partial \varphi_{a_j, \l_j}}{\partial a_j}>&=-\frac{8\pi^2}{\l_j}\frac{\n_g\mathcal{F}_j^{A}(a_j)}{\mathcal{F}_j^{A}(a_j)}\\&+O\left(\sum_{i=1}^m|\alpha_i-1|^2+\sum_{i=1}^m\frac{1}{\l_i^2}+\sum_{r=1}^{\bar m}|\beta_r|^2+\sum_{i=1}^m\tau_i^2\right).
\end{split}
\end{equation*}
thereby ending the proof of the Lemma.
\end{pf}
\vspace{6pt}

\noindent
Finally, we are going to establish a gradient estimate for $J$ in the direction of the $\beta_r$'s. Precisely, we will prove:
\begin{lem}\label{eq:gradientbeta}(Gradient estimate with respect to $\bar\beta$)\\
Assuming that $\eta$ is a small positive real number with $0<2\eta<\varrho$ where $\varrho$ is as in \eqref{eq:cutoff}, and $0<\epsilon\leq \epsilon_0$ where $\epsilon_0$ is as in \eqref{eq:mini}, then for $a_i\in M$ concentration points,  $\alpha_i$ masses, $\l_i$ concentration parameters ($i=1,\cdots,m$), ad $\beta_r$ negativity parameters ($r=1, \cdots, \bar m$) satisfying \eqref{eq:afpara}, we have that for every $l=1, \cdots, \bar m$, there holds
\begin{equation*}
\begin{split}
<\n J(\sum_{i=1}^m\alpha_i\varphi_{a_i, \l_i}+\sum_{r=1}^{\bar m}\beta_r (v_r-\ov{(v_r)}_Q), v_l-\ov{(v_l)}_Q>=&2\mu_l\beta_l+O\left(\sum_{i=1}^m|\alpha_i-1|+\sum_{i=1}^m |\tau_i| +\sum_{i=1}^m\frac{1}{\l_i^2}\right),
\end{split}
\end{equation*}
where $O(1)$ is as in Lemma \ref{eq:energyest} and for  $i=1, \cdots, m$, \;$\tau_i$ is as in Lemma \ref{eq:gradientlambdaest}
\end{lem}

\begin{pf}
Arguing as in Lemma \ref{eq:gradientlambdaest}, we obtain
\begin{equation}\label{eq:der6}
<\n J(\sum_{i=1}^m\alpha_i\varphi_{a_i, \l_i}+\sum_{r=1}^{\bar m}\beta_r (v_r-\ov{(v_r)}_Q),  v_l-\ov{(v_l)}_Q>=2(E_{1, \beta_l}+E_{0, \beta_l}-16\pi^2 mE_{2, \beta_l}),
\end{equation}
where
\begin{equation}\label{eq:derb3}
E_{1, \beta_l}:=<P_g(\sum_{i=1}^m\alpha_i \varphi_{a_i,\l_i}), v_l>+16\pi^2m \ov{(v_l)}_Q,
\end{equation}
\begin{equation}\label{eq:derb4}
E_{0, \beta_l}:=<P_g(\sum_{r=1}^{\bar m} \beta_r) v_r, v_l>,
\end{equation}
and
\begin{equation}\label{eq:der5}
E_{2, \beta_l}:=\frac{\int_M e^{4\sum_{i=1}^m\alpha_i\varphi_{a_i, \l_i}+4\sum_{r=1}^{\bar m}\beta_r v_r(x)}v_l dV_g}{\int_{M}e^{4\sum_{i=1}^m\alpha_i\varphi_{a_i, \l_i}+4\sum_{r=1}^{\bar m}\beta_r v_r(z)}dV_g}.
\end{equation}
Now, we are going to estimate $E_{1, \beta_l}$, $E_{0, \beta_l}$ and $E_{2, \beta_l}$. We start with $E_{1, \beta_l}$. For this, we use again the same strategy as in Lemma \ref{eq:gradientlambdaest} to get
\begin{equation}\label{eq:derb9}
E_{1, \beta_l}=16\pi^2\sum_{i=1}^m(\alpha_i-1)(v_l- \ov{(v)}_Q)(a_i)+16\pi^2\sum_{i=1}^m v_l(a_i)+O\left(\sum_{i=1}^m\frac{1}{\l_i^2}\right).
\end{equation}
Next, we estimate $E_{0, \beta_l}$. In order to do that, we use \eqref{eq:derb4} and the fact that the $v_r$'s are $L^2$-orthonormal basis of $P_g$ to obtain
\begin{equation}\label{eq:derb10}
E_{0, \beta_l}=\mu_l\beta_l.
\end{equation}
Now, we are going to estimate $E_{2, \beta_l}$. For this,  we set $E_{3, \beta_l}:=\int_M e^{4\sum_{i=1}^m\alpha_i\varphi_{a_i, \l_i}+4\sum_{r=1}^{\bar m}\beta_r v_r(x)}v_ldV_g$ and obtain
\begin{equation}\label{eq:derb11}
E_{2, \beta_l}=\frac{E_{3, \beta_l}}{D}.
\end{equation}
On the other hand, arguing  again as in Lemma \ref{eq:gradientlambdaest}, we get
\begin{equation}\label{eq:derb12}
E_{3, \beta_l}=\sum_{i=1}^m\gamma_i(v_l(a_i)+O\left(\frac{1}{\l_i^2}\right)).
\end{equation}
Combining \eqref{eq:der6}, \eqref{eq:derb9}, \eqref{eq:derb11}, \eqref{eq:derb10}, \eqref{eq:derb12}, and the smoothness of the $v_l$'s we obtain
\begin{equation*}
\begin{split}
<\n J(\sum_{i=1}^m\alpha_i\varphi_{a_i, \l_i}+\sum_{r=1}^{\bar m}\beta_r (v_r-\ov{(v_r)}_Q),  v_l-\ov{(v_l)}_Q>=&2\mu_l\beta_l+O\left(\sum_{i=1}^m|\alpha_i-1|+\sum_{i=1}^m |\tau_i| +\sum_{i=1}^m\frac{1}{\l_i^2}\right),
\end{split}
\end{equation*}
as desired.
\end{pf}

\section{Finite dimensional reduction}\label{eq:findim}
In this section, we perform a finite dimensional reduction. Indeed, we are going to show that in order to understand the variational analysis of $J$ in $V(m, \epsilon, \eta)$ ,  one can roughly speaking reduce it self to the finite dimensional problem $J(\sum_{i=1}^m\alpha_i\varphi_{a_i, \l_i}+\sum_{l=1}^{\bar m}\beta_l (v_l-\ov{(v_l)}_Q)+\bar w(\bar\alpha, A, \bar \beta, \bar \l))$ where $\bar w \in C^1$. For this end, we start with the following Proposition.
\begin{pro}\label{eq:expansionJ1}
Assuming that $\eta$ is a small positive real number with $0<2\eta<\varrho$ where $\varrho$ is as in \eqref{eq:cutoff}, and $0<\epsilon\leq \epsilon_0$ where $\epsilon_0$ is as in \eqref{eq:mini} and $u=\ov{u}_Q+\sum_{i=1}^m\alpha_i\varphi_{a_i, \l_i}+\sum_{r=1}^{\bar m}\beta_r (v_r-\ov{(v_r)}_Q)+w\in V(m, \epsilon, \eta)$ with $w$, the concentration points $a_i$,  the masses $\alpha_i$, the concentrating parameters $\l_i$  ($i=1, \cdots, m$), and the negativity parameters $\beta_r$ ($r=1, \cdots, \bar m$) verifying \eqref{eq:ortho}-\eqref{eq:afpara}, then we have
\begin{equation}\label{eq:exparoundbubble}
J(u)=J(\sum_{i=1}^m\alpha_i\varphi_{a_i, \l_i}+\sum_{r=1}^{\bar m}\beta_r (v_r-\ov{(v_r)}_Q)-f(w)+Q(w)+o(||w||^2),
\end{equation}
where
\begin{equation}\label{eq:linear}
f(w):=8m\pi^2\frac{\int_M Ke^{4\sum_{i=1}^m\alpha_i\varphi_{a_i, \l_i}+4\sum_{r=1}^{\bar m}\beta_r v_r}wdV_g}{\int_M Ke^{4\sum_{i=1}^m\alpha_i\varphi_{a_i, \l_i}+4\sum_{r=1}^{\bar m}\beta_r v_r}dV_g},
\end{equation}
and
\begin{equation}\label{eq:quadratic}
Q(w):=||w||^2-64m\pi^2\frac{\int_M Ke^{4\sum_{i=1}^m\alpha_i\varphi_{a_i, \l_i}+4\sum_{r=1}^{\bar m}\beta_r v_r}w^2dV_g}{\int_M Ke^{4\sum_{i=1}^m\alpha_i\varphi_{a_i, \l_i}+4\sum_{r=1}^{\bar m}\beta_r v_r}dV_g}.
\end{equation}
Moreover, setting
\begin{equation}\label{eq:eali}
\begin{split}
E_{a_i, \l_i}:=\{w\in W^{2, 2}(M): \;\;<\varphi_{a_i, \l_i}, w>_P=<\frac{\partial\varphi_{a_i, \l_i}}{\partial \l_i}, w>_P=<\frac{\partial\varphi_{a_i, \l_i}}{\partial a_i}, w>_P=0,\\\,\;\,\;<w, Q_g>=<v_k, w>=0, \;k=1, \cdots, \bar m,\,\;\text{and}\;\;||w||=O\left(\sum_{i=1}^m\frac{1}{\l_i^{1-\gamma}}\right)\},
\end{split}
\end{equation}
and
\begin{equation}\label{eq:eal}
A:=(a_1, \cdots, a_m), \;\;\bar \l=(\l_1, \cdots, \l_m), \;\;E_{A, \bar \l}:=\cap_{i=1}^m E_{a_i, \l_i},
\end{equation}
we have that, the quadratic form $Q$ is positive definite in $E_{A, \bar \l}$. Furthermore, the linear par $f$ verifies that, for every $w\in E_{A, \bar \l}$, there holds
\begin{equation}\label{eq:estlinear}
f(w)=O\left( ||w||\left(\sum_{i=1}^m\frac{|\n_g \mathcal{F}^{A}_i(a_i)|}{\l_i}+\sum_{i=1}^m|\alpha_i-1|\log \l_i+\sum_{r=1}^{\bar m}|\beta_r|+\sum_{i=1}^m\frac{\log \l_i}{\l_i^{2}}\right)\right).
\end{equation}
where here $o(1)=o_{ \bar\alpha, A, \bar \beta, \bar\l, w, \epsilon}(1)$ and $O\left(1\right):=O_{\bar\alpha, A, \bar \beta, \bar\l, w, \epsilon}\left(1\right)$ and for their meaning see section \ref{eq:notpre}.
\end{pro}
\vspace{4pt}

\noindent
To prove Proposition \ref{eq:expansionJ1}, we will need the following three coming Lemmas. We start with the following one.
\begin{lem}\label{eq:holmos}
Assuming the assumptions of Proposition \ref{eq:expansionJ1}, then for every $q\geq 1$, there holds the following estimates\\
\begin{equation}\label{eq:esthm1}
\int_M Ke^{4\sum_{i=1}^m\alpha_i\varphi_{a_i, \l_i}+4\sum_{r=1}^{\bar m}\beta_r v_r}|w|^q=O\left(||w||^q(\sum_{i=1}^m\l_i^{4+\gamma})\right),
\end{equation}
\begin{equation}\label{eq:esthm2}
\int_M Ke^{4\sum_{i=1}^m\alpha_i\hat\d_{a_i, \l_i}+4\sum_{r=1}^{\bar m}\beta_r v_r}|w|^q=O\left(||w||^q(\sum_{i=1}^m\l_i^\gamma)\right),
\end{equation}
\begin{equation}\label{eq:esthm3}
\int_M e^{4\hat\d_{a_i, \l_i}+4\sum_{r=1}^{\bar m}\beta_r v_r}d_{g_{a_i}}(a_i, \cdot)|w|^q=O\left(||w||^q\frac{1}{\l_i^{1-\gamma}}\right),\;\; i=1, \cdots, m,
\end{equation}
\begin{equation}\label{eq:esthm4}
\int_M Ke^{4\sum_{i=1}^m\alpha_i \varphi_{a_i, \l_i}+4\sum_{r=1}^{\bar m}\beta_r v_r}e^{4\theta_w w}|w|^q=O\left(||w||^q(\sum_{i=1}^m\l_i^{4+\gamma})\right),
\end{equation}
where $\theta_w\in [0, 1]$, and
\begin{equation}\label{eq:esthm5}
\int_M Ke^{4\sum_{i=1}^m\alpha_i\varphi_{a_i, \l_i}+4\sum_{r=1}^{\bar m}\beta_r v_r}\left(e^{4w}-1-4w-8w^2\right)dV_g=o\left(||w||^2(\sum_{i=1}^m\l_i^{4})\right).
\end{equation}
where here $o(1)$ and $O\left(1\right)$ are as in Proposition \ref{eq:expansionJ1}.
\end{lem}
\begin{pf}
First of all,  it is easy to see that  \eqref{eq:esthm1}-\eqref{eq:esthm3} follow from Lemma \ref{eq:intbubbleest} and the explicit expression of $\hat\d_{a_i, \l_i}$ combined with a trivial application H\"older inequality and the continuous embedding of $W^{2, 2}(M)$ in any  $L^p(M)$ with $1\leq p<+\infty$. Furthermore, \eqref{eq:esthm4} follows from the same reasons a as above combined with a direct application of Moser-Trudinger inequality. Finally \eqref{eq:esthm5} follows from Taylor expansion, the following estimate $|e^{4w}-1-4w-8w^2|\leq C_{\nu}e^{4|w|}|w|^3$ for $|w|\geq \nu$ with $\nu>0$ and very small, \eqref{eq:estwmin}, \eqref{eq:afpara}, and the arguments of \eqref{eq:esthm4}.
\end{pf}
\vspace{4pt}

\noindent
The second Lemma that we need for the proof of Proposition \ref{eq:expansionJ1} read as follows:
\begin{lem}\label{eq:holmos1}
Assuming the assumptions of Proposition \ref{eq:expansionJ1}, then there holds the following estimate
\begin{equation}\label{eq:esthm6}
\begin{split}
&\frac{\int_M K e^{4\sum_{i=1}^m\alpha_i\varphi_{a_i, \l_i}+\sum_{r=1}^{\bar m}\beta_r v_r}wdV_g}{\int_M K e^{4\sum_{i=1}^m\alpha_i\varphi_{a_i, \l_i}+4\sum_{r=1}^{\bar m}\beta_r v_r}dV_g}=\\&O\left(||w||\left(\sum_{i=1}^m\frac{|\n_g \mathcal{F}^{A}_i(a_i)|}{\l_i}+\sum_{i=1}^m|\alpha_i-1|\log\l_i+\sum_{r=1}^{\bar m}|\beta_r|+\sum_{i=1}^m\frac{\log \l_i}{\l_i^{2}}\right)\right).
\end{split}
\end{equation}
\end{lem}

\begin{pf}
Using Lemma \ref{eq:intbubbleest}, the fact that $\sum_{r=1}^{\bar m}|\beta_r|+\sum_{i=1}^m|\alpha_i-1|\sqrt{\log \l_i}=O\left(\sum_{i=1}^m\frac{1}{\l_i^{1-\gamma}}\right)$ and Sobolev embedding theorem, we infer that
\begin{equation}\label{eq:linear1}
\int_M  K e^{4\sum_{j=1}^m\alpha_j\varphi_{a_j, \l_j}+4\sum_{r=1}^{\bar m}\beta_r v_r}wdV_g=\sum_{i=1}^m\int_{B_{a_i}^{a_i}(\eta)} K e^{4\sum_{j=1}^m\alpha_j\varphi_{a_j, \l_j}+4\sum_{k=1}^{\bar m}\beta_k v_k}wdV_g +O(||w||).
\end{equation}
On the other hand, using again Lemma \ref{eq:intbubbleest}, and the fact that $\sum_{k=1}^{\bar m}|\beta_k|+\sum_{i=1}^m|\alpha_i-1|\sqrt{\log \l_i}=O\left(\sum_{i=1}^m\frac{1}{\l_i^{1-\gamma}}\right)$, the explicit expression of $\hat \d_{a_i, \l_i}$,  and arguing as in Lemma \ref{eq:energyest}, we have  in $B_{a_i}^{a_i}(\eta)$ that the following estimate holds
\begin{equation}\label{eq:linear2}
\begin{split}
& K e^{4\sum_{i=1}^m\alpha_j\varphi_{a_j, \l_j}+\sum_{k=1}^{\bar m}\beta_k v_k}=\left(\frac{\l_i^{2}}{1+\l_i^2d^2_{g_{a_i}}(\cdot,  a_i)}\right)^{4\alpha_i}\mathcal{F}_i^{A}\mathcal{G}_i\left(1+O\left(\sum_{i=1}^m\frac{1}{\l_i^3}\right)\right)\\&=\frac{\l_i^4}{16}\left(\frac{\l_i^{2}}{1+\l_i^2d^2_{g_{a_i}}(\cdot,  a_i)}\right)^{4\alpha_i-4}e^{4\hat\d_{a_i, \l_i}}\mathcal{F}_i^{A}\left(1+O\left(\sum_{j=1}^m|\alpha_j-1|+\sum_{r=1}^{\bar m}|\beta_r|+\sum_{j=1}^m\frac{1}{\l_j^2}\right)\right)\\&=\frac{\l_i^4}{16}e^{4\hat\d_{a_i, \l_i}}\mathcal{F}_i^{A}\left(1+O\left(\sum_{j=1}^m|\alpha_j-1|\log\l_j+\sum_{r=1}^{\bar m}|\beta_r|+\sum_{j=1}^m\frac{1}{\l_j^2}\right)\right)
 \end{split}
\end{equation}
Hence, combining \eqref{eq:linear1} and \eqref{eq:linear2}, we get
\begin{equation}\label{eq:linear3}
\begin{split}
&\int_M  K e^{4\sum_{j=1}^m\alpha_j\varphi_{a_j, \l_j}+\sum_{r=1}^{\bar m}\beta_r v_r}wdV_g=\frac{1}{16}\sum_{i=1}^m\l_i^4\int_{B_{a_i}^{a_i}(\eta)}e^{4\hat\d_{a_i, \l_i}}\mathcal{F}_i^{A}wdV_g\\&+O\left(\sum_{j=1}^m|\alpha_j-1|\log\l_j+\sum_{r=1}^{\bar m}|\beta_r|+\sum_{j=1}^m\frac{1}{\l_j^2}\right)\left(\sum_{i=1}^m\l_i^4\int_{B_{a_i}^{a_i}(\eta)}e^{4\hat\d_{a_i, \l_i}}|w|dV_g\right)+O(||w||).
\end{split}
\end{equation}
Now, using Taylor expansion at $a_i$ and Lemma \ref{eq:positive} combined with H\"older inequality, we deduce that
\begin{equation}\label{eq:linear4}
\int_{B_{a_i}^{a_i}(\eta)}e^{4\hat\d_{a_i, \l_i}}\mathcal{F}_i^{A}wdV_g=\mathcal{F}_i^{A}(a_i)\int_{B_{a_i}^{a_i}(\eta)}e^{4\hat\d_{a_i, \l_i}}wdV_{g_{a_i}}+O\left(\frac{|\n_g\mathcal{F}_i^{A}(a_i)|}{\l_i}||w||\right)+O\left(\frac{\log \l_i}{\l_i^2}||w||\right).
\end{equation}
On the other hand, since $w\in E_{A, \bar\l}$, then we deduce that
\begin{equation}\label{eq:cortho}
\int_{M}e^{4\hat\d_{a_i, \l_i}}wdV_{g_{a_i}}=0.
\end{equation}
Hence, we have that \eqref{eq:cortho} implies that
\begin{equation}\label{eq:cortho1}
\int_{B_{a_i}^{a_i}(\eta)}e^{4\hat\d_{a_i, \l_i}}wdV_{g_{a_i}}=-\int_{M\setminus B_{a_i}^{a_i}(\eta)}e^{4\hat\d_{a_i, \l_i}}wdV_{g_{a_i}}=O\left(\frac{1}{\l_i^4}||w||\right).
\end{equation}
On the other hand, using H\"older inequality and Lemma \ref{eq:positive}, we get
\begin{equation}\label{eq:cortho2}
\int_{B_{a_i}^{a_i}(\eta)}e^{4\hat\d_{a_i, \l_i}}|w|V_g=O\left(||w||\right).
\end{equation}
Now, combining \eqref{eq:linear3}, \eqref{eq:linear4}, \eqref{eq:cortho1}, and \eqref{eq:cortho2}, we obtain
\begin{equation}\label{eq:linear5}
\begin{split}
&\int_M  K e^{4\sum_{j=1}^m\alpha_j\varphi_{a_j, \l_j}+\sum_{r=1}^{\bar m}\beta_r v_r}wdV_g=\\&O\left(||w||(\sum_{i=1}^m\l_i^{4})\left(\sum_{i=1}^m\frac{|\n_g \mathcal{F}^{A}_i(a_i)|}{\l_i}+\sum_{i=1}^m|\alpha_i-1|\log\l_i+\sum_{r=1}^{\bar m}|\beta_r|+\sum_{i=1}^m\frac{\log \l_i}{\l_i^{2}}\right)\right).
\end{split}
\end{equation}
On the other hand, using \eqref{eq:5imp31a}, \eqref{eq:5im38}, \eqref{eq:6imp6}, \eqref{eq:coreql3} and recalling that the $\l_i$'s are comparable, we derive
\begin{equation}\label{eq:linear6}
||K e^{4\sum_{j=1}^m\alpha_j\varphi_{a_j, \l_j}+4\sum_{r=1}^{\bar m}\beta_r v_r}||_{L^1(M)}\geq C^{-1}\sum_{i=1}^m\l_i^{4},
\end{equation}
where $C$ is a large positive constant independent of $\bar \alpha$, $\bar \l$, $A$, $\bar \beta$ and of course of $w$. Hence, we have that the Lemma follows from \eqref{eq:linear5} and \eqref{eq:linear6}.
\end{pf}
\vspace{4pt}

\noindent
Finally, the third and last Lemma that we need for the proof of Proposition \ref{eq:expansionJ1} is the following one.
\begin{lem}\label{eq:tauiest}
Assuming the assumptions of Proposition \ref{eq:expansionJ1}, then for every $i=1, \cdots, m$,  there holds
\begin{equation}\label{eq:esttauieq}
\tau_i=O\left(\sum_{j=1}^m\frac{1}{\l_j^{1-2\gamma}}\right).
\end{equation}
\end{lem}
\begin{pf}
First of all, using Lemma \ref{eq:gradientlambdaest}, we infer that
\begin{equation}\label{eq:tauiest1}
32\pi^2\alpha_i\tau_i=<\n^{W^{2, 2}} J(\sum_{i=1}^m\alpha_j\varphi_{a_j, \l_j}), \l_i\frac{\partial \varphi_{a_i, \l_i}}{\partial \l_i}>_{W^{2, 2}(M)}+O\left(\sum_{j=1}^m\frac{1}{\l_j^{1-\gamma}}\right).
\end{equation}
Now, using \eqref{eq:tauiest1}, we derive that
\begin{equation}\label{eq:tauiest2}
\tau_i=O\left(||\n^{W^{2, 2}} J(\sum_{i=1}^m\alpha_i\varphi_{a_i, \l_i})||_{W^{2, 2}(M)}||\l_i\frac{\partial \varphi_{a_i, \l_i}}{\partial \l_i}||_{W^{2, 2}(M)}+\sum_{j=1}^m\frac{1}{\l_j^{1-\gamma}}\right)
\end{equation}
On the other hand, recalling that $u\in V(m, \epsilon, \eta)$ ( hence  $||\n^{W^{2, 2}} J(u)||_{W^{2, 2}(M)}=O\left(\sum_{j=1}^m\frac{1}{\l_j^{1-\gamma}}\right)$ and using the same arguments as in the proof of \eqref{eq:diff18} combined with the fact that $||\l_i\frac{\partial \varphi_{a_i, \l_i}}{\partial \l_i}||_{W^{2, 2}(M)}=O(1)$ (see \eqref{eq:normpr} and the discussion along it and Lemma \ref{eq:auxibubbleest1}), we infer that \eqref{eq:tauiest2} implies
\begin{equation}\label{eq:tauiest3}
\tau_i=O\left(\sum_{j=1}^m\frac{1}{\l_j^{1-2\gamma}}\right),
\end{equation}
thereby ending the proof of the Lemma.
\end{pf}
\vspace{6pt}

\noindent
Now, we are ready to give the proof of Proposition \ref{eq:expansionJ1}.\\\\
\begin{pfn} {of Proposition \ref{eq:expansionJ1}}\\
First of all, setting $v=\sum_{i=1}^m\alpha_i\varphi_{a_i, \l_i}+\sum_{r=1}^{\bar m}\beta_r(v_r-\ov{(v_r)}_Q)$  and using the fact that $J$ is invariant by translation by constant , $w\in E_{A, \bar \l}$, $\ov{v}_Q=0$ and expanding $J$, we get
\begin{equation}\label{eq:expaja1}
\begin{split}
&J(u)=J(u-\ov{u}_Q)=J(v)+||w||^2\\&-8m\pi^2\log\left(1+4\frac{\int_M Ke^{4v}wdV_g}{\int_M Ke^{4v}dV_g}+8\frac{\int_M Ke^{4v}w^2dV_g}{\int_M Ke^{4v}dV_g}+\frac{\int_M Ke^{4v}(e^{4w}-1-4w-8w^2)dV_g}{\int_M Ke^{4v}dV_g}\right).
\end{split}
\end{equation}
Thus, appealing to Lemma \ref{eq:holmos} , Lemma \ref{eq:holmos1} and \eqref{eq:expaja1}, we infer that
\begin{equation}\label{eq:expansionja}
J(\sum_{i=1}^m\alpha_i\varphi_{a_i, \l_i}+\sum_{r=1}^{\bar m}\beta_r(v_r-\ov{(v_r)}_Q)+w)=J(\sum_{i=1}^m\alpha_i\varphi_{a_i, \l_i}+\sum_{r=1}^{\bar m}\beta_r(v_r-\ov{(v_r)}_Q))-f(w)+Q(w)+o\left(||w||^2\right),
\end{equation}
and
\begin{equation}\label{eq:estfw}
f(w)=O\left(||w||\left(\sum_{i=1}^m\frac{|\n_g \mathcal{F}^{A}_i(a_i)|}{\l_i}+\sum_{i=1}^m|\alpha_i-1|\log\l_i+\sum_{r=1}^{\bar m}|\beta_r|+\sum_{i=1}^m\frac{\log \l_i}{\l_i^{2}}\right)\right).
\end{equation}
Now, we are going to show that $Q$ is positive definite in $E_{A, \bar \l}$.. For this end, we first use Lemma \ref{eq:intbubbleest} to obtain
\begin{equation}\label{eq:q1}
\int_M K e^{4\sum_{j=1}^m\alpha_j\varphi_{a_j, \l_j}+4\sum_{r=1}^{\bar m}\beta_r v_r}w^2dV_g=\sum_{i=1}^m\int_{B_{a_i}^{a_i}(\eta)} K e^{4\sum_{j=1}^m\alpha_j\varphi_{a_j, \l_j}+4\sum_{r=1}^{\bar m}\beta_r v_r}w^2+ O\left(||w||^2\right).
\end{equation}
Next, we recall that on each $B_{a_i}^{a_i}(\eta)$ the following estimate holds (see \eqref{eq:linear2})
\begin{equation}\label{eq:q2}
 K e^{4\sum_{j=1}^m\alpha_j\varphi_{a_j, \l_j}+4\sum_{r=1}^{\bar m}\beta_r v_r}=\frac{\l_i^4}{16} e^{4\hat\d_{a_i, \l_i}}\mathcal{F}^{A}_i\left(1+O\left(\sum_{j=1}^m|\alpha_j-1|\log\l_j+\sum_{r=1}^{\bar m}|\beta_r|+\sum_{j=1}^m\frac{1}{\l_j^2}\right)\right).
\end{equation}
Thus, appealing to \eqref{eq:q2}, we obtain
\begin{equation}\label{eq:q3}
\begin{split}
&\int_{B_{a_i}^{a_i}(\eta)}K e^{4\sum_{j=1}^m\alpha_j\varphi_{a_j, \l_j}+4\sum_{r=1}^{\bar m}\beta_r v_r}w^2dV_g=\frac{\l_i^4}{16} \int_{B_{a_i}^{a_i}(\eta)}e^{4\hat\d_{a_i, \l_i}}\mathcal{F}^{A}_iw^2dV_g\\&+O\left(\left(\sum_{j=1}^m|\alpha_j-1|\log\l_j+\sum_{r=1}^{\bar m}|\beta_r|+\sum_{i=1}^m\frac{1}{\l_j^2}\right)\l_i^4\int_{B_{a_i}^{a_i}(\eta)}e^{4\hat\d_{a_i, \l_i}}w^2\right).
\end{split}
\end{equation}
Now, making Taylor expansion at $a_i$ and using again \eqref{eq:esthm3}, we get
\begin{equation}\label{eq:q4}
\frac{\l_i^4}{16} \int_{B_{a_i}^{a_i}(\eta)}e^{4\hat\d_{a_i, \l_i}}\mathcal{F}^{A}_iw^2dV_g=\frac{\l_i^4}{16} \mathcal{F}^{A}_i(a_i)\int_{B_{a_i}^{a_i}(\eta)}e^{4\hat\d_{a_i, \l_i}}w^2dV_{g_{a_i}}+o\left(\l_i^4||w||^2\right).
\end{equation}
Thus, using again Lemma \ref{eq:positive} combined with \eqref{eq:q1}, \eqref{eq:q3}, and  \eqref{eq:q4},  we derive that the following estimate holds
\begin{equation}\label{eq:q5}
\int_M K e^{4\sum_{j=1}^m\alpha_j\varphi_{a_j, \l_j}+4\sum_{r=1}^{\bar m}\beta_r v_r}w^2dV_g=\sum_{i=1}^m\frac{\l_i^4}{16} \mathcal{F}^{A}_i(a_i)\int_{M}e^{4\hat\d_{a_i, \l_i}}w^2dV_{g_{a_i}}+o\left((\sum_{i=1}^m\l_i^4)||w||^2\right).
\end{equation}
Therefore, combining \eqref{eq:linear1} and \eqref{eq:q5}, we obtain
\begin{equation}\label{eq:q6}
\begin{split}
&64m\pi^2\frac{\int_M K e^{4\sum_{j=1}^m\alpha_j\varphi_{a_j, \l_j}+4\sum_{r=1}^{\bar m}\beta_r v_r}w^2dV_g}{\int_M K e^{4\sum_{j=1}^m\alpha_j\varphi_{a_j, \l_j}+4\sum_{r=1}^{\bar m}\beta_r v_r}dV_g}=\\&8m\pi^2\sum_{i=1}^m\frac{\l_i^4}{\int_M K e^{4\sum_{j=1}^m\alpha_j\varphi_{a_j, \l_j}+4\sum_{r=1}^{\bar m}\beta_r v_r}dV_g} \mathcal{F}^{A}_i(a_i)\int_{M}e^{4\hat\d_{a_i, \l_i}}w^2dV_{g_{a_i}}+o\left(||w||^2\right).
\end{split}
\end{equation}
Now, using \eqref{eq:tildegammai}, we derive
\begin{equation}\label{eq:critgammai}
\gamma_i==\frac{\pi^2}{6}\l_i^4\mathcal{F}_i^{A}(a_i)(1+O\left(\sum_{j=1}^m|\alpha_j-1|\log\l_j+\sum_{k=1}^{\bar m}|\beta_k|+\sum_{j=1}^m\frac{1}{\l_j^2}\right))
\end{equation}
Thus, using again Lemma \ref{eq:positive} combined with \eqref{eq:linear6} and \eqref{eq:critgammai}, we have that \eqref{eq:q6} implies
\begin{equation}\label{eq:q7}
\begin{split}
&64m\pi^2\frac{\int_M K e^{4\sum_{j=1}^m\alpha_j\varphi_{a_j, \l_j}+4\sum_{r=1}^{\bar m}\beta_r v_r}w^2dV_g}{\int_M K e^{4\sum_{j=1}^m\alpha_j\varphi_{a_j, \l_j}+4\sum_{r=1}^{\bar m}\beta_r v_r}dV_g}=\\&48m\sum_{i=1}^m\frac{\gamma_i}{\int_M K e^{4\sum_{j=1}^m\alpha_j\varphi_{a_j, \l_j}+4\sum_{r=1}^{\bar m}\beta_r v_r}dV_g} \int_{M}e^{4\hat\d_{a_i, \l_i}}w^2dV_{g_{a_i}}+o\left(||w||^2\right).
\end{split}
\end{equation}
On the other hand, using the definition of $\tau_i$ (see \eqref{eq:6imp53}), we infer from \eqref{eq:q7} that
\begin{equation}\label{eq:q8}
\begin{split}
&64m\pi^2\frac{\int_M K e^{4\sum_{j=1}^m\alpha_j\varphi_{a_j, \l_j}+4\sum_{r=1}^{\bar m}\beta_r v_r}w^2dV_g}{\int_M K e^{4\sum_{j=1}^m\alpha_j\varphi_{a_j, \l_j}+4\sum_{r=1}^{\bar m}\beta_r v_r}dV_g}=\\&48\sum_{i=1}^m \int_{M}e^{4\hat\d_{a_i, \l_i}}w^2dV_g+48\sum_{i=1}^m\tau_i \int_{M}e^{4\hat\d_{a_i, \l_i}}w^2dV_{g_{a_i}}+o\left(||w||^2\right).
\end{split}
\end{equation}
Thus, using again  Lemma \ref{eq:positive} combined with Lemma \ref{eq:tauiest}, we have that \eqref{eq:q8} implies
\begin{equation}\label{eq:q9}
\begin{split}
64m\pi^2\frac{\int_M K e^{4\sum_{j=1}^m\alpha_j\varphi_{a_j, \l_j}+4\sum_{r=1}^{\bar m}\beta_r v_r}w^2dV_g}{\int_M K e^{4\sum_{j=1}^m\alpha_j\varphi_{a_j, \l_j}+4\sum_{r=1}^{\bar m}\beta_r v_r}dV_g}=48\sum_{i=1}^m \int_{M}e^{4\hat\d_{a_i, \l_i}}w^2dV_{g_{a_i}}+o\left(||w||^2\right).
\end{split}
\end{equation}
Now, using the definition of $Q$ (see \eqref{eq:quadratic}), we have that \eqref{eq:q9} gives
\begin{equation}\label{eq:q10}
\begin{split}
Q(w)=||w||^2-48\sum_{i=1}^m \int_{M}e^{4\hat\d_{a_i, \l_i}}w^2dV_{g_{a_i}}+o\left(||w||^2\right).
\end{split}
\end{equation}
Therefore, using \eqref{eq:q10}, Lemma \ref{eq:positiveg}, and the obvious fact that $E_{A, \bar\l}\subset E^{*}_{A, \bar \l}$ (where $E^{*}_{A, \bar \l}$ is as in Lemma \ref{eq:positiveg}), we infer that
\begin{equation}\label{eq:Qpositive}
Q\,\;\;\text{is positive definite on }\;\;E_{A, \l}.
\end{equation}
Hence the Lemma follows from \eqref{eq:expansionja}, \eqref{eq:estfw} and \eqref{eq:Qpositive}.
\end{pfn}
\vspace{6pt}

\noindent
Like in \cite{bcch}, we have that Proposition \ref{eq:expansionJ1} implies the following direct corollaries.
\begin{cor}\label{eq:cexpansionj1}
Assuming that $\eta$ is a small positive real number with $0<2\eta<\varrho$ where $\varrho$ is as in \eqref{eq:cutoff}, $0<\epsilon\leq \epsilon_0$ where $\epsilon_0$ is as in \eqref{eq:mini} and $u:=\sum_{i=1}^m\alpha_i\varphi_{a_i, \l_i}+\sum_{r=1}^{\bar m}\beta_r(v_r-\ov{(v_r)}_Q$ with the concentration  points $a_i$,  the masses $\alpha_i$, the concentrating parameters $\l_i$ ($i=1, \cdots, m$)  and the negativity parameters $\beta_r$ ($r=1, \cdots, \bar m$) satisfying \eqref{eq:afpara}, then there exists a unique $\bar w(\bar \alpha, A, \bar \l, \bar \beta)\in E_{A, \bar \l}$ such that
\begin{equation}\label{eq:minj}
J(u+\bar w(\bar \alpha, A, \bar \l, \bar \beta))=\min_{w\in E_{A, \bar \l}, u+w\in V(m, \epsilon, \eta)} J(u+w),
\end{equation}
where $\bar \alpha:=(\alpha_1, \cdots, \alpha_m)$, $A:=(a_1, \cdots, a_m)$, $\bar \l:=(\l_1, \cdots, \l_m)$ and $\bar \beta:=(\beta_1, \cdots, \beta_m)$.\\
Furthermore, $(\bar \alpha, A, \bar \l, \bar \beta)\longrightarrow \bar w(\bar \alpha, A, \bar \l, \bar \beta)\in C^1$ and satisfies the following estimate
\begin{equation}\label{eq:linminqua}
\frac{1}{C}||\bar w(\bar \alpha, A,  \bar \l, \bar \beta)||^2\leq |f(\bar w(\bar \alpha, A, \bar \l, \bar \beta))|\leq C|| \bar w(\bar \alpha, A, \bar \l, \bar \beta)||^2,
\end{equation}
for some large positive constant $C$ independent of $\bar \alpha$, $A$, $\bar \l$, and $\bar \beta$, hence
\begin{equation}\label{eq:estbarw}
||\bar w(\bar \alpha, A, \bar \l, \bar \beta)||=O\left(\sum_{i=1}^m\frac{|\n_g \mathcal{F}^{A}_i(a_i)|}{\l_i}+\sum_{i=1}^m|\alpha_i-1|\log\l_i+\sum_{r=1}^{\bar m}|\beta_r|+\sum_{i=1}^m\frac{\log \l_i}{\l_i^{2}}\right) .
\end{equation}
\end{cor}
\vspace{4pt}

\noindent
\begin{cor}\label{eq:c1expansionj1}
Assuming that $\eta$ is a small positive real number with $0<2\eta<\varrho$ where $\varrho$ is as in \eqref{eq:cutoff}, $0<\epsilon\leq \epsilon_0$ where $\epsilon_0$ is as in \eqref{eq:mini}, and $u_0:=\sum_{i=1}^m\alpha_i^0\varphi_{a_i^0, \l_i^0}+\sum_{r=1}^{\bar m}\beta_k^0(v_r-\ov{(v_r)}_Q)$ with the concentration  points $a_i^0$,  the masses $\alpha_i^0$, the concentrating parameters $\l_i^0$ ($i=1, \cdots, m$) and the negativity parameters $\beta_r^0$ ($r=1, \cdots, \bar m$) satisfying \eqref{eq:afpara}, then there exists an open neighborhood  $U$ of $(\bar\alpha^0, A^0, \bar \l^0, \bar \beta^0)$  (with $\bar \alpha^0:=(\alpha^0_1, \cdots, \alpha^0_m)$, $A^0:=(a_1^0, \cdots, a^0_m)$,  $\bar \l:=(\l_1^0, \cdots, \l_m^0)$ and  $\bar \beta^0:=(\beta_1^0, \cdots, \beta_{\bar m}^0)$) such that for every $(\bar \alpha, A, \bar \l, \bar \beta)\in U$ with $\bar \alpha:=(\alpha_1, \cdots, \alpha_m)$, $A:=(a_1, \cdots, a_m)$,  $\bar \l:=(\l_1, \cdots, \l_m)$, $\bar \beta:=(\beta_1, \cdots, \beta_{\bar m})$, and the $a_i$,  the $\alpha_i$, the $\l_i$ ($i=1, \cdots, m$)  and the $\beta_r$ ($r=1, \cdots, \bar m$) satisfying \eqref{eq:afpara}, and $w$ satisfying \eqref{eq:afpara}  with $\sum_{i=1}^m\alpha_i\varphi_{a_i, \l_i}+\sum_{r=1}^{\bar m}\beta_r(v_r-\ov{(v_r)}_Q+w\in V(m, \epsilon, \eta)$, we have the existence of a change of variable
\begin{equation}\label{eq:changev}
w\longrightarrow V
\end{equation}
from a neighborhood of $ \bar w(\bar \alpha, A, \bar \l, \bar \beta)$ to a neighborhood of $0$ such that
\begin{equation}\label{eq:expjv}
\begin{split}
&J(\sum_{i=1}^m\alpha_i\varphi_{a_i, \l_i}+\sum_{r=1}^{\bar m}\beta_r(v_r-\ov{(v_r)}_Q)+w)=\\&J(\sum_{i=1}^m\alpha_i\varphi_{a_i, \l_i}+\sum_{r=1}^{\bar m}\beta_r(v_r-\ov{(v_r)}_Q)+\bar w (\bar \alpha, A, \bar \l, \bar \beta))+\frac{1}{2}\partial^2 J(\sum_{i=1}^m\alpha_i^0\varphi_{a_i^0, \l_i^0}+\bar w(\bar \alpha^0, A^0, \bar \l^0, \bar \beta^0))(V, V),
\end{split}
\end{equation}
\end{cor}
\vspace{4pt}

\noindent
Thus, with this new variable, it is easy to see that in $V(m, \epsilon, \eta)$ we have splitting of the variable $(\bar \alpha, A, \bar \l, \bar \beta)$ and $V$, namely that one can decrease the functional $J$ in the variable $V$  without touching the variable $(\bar \alpha, A, \bar \l, \bar \beta)$  by considering just the flow
\begin{equation}\label{eq:vflow}
\frac{dV}{dt}=-V.
\end{equation}
So, since $J$ is invariant by translations by constants, then the variational study of $J$ in $V(m, \epsilon, \eta)$ (in the spirit of Morse theory) is equivalent to the one of  \begin{equation}\label{eq:finitedj}
J(\bar \alpha, A, \bar \l, \bar \beta):=J(\sum_{i=1}^m\alpha_i\varphi_{a_i, \l_i}+\sum_{r=1}^{\bar m}\beta_r(v_r-\ov{(v_r)}_Q)+\bar w(\bar \alpha, A, \bar \l, \bar \beta)),
\end{equation}
where $\bar \alpha=(\alpha_1, \cdots, \alpha_m)$, $A=(a_1, \cdots, a_m)$, $\bar \l=(\l_1, \cdots, \l_m)$ and $\bar \beta=\beta_1, \cdots, \beta_{\bar m}$ with the concentration  points $a_i$,  the masses $\alpha_i$, the concentrating parameters $\l_i$ ($i=1, \cdots, m$)  and the negativity parameters $\beta_r$ ($r=1, \cdots, \bar m$) satisfying \eqref{eq:afpara}, and $\bar w(\bar \alpha, A, \bar \l, \bar \beta)$ is as in Corollary \ref{eq:cexpansionj1}. Hence, the goal of this section is achieved.

\section{Topology of some appropriate sublevels of $J$}\label{eq:topsublevel}
In this section, we characterize the topology  of very high sublevels of $J$ and its every negative ones. In order to do that, we start with the following improvement of Lemma \ref{eq:energyest}.
\begin{lem}\label{eq:expansionj}
Under  the assumptions of Proposition \ref{eq:expansionJ1}, there holds
\begin{equation*}
\begin{split}
&J(\sum_{i=1}^m\alpha_i\varphi_{a_i, \l_i}+\sum_{r=1}^{\bar m}\beta_r(v_r-\ov{(v_r)}_Q)+w)=-\frac{40}{3}m\pi^2-8m\pi^2\log(\frac{m\pi^2}{6})-8\pi^2\mathcal{F}_K(a_1, \dots, a_m)\\&+\sum_{i=1}^m(\alpha_i-1)^2\left(32\pi^2\log\l_i+16\pi^2 H(a_i, a_i)+\frac{472\pi^2}{9}\right)+\sum_{r=1}^{\bar m}\mu_r\beta_r^2\\&+16\pi^2\sum_{i=1}^m(\alpha_i-1)\left[\sum_{r=1}^{\bar m}2\beta_r(v_r-\ov{(v_r)}_Q)(a_i)-\sum_{ j=1, j\neq i}^m(\alpha_j-1) G(a_i, a_j)\right]\\&-2\pi^2\sum_{i=1}^m\frac{1}{\l_i^2}\left(\frac{\D_{g_{a_i}}\mathcal{F}^{A}_i(a_i)}{\mathcal{F}^{A}_i(a_i)}-\frac{2}{3}R_g(a_i)\right)+2\pi^2\sum_{i=1}^m\frac{\tilde \tau_i}{\l_i^2}\left(\frac{\D_{g_{a_i}}\mathcal{F}^{A}_i(a_i)}{\mathcal{F}^{A}_i(a_i)}-\frac{2}{3}R_g(a_i)\right)\\&+8\pi^2\sum_{i=1}^m\log(1-\tilde\tau_i)+O\left(\sum_{i=1}^m|\alpha_i-1|^3+\sum_{r=1}^{\bar m}|\beta_r|^3+\sum_{i=1}^m\frac{1}{\l^3_i}+|f(w)|+||w||^2\right),
\end{split}
\end{equation*}
where \;$O\left(1\right)$ means here \;$O_{\bar\alpha, A, \bar \l, \bar \beta, w, \epsilon}\left(1\right)$ \;with \;$\bar\alpha=(\alpha_1, \cdots, \alpha_m)$, $A:=(a_1, \cdots, a_m)$ $\bar \l:=(\l_1, \cdots, \l_m)$, $\bar \beta:=(\beta_1, \cdots, \beta_{\bar m})$ and for $i=1, \cdots, m$, $\tilde \tau_i$ is as in Lemma \ref{eq:energyest}.  Thus after minimization in the $w$ variable, we have
\begin{equation*}
\begin{split}
&J(\sum_{i=1}^m\alpha_i\varphi_{a_i, \l_i}+\sum_{r=1}^{\bar m}\beta_r(v_r-\ov{(v_r)}_Q)+\bar w(\bar \alpha, A, \bar \l, \bar \beta))=-\frac{40}{3}m\pi^2-8m\pi^2\log(\frac{m\pi^2}{6})-8\pi^2\mathcal{F}_K(a_1, \dots, a_m)\\&+\sum_{i=1}^m(\alpha_i-1)^2\left(32\pi^2\log\l_i+16\pi^2 H(a_i, a_i)+\frac{472\pi^2}{9}\right)+\sum_{r=1}^{\bar m}\mu_r\beta_r^2\\&+16\pi^2\sum_{i=1}^m(\alpha_i-1)\left[\sum_{r=1}^{\bar m}2\beta_r(v_r-\ov{(v_r)}_Q)(a_i)-\sum_{ j=1, j\neq i}^m (\alpha_j-1)G(a_i, a_j)\right]\\&-2\pi^2\sum_{i=1}^m\frac{1}{\l_i^2}\left(\frac{\D_{g_{a_i}}\mathcal{F}^{A}_i(a_i)}{\mathcal{F}^{A}_i(a_i)}-\frac{2}{3}R_g(a_i)\right)+2\pi^2\sum_{i=1}^m\frac{\tilde \tau_i}{\l_i^2}\left(\frac{\D_{g_{a_i}}\mathcal{F}^{A}_i(a_i)}{\mathcal{F}^{A}_i(a_i)}-\frac{2}{3}R_g(a_i)\right)\\&+8\pi^2\sum_{i=1}^m\log(1-\tilde\tau_i)+O\left(\sum_{i=1}^m|\alpha_i-1|^3+\sum_{r=1}^{\bar m}|\beta_r|^3+\sum_{i=1}^m\frac{1}{\l^3_i}++||\bar w(\bar \alpha, A, \bar \l, \bar \beta)||^2\right),
\end{split}
\end{equation*}
where \;$O\left(1\right)$ means here \;$O_{\bar\alpha, A, \bar \l, \bar \beta, \epsilon}\left(1\right)$ \;with \;$\bar\alpha=(\alpha_1, \cdots, \alpha_m)$, $A:=(a_1, \cdots, a_m)$ $\bar \l:=(\l_1, \cdots, \l_m)$, $\bar \beta:=(\beta_1, \cdots, \beta_{\bar m})$ and for $i=1, \cdots, m$, $\tilde \tau_i$ is as in Lemma \ref{eq:energyest}.
where $\bar w(\bar \alpha, A, \bar \l, \bar \beta)$ is as in Corollary \ref{eq:cexpansionj1}.
\end{lem}
\begin{pf}
It follows directly from Lemma \ref{eq:energyest} and Proposition \ref{eq:expansionJ1}.
\end{pf}
\vspace{4pt}

\noindent
Lemma \ref{eq:expansionj} imply the following two Corollaries.
\begin{cor}\label{eq:energybddinf}(Energy bound at infinity)\\
Assuming that $\eta$ is a small positive real number with $0<2\eta<\varrho$ where $\varrho$ is as in \eqref{eq:cutoff},  then there exists $C_0^m:=C^m_0(\eta)$ such that  for every  $0<\epsilon\leq \epsilon_0$ where $\epsilon_0$ is as in \eqref{eq:mini}, there holds
\begin{equation*}
V(m, \epsilon, \eta)\subset J^{C_0^m}\setminus J^{-C_0^m}.
\end{equation*}
\end{cor}
\begin{pf}
It follows directly from \eqref{eq:para}-\eqref{eq:afpara}, \eqref{eq:coreql7}, Lemma \ref{eq:tauiest}, Proposition \ref{eq:expansionJ1} and Lemma \ref{eq:expansionj}.
\end{pf}
\begin{cor}\label{eq:energybddcrit} (Energy bound on critical set)\\
There exists a large positive constant \;$C_1^m$\; such that
\begin{equation*}
Crit(J)\subset J^{C_1^m}\setminus J^{-C_1^m}.
\end{equation*}
\end{cor}
\begin{pf}
It follows, via a contradiction argument, from the the fact that $J$ is invariant by translation by constants, Proposition \ref{eq:escape}, and  Lemma \ref{eq:energybddinf}.
\end{pf}
\vspace{6pt}

\noindent
Now, we are ready to characterize the topology of very high and very negative sublevels of $J$. We start with the very high ones. Indeed, we will combine Lemma \ref{eq:deformlemr}, Corollary \ref{eq:energybddinf} and Corollary \ref{eq:energybddcrit} with an argument of Malchiodi\cite{maldeg} to show the following lemma. We anticipate that with the above results at hand, the argument  of Malchiodi\cite{maldeg} carries over without any modifications, and for this reason we will not give the proof of the Lemma.
\begin{lem}\label{eq:tophigh}(Contractibility of high sublevels)\\
Assuming that $\eta$ is a small positive real number with $0<2\eta<\varrho$ where $\varrho$ is as in \eqref{eq:cutoff}, then there exists a large positive constant  $L^m:=L^m(\eta)$ with \;$L^m>2\max\{C_0^m, C_1^m\}$\;  such that for every \;$L\geq L^m$, we have  that \;$J^L$\; is a deformation retract of $W^{2, 2}(M)$, and hence it has the homology of a point, where $C^m_0$ is as in Lemma \ref{eq:energybddinf} and $C^m_1$ as in Lemma \ref{eq:energybddcrit}.
\end{lem}
\vspace{6pt}

\noindent
Now, to finish this section, we characterize the topology of very negative sublevels of $J$ when $m\geq 2$ or $\bar m\geq 1$. Indeed, we have:
\begin{lem}\label{eq:topnegative} (Characterization of topology of very negative sublevels)\\
Assuming that $m\geq 2$ or $\bar m\geq 1$, and $\eta$ is a small positive real number with $0<2\eta<\varrho$ where $\varrho$ is as in \eqref{eq:cutoff}, then there exists a large positive constant $L_{m, \bar m}:=L_{m, \bar m}(\eta)$ with  $L_{m, \bar m}>2\max\{C_0^m, C_1^m\}$ such that for every \;$L\geq L_{m, \bar m}$, we have  that \;$J^{-L}\;$  has the same homotopy type as  $B_{m-1}(M)$ if $m\geq 2$ and $\bar m=0$,  as\;$A_{m-1, \bar m}$ if $m\geq 2$ and $\bar m\geq 1$ and  as $S^{\bar m-1}$ if $m=1$ and $\bar m\geq 1$, where $C^m_0$ is as in Lemma \ref{eq:energybddinf} and $C^m_1$ as in Lemma \ref{eq:energybddcrit} .
\end{lem}
\begin{pf}
Using the work of Djadli-Malchiodi\cite{dm} and the one of Malchiodi\cite{maldeg}, we infer that  there exists $\tilde L_{m, \bar m}>0$ such that for every $L\geq\tilde L_{m, \bar m}$ and  for every $t\in ]\frac{m-1}{m}, 1[, $ $J^{-L}_t$ (for the definition see $\eqref{eq:jt}$) is homotopy equivalent to   $B_{m-1}(M)$ if $m\geq 2$ and $\bar m=0$,  to \;$A_{m-1, \bar m}$ if $m\geq 2$ and $\bar m\geq 1$ and to  $S^{\bar m-1}$ if $m=1$ and $\bar m\geq 1$. Next, applying Proposition \ref{eq:escape}, Lemma \ref{eq:energybddinf}, and Lemma \ref{eq:energybddcrit}, we have that for every $L\geq 2\max\{\tilde L_{m, \bar m}, C^m_0, C_1^m\}$ (with $C^m_0$ given by Lemma \ref{eq:energybddinf}, and $C_1^m$ given by Lemma \ref{eq:energybddcrit}), there exists $\tau=\tau(L)$ such that
\begin{equation}\label{eq:criticinf}
\n_uJ_t(u)\neq 0 \;\;\forall (t, u)\in [1-\tau, 1]\times J^{-L},
\end{equation}
where $J_t$ is as in \eqref{eq.jt1}. Thus, using the implicit function theorem, we infer that for every $s\in[1-\tau, 1]$, $\partial J^{-L}_s$ is a smooth submanifold in $W^{2, 2}(M)$ of codimension $1$. On the other hand, it is easy to see that the following inclusion $ J^{-L}_s\subset J^{-L}_{s^{'}}$ holds for every $s, s^{'}\in [1-\tau, 1]$ verifying $s\leq s^{'}$ . Hence the $\partial J^{-L}_s$'s foliate $J^{-L}\setminus J^{-L}_{1-\tau}$.  Now we are going to show that $J^{-L}$ retracts by deformation onto $J^{-L}_{1-\tau}$. We recall that the $J_t$'s are invariant by translations by constants. Thus from now until the end  of the proof of the fact $J^{-L}$ retracts by deformation onto $J^{-L}_{1-\tau}$, we will be working without lost of generality with functions $u$ satisfying $\ov{(u)}_Q=0$. Now, in order to achieve our goal, for every $u\in \ov{J^{-L}\setminus J^{-L}_{1-\tau}}$, we call $s(u)$ the unique element in $[1-\tau, 1]$ such that $u\in \partial J^{-L}_{s(u)}$. By the above discussions, we have that the map $u\rightarrow s(u)$ is continuous. To continue, we fix some notations. First, using \eqref{eq:jt1} and \eqref{eq:normpr}, we write
$$
J_t=||u||^2-tI(u),
$$
with
$$
I(u)=8m\pi^2\log \int_M Ke^{4u}dV_g-4\int_MQ_gudV_g-2\sum_{k=1}^{\bar m}\mu_k (u^k)^2,
$$
Next, we define a vector field which is a modification of the vector field defined by Lucia\cite{lu} (and inspired from the work of Bahri on contact forms, see \cite{bah1} or \cite{bah2}), see also \cite{maldeg}. Precisely, we set
$$
W_{ma}(u):=-\omega_{\tau}\left(\frac{|\n J_{s(u)}|}{|\n I(u)|}\right)\n J_{s(u)}+Z_{s(u)}(u),
$$
where
$$
Z_{s(u)}=-\lbrace\lbrace|\n I|\n J_{s(u)}+|\n J_{s(u)}|\n I(u)\rbrace,
$$
and $\omega_{\tau}$ a smooth non-decreasing cut-off function satisfying $\omega_{\tau}:\R\rightarrow[0,1]$ $\omega_{\tau}(t)=0$ for $t\leq \tau$ and $\omega_{\tau}(t)=1$ for $t\geq 2\tau$.
Now, using arguments similar to the ones \cite{lu}, see also \cite{maldeg}, we have that every trajectory of $W$ stays bounded as long as $s(u)\in[1-\tau, 1]$. To see this, let us consider the following Cauchy problem
$$
\begin{cases}
\frac{d}{dt}u=W_{ma}(u(t)),\\
u(0)=u_0.
\end{cases}
$$
Like in \cite{maldeg}, by arguing as in \cite{lu}, we obtain
\begin{equation}\label{eq:diffj}
\frac{d}{dt}|_{t=0}[I(u(t)]\leq -\frac{1}{\tau}\frac{d}{dt}|_{t=0}[J_{s(u_0)}(u(t))].
\end{equation}
Now, by direct calculations, we get
$$
\frac{d}{dt}|_{t=0}[J_{s(u(t))}]=\frac{d}{dt}|_{t=0}[J_{s(u_0)}]-\frac{d}{dt}|_{t=0}s(u(t))I(u(t)).
$$
Next, using the definition of $s(u)$, we derive
$$
\frac{d}{dt}|_{t=0}[J_{s(u(t))}]=0.
$$
So, we obtain
\begin{equation}\label{eq:derivative}
\frac{d}{dt}|_{t=0}[J_{s(u_0)}]=\frac{d}{dt}|_{t=0}s(u(t))I(u(t)).
\end{equation}
Thus, since the last three formulas holds for each initial data, then we have that $\eqref{eq:diffj}$ implies
$$
I(u(t))\leq I(u(0))e^{\frac{1}{\tau}(s(u(t))-s(u(0)))}.
$$
Hence, since $J(u(t))\leq -L$, then  as soon as $s(u(t))\in [1-\tau, 1]$, we have uniform bounds on the trajectory $u(t)$, depending on the initial value.  Next,  it is easy to see that $s(u(t))$ is non-increasing. But we have even more, we claim that for each trajectory $u(t)$, $\frac{d}{dt}s(u(t))$ is bounded below by a negative constant. Assuming the contrary, and recalling that $I(u(t))$ is bounded, then  from $\eqref{eq:derivative}$, we obtain a bounded  sequence $u_l$ such that $\n J_{s(u_l)}[W_{ma}(u_l)]\rightarrow 0$. Now, using the same arguments as in \cite{lu}, page 11 (see also \cite{maldeg}), we get a contradiction to $\eqref{eq:criticinf}$. In summary, we have that each point of $\partial J^{-L}$ reaches $\partial J^{-L}_{1-\tau}$ in finite time. Viceversa, using the inverse flow, we have also that each point in $\partial J^{-L}_{1-\tau}$ reaches $\partial J^{-L}$ in finite.. Hence, by using a suitable reparametrization in $t$, one has a deformation retract of $J^{-L}$ onto $J^{-L}_{1-\tau}$. Now we set  $L_{m, \bar m}:=2\max\{\tilde L_{m ,\bar m}, 2C^m_0, 2C_1^m\}>2\max\{C_0^m, C_1^m\}$  and have from the above analysis that for every $L\geq L_{m, \bar m}$, the exists $\tau=\tau(L)$ such that $J^{-L}_{1-\tau} $ has the same homotopy type as  $B_{m-1}(M)$ if $m\geq 2$ and $\bar m=0$,  as\;$A_{m-1, \bar m}$ if $m\geq 2$ and $\bar m\geq 1$ and  as $S^{\bar m-1}$ if $m=1$ and $J^{-L}$ retracts by deformation onto $J^{-L}_{1-\tau}$. Hence, we have $J^{-L}$ has the same homotopy type as  $B_{m-1}(M)$ if $m\geq 2$ and $\bar m=0$,  as\;$A_{m-1, \bar m}$ if $m\geq 2$ and $\bar m\geq 1$ and  as $S^{\bar m-1}$ if $m=1$ and $\bar m\geq 1$, thereby ending the proof of the Lemma.
\end{pf}

\section{Morse Lemma at infinity for $J$}\label{eq:morseinfinity}\label{eq:morseleminf}
In this section, we derive a Morse Lemma at infinity for $J$. In order to do that, we start by constructing a pseudo-gradient for $J(\bar \alpha, A, \bar \l, \bar \beta)$, where $J(\bar \alpha, A, \bar \l, \bar \beta)$ is defined as in \eqref{eq:finitedj} . Indeed, we have:
\begin{pro}\label{eq:conspseudograd}
Assuming that $\eta$ is a small positive real number with $0<2\eta<\varrho$ where $\varrho$ is as in \eqref{eq:cutoff}, and $0<\epsilon\leq \epsilon_0$ where $\epsilon_0$ is as in \eqref{eq:mini}, then there exists a pseudogradient $W_{no}$ of $J(\bar \alpha, A, \bar \l, \bar \beta)$ such that \\
1) For every $u:=\sum_{i=1}^m\alpha_i\varphi_{a_i, \l_i}+\sum_{r=1}^{\bar m}\beta_r (v_r-\ov{(v_r)}_Q)\in V(m, \epsilon, \eta)$ with the concentration  points $a_i$, the masses $\alpha_i$, the concentrating parameters $\l_i$ ($i=1, \cdots, m$)  and the negativity parameters $\beta_r$ ($r=1, \cdots, \bar m$) satisfying \eqref{eq:afpara}, there holds
\begin{equation}\label{eq:pseudoexact}
<-\n J(u), W_{no}>\geq c\left(\sum_{i=1}^m\frac{1}{\l_i^2}+\sum_{i=1}^m\frac{|\n_g\mathcal{F}^{A}_i(a_i)|}{\l_i}+\sum_{i=1}^m|\alpha_i-1|+\sum_{i=1}^m|\tau_i|+\sum_{r=1}^{\bar m}|\beta_r|)\right),
\end{equation}
and  for every $u:=\sum_{i=1}^m\alpha_i\varphi_{a_i, \l_i}+\sum_{r=1}^{\bar m}\beta_r (v_r-\ov{(v_r)}_Q)+\bar w(\bar \alpha, A, \bar \l, \bar \beta)\in V(m, \epsilon, \eta)$ with the concentration  points $a_i$, the masses $\alpha_i$, the concentrating parameters $\l_i$ ($i=1, \cdots, m$)  and the negativity parameters $\beta_r$ ($r=1, \cdots, \bar m$) satisfying \eqref{eq:afpara}, and $\bar w(\bar \alpha, A, \bar \l, \bar \beta)$ is as in \eqref{eq:minj},  there holds
\begin{equation}\label{eq:pseudoperturb}
<-\n J(u+\bar w), W_{no}+\frac{\partial \bar w(W_{no})}{\partial (\bar\alpha, A, \bar \l, \bar \beta)}>\geq c\left(\sum_{i=1}^m\frac{1}{\l_i^2}+\sum_{i=1}^m\frac{|\n_g\mathcal{F}^{A}_i(a_i)|}{\l_i}+\sum_{i=1}^m|\alpha_i-1|+\sum_{i=1}^m|\tau_i|+\sum_{r=1}^{\bar m}|\beta_r|\right),
\end{equation}
where $c$ is a small positive constant independent of $A:=(a_1, \cdots, a_m)$, $\bar\alpha=(\alpha_1, \cdots, \alpha_m)$, $\bar\l=(\l_1, \cdots, \l_m)$, $\bar \beta=(\beta_1, \cdots, \beta_{\bar m})$ and $\epsilon$.\\
2) $W_{no}$ is a $W^{2, 2}$-bounded vector field and is compactifying outside the region where $A$ is very close to a critical point $B$ of $\mathcal{F}_K$ satisfying $l_K(B)<0$.
\end{pro}
\begin{pf}
To construct the required vector field, we will make use of the following ones whose description are as follows. First of all, to move the $\alpha_i's$ we use the following vector fields
\begin{equation}\label{eq:walphai}
W_{\alpha_i}:=-\frac{\alpha_i-1}{|\alpha_i-1|}\left(1-\chi(\l_i^2|\alpha_i-1|)\right)Z_i, \;\;\;i=1,\cdots, m,
\end{equation}
with
\begin{equation}\label{eq:zi}
Z_i=\frac{1}{\log \l_i}\left(\varphi_{a_i, \l_i}-\sum_{j=1, j\neq i}^m G(a_i, a_j)\l_j \frac{\partial \varphi_{a_j, \l_j}}{\partial \l_j}-\frac{1}{\alpha_i}\left[2\log \l_i-\frac{5}{6}+H(a_i, a_i)\right]\l_i\frac{\partial \varphi_{a_i, \l_i}}{\partial \l_i}\right),
\end{equation}
where $\chi$ is a  smooth cut-off function verifying
\begin{equation}
\chi(t)=1\;\,\text{for}\;\; |t|\leq 1,\;\;\;\,\,\chi(t)=0\;\;\;\text{for}\;\; |t|\geq 2\;\;\;\text{and}\;\;\;0\leq\chi(t)\leq 1\;\;\text{for} \;\;|t|\geq 0.
\end{equation}
Now, setting
\begin{equation}\label{eq:wbaralpha}
W_{\bar \alpha}=\sum_{i=1}^m W_{\alpha_i},
\end{equation}
and using Lemma \ref{eq:gradientalpha}, we obtain
\begin{equation}\label{eq:estwbaralpha}
<-\n J(u), W_{\bar \alpha}>\geq c\sum_{i=1}^m|\alpha_i-1|(\left(1-\chi(\l_i^2|\alpha_i-1|\right)+O\left(\sum_{i=1}^m\frac{|\alpha_i-1|}{\log\l_i}+\sum_{i=1}^m\frac{1}{\l_i^2}+(\sum_{r=1}^{\bar m}|\beta_r|)(\sum_{i=1}^{m}\frac{1}{\log\l_i})\right),
\end{equation}
 where $c$ is a small positive constant independent of $A:=(a_1, \cdots, a_m)$, $\bar\alpha=(\alpha_1, \cdots, \alpha_m)$, $\bar\l=(\l_1, \cdots, \l_m)$, $\bar \beta:=(\beta_1, \cdots, \beta_{\bar m})$, and $\epsilon$.. Next, to move the points $a_i$'s, we will make use of the following vector field
 \begin{equation}\label{eq:wai}
 W_{a_i}:=\frac{\n_g\mathcal{F}_i^{A}(a_i)}{|\n_g \mathcal{F}_i^{A}(a_i)|}\frac{1}{\l_i}\frac{\partial \varphi_{a_i, \l_i}}{\partial a_i}\left(1-\zeta(\l_i |\n_g\mathcal{F}^{A}(a_i)|)\right),
 \end{equation}
 where $\zeta$ is a smooth cutt-off function verifying
 \begin{equation}\label{eq:zeta}
 \zeta(t)=1\;\;\;\;\text{for}\;\; |t|\leq C_0\;\;\;\;\;\zeta(t)=0\,\;\;\;\;\;\;\text{for} \;\;|t|\geq 2C_0,\;\;\text{and}\;\;\;0\leq \zeta(t)\leq 1\;\;\text{for} \;\;|t|\geq0.
 \end{equation}
 where $C_0$ is a large positive constant to be fixed later. Now, setting
 \begin{equation}\label{eq:wa}
 W_{A}:=\sum_{i=1}^mW_{a_i}
 \end{equation}
 with $A:=(a_1, \cdots, a_m)$, and using Lemma \ref{eq:gradientaest}, we obtain the following estimate
 \begin{equation}\label{eq:estwa}
 \begin{split}
 <-\n J(u), W_A>\geq c\sum_{i=1}^m\frac{|\n_g \mathcal{F}^{A}_{i}(a_i)|}{\l_i}\left(1-\zeta(\l_i |\n_g\mathcal{F}^{A}(a_i)|)\right)\\+O\left(\sum_{i=1}^m\frac{1}{\l_i^2}+\sum_{i=1}^m\tau_i^2+\sum_{i=1}^m |\alpha_i-1|^2+\sum_{r=1}^{\bar m}|\beta_r|^2\right).
 \end{split}
 \end{equation}
 To move the $\beta_r$'s, we will consider the following vector fields
 \begin{equation}\label{eq:wbi}
 W_{\beta_r}:=v_r-\ov{(v_r)}_Q.
 \end{equation}
 Next, using Lemma \ref{eq:gradientbeta}, we get
 \begin{equation}\label{eq:estwb}
 <-\n J(u), W_{\beta_r}>=2|\mu_r|\beta_r+O\left(\sum_{i=1}^m|\alpha_i-1|+\sum_{i=1}^m|\tau_i|+\sum_{i=1}^m\frac{1}{\l_i^2}\right)
 \end{equation}
 Finally to move the $\l_i$'s, we will make use of the following vector fields
 \begin{equation}\label{eq:wli}
 W_{\l_i}:=-\frac{1}{\alpha_i}\l_i\frac{\partial \varphi_{a_i, \l_i}}{\partial \l_i}.
 \end{equation}
 Now, setting
 \begin{equation}\label{eq:wlambda}
 W_{\bar \l}:=\sum_{i=1}^mW_{\l_i},
 \end{equation}
 with $\bar \l:=(\l_1, \cdots, \l_m)$, and using  Lemma \ref{eq:gradientlambdaest} and Corollary \ref{eq:cgradientlambdaest}, we get
 \begin{equation}\label{eq:estwlambdai}
 <-\n J(u), W_{ \l_i}>=32\pi^2\tau_i+O\left(\sum_{i=1}^m|\alpha_i-1|^2+\sum_{r=1}^{\bar m}|\beta_r|^3+\sum_{i=1}^{m}\frac{1}{\l_i^2}\right),
 \end{equation}
 and
 \begin{equation}\label{eq:estwbarl}
 \begin{split}
 <-\n J(u), W_{\bar \l}>=4\pi^2\sum_{i=1}^m\frac{1}{\l_i^2}\left(\frac{\D_{g_{a_i}} \mathcal{F}^{A}_i(a_i)}{\mathcal{F}_i^{A}(a_i)}-\frac{2}{3}R_g(a_i)\right)\\+O\left(\sum_{i=1}^m|\alpha_i-1|^2+\sum_{i=1}^m\tau_i^3+\sum_{r=1}^{\bar m}|\beta_r|^3+\sum_{i=1}^m\frac{1}{\l_i^3}\right).
 \end{split}
 \end{equation}
 Next, in order to make use of the vector fields $W_{\alpha_i}$, $W_{a_i}$ and $W_{\lambda_i}$ and $W_{\beta_r}$  to define the desired pseudogradient for $J(\bar \alpha, A, \bar \l, \bar \beta)$, we first divide $V(m, \epsilon, \eta)$ into five regions (with non empty intersections). After, we construct an appropriate vector field on each region and make convex combination of all of them to obtain the desired pseudogradient. To begin,  we first divide $V(m, \epsilon, \eta)$ in two regions $\hat{V}(m, \epsilon, \eta)$ and $\tilde{V}(m, \epsilon, \eta)$ (with non empty intersection) as follows
 \begin{equation}\label{eq:hatw}
 \begin{split}
& \hat{V}(m, \epsilon, \eta)=\{u\in V(m, \epsilon, \eta):\:\;\;\;\exists r\in \{1, \cdots, \bar m\}\,\,\text{such that}\;\,|\beta_r|\geq\\& \tilde C_0\left(\sum_{i=1}^m\frac{|\n_g \mathcal{F}^{A}(a_i)|}{\l_i}+\sum_{i=1}^m|\alpha_i-1|+\sum_{i=1}^m|\tau_i|+\sum_{i=1}^m\frac{1}{\l_i^2}\right)\}
 \end{split}
 \end{equation}
 where $\tilde C_0$ is a large positive constant to be chosen later and
 and
 \begin{equation}\label{eq:tildew}
 \begin{split}
 &\tilde{V}(m, \epsilon, \eta)=\{u\in V(m, \epsilon, \eta):\:\;\;\;\forall  r\in \{1, \cdots, \bar m\}\,\,\text{there holds}\;\,|\beta_r|\leq \\&2\tilde C_0\left(\sum_{i=1}^m\frac{|\n_g \mathcal{F}^{A}(a_i)|}{\l_i} +\sum_{i=1}^m|\alpha_i-1|+\sum_{i=1}^m|\tau_i|+\sum_{i=1}^m\frac{1}{\l_i^2}\right\}
 \end{split}
 \end{equation}
Next, we are going to divide again $V(m, \epsilon, \eta)$ into four regions $V_1(m, \epsilon, \eta)$, $V_2(m, \epsilon, \eta)$, $V_3(m, \epsilon, \eta)$,  and $V_4(m, \epsilon, \eta)$ (with non empty intersections) as follows
 \begin{equation}\label{eq:v1}
 V_1(m, \epsilon, \eta):=\{u\in V(m, \epsilon, \eta): \;\;\exists \;i\in\{1, \cdots, m\}: \;\;\;|\tau_i|\geq \frac{\hat C_0}{\l_i^2}\}.
 \end{equation}
 where $\hat C_0$ is a large positive constant to be chosen later,
 \begin{equation}\label{eq:v2}
\begin{split}
V_2(m, \epsilon, \eta):=\{u\in V(m, \epsilon, \eta): \;\;\forall \;j\in \{1, \cdots, m\}\;\,\;|\tau_j|\leq 2\frac{\hat C_0}{\l_j^2},  \;\;\exists i\in \{1, \cdots, m\}\;\; \text{such that}\\\;\;\;\frac{|\n_g \mathcal{F}^{A}_i(a_i)|}{\l_i}\geq 2\frac{C_0}{\l_i^2}\;\}.
\end{split}
\end{equation}
\begin{equation}\label{eq:v3}
 \begin{split}
 V_3(m, \epsilon, \eta):=\{u\in V(m, \epsilon, \eta): l_K(A)>0, \;\;\forall \;i\in \{1, \cdots, m\}\;\,\;|\tau_i|\leq 2\frac{\hat C_0}{\l_i^2},  \;\; \text{and}\\\;\;\;\frac{|\n_g \mathcal{F}^{A}_i(a_i)|}{\l_i}\leq 4\frac{C_0}{\l_i^2}\;\}
 \end{split}
 \end{equation}
 and
 \begin{equation}\label{eq:v4}
\begin{split}
 V_4(m, \epsilon, \eta):=\{u\in V(m, \epsilon, \eta): l_K(A)<0, \;\;\forall \;i\in \{1, \cdots, m\}\;\,\;|\tau_i|\leq 2\frac{\hat C_0}{\l_i^2},  \;\; \text{and}\\\;\;\;\frac{|\n_g \mathcal{F}^{A}_i(a_i)|}{\l_i}\leq 4\frac{C_0}{\l_i^2}\;\}.
 \end{split}
 \end{equation}
 On the other hand,  using \eqref{eq:auxiindexa1} and our  non degeneracy assumption,  it is easy to see that
 \begin{equation}\label{eq:recover}
 V(m, \epsilon, \eta)=\hat{V}(m, \epsilon, \eta)\cup V_s(m, \epsilon, \eta)\cap \tilde {V}(m, \epsilon, \eta).
 \end{equation}
 The five regions in question are the following ones $\hat{V}(m, \epsilon, \eta)$, $V_s(m, \epsilon, \eta)\cap \tilde {V}(m, \epsilon, \eta)$, $s=1, \cdots, 4$.
We start with the region $ \hat{V}(m, \epsilon, \eta)$. First of all,  for every $u\in  \hat{V}(m, \epsilon, \eta)$, we associate the set
 \begin{equation}\label{eq:hatf}
 \begin{split}
 &\hat{F}:=\{r\in \{1, \cdots, \bar m\}\,\,\text{such that}\;\,|\beta_r|\geq  \tilde C_0\left(\sum_{i=1}^m|\alpha_i-1|+\sum_{i=1}^m|\tau_i|+\sum_{i=1}^m\frac{1}{\l_i^2}+\sum_{i=1}^m\frac{|\n_g \mathcal{F}^{A}(a_i)|}{\l_i}\right)\},
 \end{split}
 \end{equation}
 and define the vector field that we are going to consider in $\hat{V}(m, \epsilon, \eta)$ as follows
 \begin{equation}\label{eq:whatf}
 W_{\hat{F}}:=\sum_{r\in \hat{F}} \frac{\beta_r}{|\beta_r|}W_{\beta_r}.
 \end{equation}
Thus, combining \eqref{eq:estwb}, \eqref{eq:hatf}, \eqref{eq:whatf}, and taking $\tilde C_0$ very large, we derive that $u:=\sum_{i=1}^m\alpha_i\varphi_{a_i, \l_i}+\sum_{r=1}^{\bar m}\beta_r (v_r-\ov{(v_r)}_Q)\in \hat{V}(m, \epsilon, \eta)$ with the concentration  points $a_i$,  the masses $\alpha_i$, the concentrating parameters $\l_i$ ($i=1, \cdots, m$)  and the negativity parameters $\beta_r$ ($r=1, \cdots, \bar m$) satisfying \eqref{eq:afpara}, there holds
 \begin{equation}\label{eq:estwhatf}
 \begin{split}
<-\n J(u), W_{\hat{F} }>&=2\sum_{r\in \hat {F}}|\mu_r||\beta_r|+O\left(\sum_{i=1}^m|\alpha_i-1|+\sum_{i=1}^m|\tau_i|+\sum_{i=1}^m\frac{1}{\l_i^2}\right)\\&\geq c\left(\sum_{i=1}^m\frac{1}{\l_i^2}+\sum_{i=1}^m\frac{|\n_g\mathcal{F}^{A}_i(a_i)|}{\l_i}+\sum_{i=1}^m|\alpha_i-1|+\sum_{i=1}^m|\tau_i|+\sum_{r=1}^{\bar m}|\beta_r|\right),
\end{split}
 \end{equation}
 where $c$ is a small positive constant independent of $\bar \alpha:=(\alpha_1, \cdots, \alpha_m)$, $A:=(a_1, \cdots, a_m)$, $\bar \l:=(\l_1, \cdots, \l_m)$, $\bar \beta:=(\beta_1, \cdots, \beta_{\bar m})$ and $\epsilon$. Now, we are going to consider the region $V_1(m, \epsilon, \eta)\cap \tilde{V}(m, \epsilon, \eta)$. For this end,   for every $u\in  V_1(m, \epsilon, \eta)$, we associate the set
 \begin{equation}\label{eq:f1}
 F_1:=\{\;i\in\{1, \cdots, m\}: \;\;\;|\tau_i|\geq \frac{\hat C_0}{\l_i^2}\},
 \end{equation}
 and define the following vector field in $V_1(m, \epsilon, \eta)$
 \begin{equation}\label{eq:wf1}
 W_{F_1}:=\sum_{i\in F_1}\frac{\tau_i}{|\tau_i|}W_{\l_i}.
 \end{equation}
 Thus, combining \eqref{eq:estwlambdai} and \eqref{eq:wf1}, we get that for $u:=\sum_{i=1}^m\alpha_i\varphi_{a_i, \l_i}+\sum_{r=1}^{\bar m}\beta_r (v_r-\ov{(v_r)}_Q)\in V_1(m, \epsilon, \eta)$ with the concentration  points $a_i$,  the masses $\alpha_i$, the concentrating parameters $\l_i$ ($i=1, \cdots, m$)  and the negativity parameters $\beta_r$ ($r=1, \cdots, \bar m$) satisfying \eqref{eq:afpara}, there holds
 \begin{equation}\label{eq:estwf1}
 \begin{split}
 <-\n J(u), W_{F_1}>=&\sum_{i\in F_1}32\pi^2|\tau_i|+O(\left(\sum_{i=1}^m|\alpha_i-1|^2+\sum_{r=1}^{\bar m}|\beta_r|^3+\sum_{i=1}^m\frac{1}{\l_i^2}\right)\\&\geq c\sum_{i\in F_1}|\tau_i|+O(\left(\sum_{i=1}^m|\alpha_i-1|^2+\sum_{r=1}^{\bar m}|\beta_r|^3+\sum_{i=1}^m\frac{1}{\l_i^2}\right),
 \end{split}
 \end{equation}
 where $c$ is a small positive constant independent of $A:=(a_1, \cdots, a_m)$, $\bar \alpha=:(\alpha_1, \cdots, \alpha_m)$, $\bar \l:=(\l_1, \cdots, \l_m)$, $\bar \beta:=(\beta_1, \cdots, \beta_{\bar m})$, and $\epsilon$.
Now, we are ready to define the vector field that we are going to consider in $V_1(m, \epsilon, \eta)\cap \tilde{V}(m, \epsilon, \eta)$ as follows
 \begin{equation}\label{eq:hatw1}
 \tilde{W}_1:=W_{\bar \alpha}+W_{A}+W_{F_1}.
\end{equation}
So, using \eqref{eq:estwbaralpha}, \eqref{eq:estwa} , \eqref{eq:tildew}, \eqref{eq:estwf1}, and taking $\hat C_0$ large enough, we get that $u:=\sum_{i=1}^m\alpha_i\varphi_{a_i, \l_i}+\sum_{r=1}^{\bar m}\beta_r (v_r-\ov{(v_r)}_Q)\in V_1(m, \epsilon, \eta)\cap\hat{V}(m, \epsilon, \eta)$ with the concentration  points $a_i$,  the masses $\alpha_i$, the concentrating parameters $\l_i$ ($i=1, \cdots, m$)  and the negativity parameters $\beta_r$ ($r=1, \cdots, \bar m$) satisfying \eqref{eq:afpara}, there holds
 \begin{equation}\label{eq:esthatw1}
 <-\n J(u), \tilde{W}_{1}>\geq c\left(\sum_{i=1}^m\frac{1}{\l_i^2}+\sum_{i=1}^m\frac{|\n_g\mathcal{F}^{A}_i(a_i)|}{\l_i}+\sum_{i=1}^m|\alpha_i-1|+\sum_{i=1}^m|\tau_i|+\sum_{r=1}^{\bar m}|\beta_r|\right),
 \end{equation}
 where $c$ is a small positive constant independent of $A:=(a_1, \cdots, a_m)$, $\bar \alpha=:(\alpha_1, \cdots, \alpha_m)$, $\bar \l:=(\l_1, \cdots, \l_m)$, $\bar \beta:=(\beta_1, \cdots, \beta_{\bar m})$, and $\epsilon$. Next, we consider the region $V_2(m, \epsilon, \eta)\cap\hat{V}(m, \epsilon, \eta)$.
For this, for  $u\in  V_2(m, \epsilon, \eta)$, we associate the set
 \begin{equation}\label{eq:f2}
 F_2:=\{\;i\in\{1, \cdots, m\}: \;\;\;\frac{|\n_g \mathcal{F}^{A}_i(a_i)|}{\l_i}\geq 2\frac{C_0}{\l_i^2}\}.
 \end{equation}
 Using \eqref{eq:estwa}, \eqref{eq:v2}, and \eqref{eq:f2}, we have that for every $u:=\sum_{i=1}^m\alpha_i\varphi_{a_i, \l_i}+\sum_{r=1}^{\bar m}\beta_r (v_r-\ov{(v_r)}_Q)\in V_2(m, \epsilon, \eta)$ with the concentration  points $a_i$,  the masses $\alpha_i$, the concentrating parameters $\l_i$ ($i=1, \cdots, m$)  and the negativity parameters $\beta_r$ ($r=1, \cdots, \bar m$) satisfying \eqref{eq:afpara}, the following estimate holds
 \begin{equation}\label{eq:estwa1}
<-\n J(u), W_A>\geq c\sum_{i\in F_2}\frac{|\n_g \mathcal{F}^{A}_i(a_i)|}{\l_i}+O\left(\sum_{i=1}^m\frac{1}{\l_i^2}+\sum_{i=1}^m|\alpha_i-1|^2+\sum_{r=1}^{\bar m}|\beta_r|^2\right),
 \end{equation}
where $c$ is a small positive constant independent of  $A:=(a_1, \cdots, a_m)$, $\bar \alpha=:(\alpha_1, \cdots, \alpha_m)$,  $\bar \l:=(\l_1, \cdots, \l_m)$, $\bar \beta:=(\beta_1, \cdots, \beta_{\bar m})$ and $\epsilon$.
  Now, we define the vector that we are going to consider in $ V_2(m, \epsilon, \eta)\cap \tilde{V}(m, \epsilon, \eta)$ as follows
 \begin{equation}\label{eq:w2}
 \tilde{W}_2:=W_{\bar \alpha}+W_A+W_{\bar \l}.
 \end{equation}
 Hence,  combining \eqref{eq:estwbaralpha}, \eqref{eq:estwbarl}, \eqref{eq:tildew}, \eqref{eq:v2}, and \eqref{eq:estwa1},  we obtain that
 we have that for every $u:=\sum_{i=1}^m\alpha_i\varphi_{a_i, \l_i}+\sum_{r=1}^{\bar m}\beta_r (v_r-\ov{(v_r)}_Q)\in V_2(m, \epsilon, \eta)\cap \tilde{V}(m, \epsilon, \eta)$ with the concentration  points $a_i$,  the masses $\alpha_i$, the concentrating parameters $\l_i$ ($i=1, \cdots, m$) and the negativity parameters $\beta_r$ ($r=1, \cdots, \bar m$) satisfying \eqref{eq:afpara}, the following estimate holds
\begin{equation}\label{eq:tildew21}
<-\n J(u), \tilde{W}_2>\geq c\sum_{i\in F_2}\frac{|\n_g \mathcal{F}^{A}_i(a_i)|}{\l_i}+c\sum_{i=1}^m|\alpha_i-1|(\left(1-\chi(\l_i^2|\alpha_i-1|\right)+O\left(\sum_{i=1}^m\frac{|\alpha_i-1|}{\log\l_i}+\sum_{i=1}^m\frac{1}{\l_i^2}\right),
\end{equation}
where $c$ is a small positive constant independent of $A:=(a_1, \cdots, a_m)$, $\bar\alpha=(\alpha_1, \cdots, \alpha_m)$, $\bar\l=(\l_1, \cdots, \l_m)$, $\bar \beta=(\beta_1, \cdots, \beta_{\bar m})$, and $\epsilon$. Thus, using again \eqref{eq:tildew}, \eqref{eq:v2}, the fact that the $\l_i$'s are comparable and taking $C_0$ very large, we infer that \eqref{eq:tildew21} implies that for every $u:=\sum_{i=1}^m\alpha_i\varphi_{a_i, \l_i}+\sum_{r=1}^{\bar m}\beta_r (v_k-\ov{(v_r)}_Q)\in V_2(m, \epsilon, \eta)\cap \tilde{V}(m, \epsilon, \eta)$ with the concentration  points $a_i$,  the masses $\alpha_i$, the concentrating parameters $\l_i$ ($i=1, \cdots, m$) and the negativity parameters $\beta_r$ ($r=1, \cdots, \bar m$) satisfying \eqref{eq:afpara}, the following estimate holds
\begin{equation}\label{eq:tildew21}
<-\n J(u), \tilde{W}_2>\geq c\left(\sum_{i=1}^m\frac{1}{\l_i^2}+\sum_{i=1}^m\frac{|\n_g\mathcal{F}^{A}_i(a_i)|}{\l_i}+\sum_{i=1}^m|\alpha_i-1|+\sum_{i=1}^m|\tau_i|+\sum_{r=1}^{\bar m}|\beta_r|\right),
\end{equation}
where $c$ is a small positive constant independent of $A:=(a_1, \cdots, a_m)$, $\bar\alpha=(\alpha_1, \cdots, \alpha_m)$, $\bar\l=(\l_1, \cdots, \l_m)$, $\bar \beta=(\beta_1, \cdots, \beta_{\bar m})$, and $\epsilon$. Now, we are going to deal with the region $V_3(m, \epsilon, \eta)\cap \tilde{V}(m, \epsilon, \eta)$. In this region, we are going to consider the following vector field
\begin{equation}\label{eq:tildew3}
\tilde{W}_3:=W_{\bar \alpha}+\Omega_3W_{\bar \l}+ W_{A},
\end{equation}
where $\Omega_3$ is a large positive constant to be chosen later. On the other hand, since in $V_3(m, \epsilon, \eta)$, we have $\forall \;i\in \{1, \cdots, m\}\;\,\;|\tau_i|\leq 2\frac{\hat C_0}{\l_i^2}$,  and the $\l_i$'s are comparable, then we derive that the following estimate holds
 \begin{equation}\label{eq:estrlij}
 \frac{\l_i^2\sqrt{\mathcal{F}^{A}_i(a_i)}}{\l_j^2\sqrt{\mathcal{F}^{A}_j(a_j)}}=1+O\left(\sum_{k=1}^m|\alpha_k-1|\log \l_k+\sum_{k=1}^m\frac{1}{\l_k^2}+\sum_{r=1}^{\bar m}|\beta_r|\right)\;\;\;\forall i, j=1, \cdots, m.
 \end{equation}
 Hence, using \eqref{eq:estrlij} and the fact that the $\l_i$'s  are comparable, we obtain
 \begin{equation}\label{eq:la}
 \begin{split}
 \sum_{i=1}^m\frac{1}{\l_i^2}\left(\frac{\D_{g_{a_i}} \mathcal{F}^{A}_i(a_i)}{\mathcal{F}_i^{A}(a_i)}-\frac{2}{3}R_g(a_i)\right)=\frac{1}{\l_1^2\sqrt{\mathcal{F}^{A}_1(a_1)}}l_K(A)\\+O\left(\sum_{i=1}^m|\alpha_i-1|\frac{\log \l_i}{\l_i^2}+\sum_{i=1}^m\frac{1}{\l_i^4}+(\sum_{r=1}^{\bar m}|\beta_r|)(\sum_{i=1}^m\frac{1}{\l_i^2})\right).
 \end{split}
 \end{equation}
Thus, combining \eqref{eq:estwbarl}, \eqref{eq:v3}, \eqref{eq:la}, and using Young's inequality, we infer that the following estimate holds in $V_3(m, \epsilon, \eta)$
\begin{equation}\label{eq:estwbarlref}
<-\n J(u), W_{\bar \l}>=\frac{4\pi^2 l_K(A)}{\l_1^2\sqrt{\mathcal{F}(a_i)}}+O\left(\sum_{i=1}^m|\alpha_i-1|^2+\sum_{i=1}^m\frac{1}{\l_i^3}+\sum_{r=1}^{\bar m}|\beta_r|^3\right)
\end{equation}
 Next, using \eqref{eq:estwbaralpha}, \eqref{eq:estwa}, \eqref{eq:v3}, and \eqref{eq:estwbarlref}, we deduce that in $V_3(m, \epsilon, \eta)$ the following estimate holds
 \begin{equation}\label{eq:esttildew3a}
 \begin{split}
 &<-\n J(u), \tilde{W}_{3}>\geq \Omega_3\frac{4\pi^2l_K(A)}{\l_1^2\sqrt{\mathcal{F}(a_i)}}+c\sum_{i=1}^m|\alpha_i-1|(\left(1-\chi(\l_i^2|\alpha_i-1|\right)\\&+c\sum_{i=1}^m\frac{|\n_g \mathcal{F}^{A}_{i}(a_i)|}{\l_i}\left(1-\zeta(\l_i |\n_g\mathcal{F}^{A}(a_i))|\right)+O\left(\sum_{i=1}^m\frac{|\alpha_i-1|}{\log \l_i}+\sum_{i=1}^m\frac{1}{\l_i^2}+\sum_{k=1}^{\bar m}|\beta_k|)(\sum_{i=1}^m\frac{1}{\log \l_i})\right),
 \end{split}
 \end{equation}
 where $c$ is a small positive constant independent of $A:=(a_1, \cdots, a_m)$, $\bar\alpha=(\alpha_1, \cdots, \alpha_m)$, $\bar\l=(\l_1, \cdots, \l_m)$, $\bar \beta=(\beta_1, \cdots, \beta_{\bar m})$ and $\epsilon$. Hence, using again \eqref{eq:v3} and \eqref{eq:tildew}, we have up to taking $\Omega_3$ large enough that \eqref{eq:esttildew3a} implies that for every $u:=\sum_{i=1}^m\alpha_i\varphi_{a_i, \l_i}+\sum_{r=1}^{\bar m}\beta_r (v_r-\ov{(v_r)}_Q)\in V_3(m, \epsilon, \eta)\cap \tilde V(m, \epsilon, \eta)$ with the concentration  points $a_i$,  the masses $\alpha_i$, the concentrating parameters $\l_i$ ($i=1, \cdots, m$) and the negativity parameters $\beta_r$ ($r=1, \cdots, \bar m$) satisfying \eqref{eq:afpara}, there holds
\begin{equation}\label{eq:esttildew3}
<-\n J(u), \tilde{W}_{3}>\geq c\left(\sum_{i=1}^m\frac{1}{\l_i^2}+\sum_{i=1}^m\frac{|\n_g\mathcal{F}^{A}_i(a_i)|}{\l_i}+\sum_{i=1}^m|\alpha_i-1|+\sum_{i=1}^m|\tau_i|+\sum_{r=1}^{\bar m}|\beta_r|\right),
\end{equation}
where $c$ is a small positive constant independent of $A:=(a_1, \cdots, a_m)$, $\bar\alpha=(\alpha_1, \cdots, \alpha_m)$, $\bar\l=(\l_1, \cdots, \l_m)$, $\bar \beta=(\beta_1, \cdots, \beta_{\bar m})$ and $\epsilon$. In the remaining region namely $V_4(m, \epsilon, \eta)\cap \tilde V(m, \epsilon, \eta)$, we will consider the following vector field
\begin{equation}\label{eq:tildew4}
\tilde W_4:==W_{\bar \alpha}-\Omega_4 W_{\bar \l}+W_{A},
\end{equation}
where $\Omega_4$ is a large positive constant to be chosen later. Arguing exactly as in the proof of \eqref{eq:esttildew3}, we derive that up to taken $\Omega_4$ large enough that we that for every $u:=\sum_{i=1}^m\alpha_i\varphi_{a_i, \l_i}+\sum_{r=1}^{\bar m}\beta_r (v_r-\ov{(v_r)}_Q)\in V_4(m, \epsilon, \eta)\cap \tilde V(m, \epsilon, \eta)$ with the concentration  points $a_i$,  the masses $\alpha_i$, the concentrating parameters $\l_i$ ($i=1, \cdots, m$)  and the negativity parameters $\beta_r$ ($r=1, \cdots, \bar m$) satisfying \eqref{eq:afpara}, the following estimate
\begin{equation}\label{eq:esttildew4}
<-\n J(u), \tilde{W}_{4}>\geq c\left(\sum_{i=1}^m\frac{1}{\l_i^2}+\sum_{i=1}^m\frac{|\n_g\mathcal{F}^{A}_i(a_i)|}{\l_i}+\sum_{i=1}^m|\alpha_i-1|+\sum_{i=1}^m|\tau_i|+\sum_{r=1}^{\bar m}|\beta_r|\right),
\end{equation}
where $c$ is a small positive constant independent of $A:=(a_1, \cdots, a_m)$, $\bar\alpha=(\alpha_1, \cdots, \alpha_m)$, $\bar\l=(\l_1, \cdots, \l_m)$, $\bar \beta=(\beta_1, \cdots, \beta_{\bar m})$ and $\epsilon$. Finally, we define the vector field  $W_{no}$ in $V(m, \epsilon, \eta)$ as a convex combination of the vectors fields $\hat{W}$, $\tilde{W}_1, \cdots, \tilde{W}_4$. Hence, using \eqref{eq:hatw}, \eqref{eq:esthatw1}, \eqref{eq:tildew21}, \eqref{eq:esttildew3}, and \eqref{eq:esttildew4}, we obtain that for every $u:=\sum_{i=1}^m\alpha_i\varphi_{a_i, \l_i}+\sum_{r=1}^{\bar m}\beta_r (v_r-\ov{(v_r)}_Q)\in V(m, \epsilon, \eta)$ with the concentration  points $a_i$,  the masses $\alpha_i$, the concentrating parameters $\l_i$ ($i=1, \cdots, m$)  and the negativity parameters $\beta_r$ ($r=1, \cdots, \bar m$) satisfying \eqref{eq:afpara}, the following estimate
\begin{equation}\label{eq:esttildew4}
<-\n J(u), W_{no}>\geq c\left(\sum_{i=1}^m\frac{1}{\l_i^2}+\sum_{i=1}^m\frac{|\n_g\mathcal{F}^{A}_i(a_i)|}{\l_i}+\sum_{i=1}^m|\alpha_i-1|+\sum_{i=1}^m|\tau_i|+\sum_{r=1}^{\bar m}|\beta_r|\right).
\end{equation}
 where $c$ is a small positive constant independent of  $A:=(a_1, \cdots, a_m)$, $\bar \alpha=:(\alpha_1, \cdots, \alpha_m)$, $\bar \l:=(\l_1, \cdots, \l_m)$, $\bar \beta:=(\beta_1, \cdots, \beta_{\bar m})$ and $\epsilon$, thereby ending the proof of \eqref{eq:pseudoexact}.
  On the other hand,  \eqref{eq:pseudoperturb} follows from \eqref{eq:pseudoexact} by using exactly the same arguments as \cite{bcch}, which is possible thanks, Lemma \ref{eq:energyest}, Lemma \ref{eq:expansionj}, \eqref{eq:estbarw},  Lemma \ref{eq:gradientlambdaest}, Lemma \ref{eq:gradientaest}, and to the minimality property of $\bar w$  . Next, using Lemma \ref{eq:selfintest} and Lemma \ref{eq:interactest}, we have that that the following estimate holds
 \begin{equation}\label{eq:wbdd}
 \sum_{i=1}^m||\l_i\frac{\partial \varphi_{a_i, \l_i}}{\partial \l_i}||_{W^{2, 2}(M)}+ \sum_{i=1}^m||\frac{1}{\l_i}\frac{\partial \varphi_{a_i, \l_i}}{a_i}||_{W^{2, 2}(M)}+ \sum_{i=1}^m||\frac{1}{\log\l_i} \varphi_{a_i, \l_i}||_{W^{2, 2}(M)}=O(1).
 \end{equation}
 So, using \eqref{eq:wbdd}, and the definition of $W_{no}$, it is easy to see that $W_{no}$ is a $W^{2, 2}$-bounded vector field and  that the only region where the $\l_i$'s increase along the flow lines associated to $W_{no}$ is the region $V_4(m, \epsilon, \eta)\cap \tilde{V}(m, \epsilon, \eta)$, thereby ending the proof of point 3) . Hence the proof of the Lemma is complete.
 \end{pf}
 \vspace{6pt}

 \noindent
Now, we are going to use  Proposition \ref{eq:conspseudograd} to characterize the  critical points at infinity of $J$ through the following corollary.
 \begin{cor}\label{eq:loccritinf}
 1) The  critical points at infinity of $J$ belongs to $V_4(m, \epsilon, \eta)\cap \tilde{V}(m, \epsilon, \eta)$  with  $\eta$ a small positive real number verifying $0<2\eta<\varrho$ where $\varrho$ is as in \eqref{eq:cutoff}, and $0<\epsilon\leq \epsilon_0$ where $\epsilon_0$ is as in \eqref{eq:mini} and correspond to the "configurations" $\alpha_i=1$, $\l_i=+\infty$,$\tau_i=0$ $i=1, \cdots, m$, $\beta_r=0$, $r=1, \cdots, \bar m$, $A$ is a critical point of $\mathcal{F}_K$ with $\mathcal{L}_K(A)<0$ and $V=0$.\\ Thus we have a clear and natural one to one correspondence between the set of  critical points at infinity of $J$ and the set $\{A\in Crit(\mathcal{F}_K), \mathcal{L}_K(A)<0\}$. Precisely, the correspondence consist just of sending each critical point at infinity of $J$ to the element $A\in\{B\in Crit(\mathcal{F}_K), \mathcal{L}_K(B)<0\}$ which correspond to the "configuration" defining the critical point at infinity and use the notation $A_{\infty}$ to denote the corresponding critical point at infinity. \\\\
 2) The $J$-energy of a  critical point at infinity $A_{\infty}$ denoted by $E_J(A_{\infty})$ is given by
 \begin{equation}\label{eq:infinitycriticallevel}
 E_J(A_{\infty})=-\frac{40}{3}m\pi^2-8m\pi^2\log(\frac{m\pi^2}{6})-8\pi^2\mathcal{F}_K(a_1, \dots, a_m)
 \end{equation}
 where $A=(a_1, \cdots, a_m)$.
 \end{cor}
 \begin{pf}
Point 1) follows from \eqref{eq:auxiindexa1}, Lemma \ref{eq:deformlem} combined with Proposition \ref{eq:conspseudograd} and the discussion just after \eqref{eq:expjv}, while point 2) follows from point 1) combined with \eqref{eq:coreql6}, \eqref{eq:estbarw}, and Lemma \ref{eq:expansionj}.
 \end{pf}
 \begin{rem}\label{eq:reminfinity}
 We would like to point out the following trivial fact
 $$
 \tau_i=0, \;\;\;\l_i=+\infty, \;\;i=1, \cdots, m\Longleftrightarrow \tau_i=0,\;\;i=2, \cdots, m, \l_1=+\infty.
 $$
 \end{rem}
\vspace{6pt}

 \noindent
Now, recalling that $V_4(m, \epsilon, \eta)\cap  \tilde{V}(m, \epsilon, \eta)$ is a neighborhood of potential  critical point at infinity (see Corollary \ref{eq:loccritinf})), we have that \eqref{eq:coreql7}, \eqref{eq:estbarw}, Corollary \ref{eq:c1expansionj1}, Lemma \ref{eq:expansionj},  \eqref{eq:la}  and the definition of $V_4(m, \epsilon, \eta)\cap  \tilde{V}(m, \epsilon, \eta)$(see \eqref{eq:v4}) combined with classical Morse Lemma implies the following Morse Lemma at infinity.
\begin{lem}\label{eq:morleminf}(Morse Lemma at infinity)\\
Assuming that $\eta$ is a small positive real number with $0<2\eta<\varrho$ where $\varrho$ is as in \eqref{eq:cutoff}, $0<\epsilon\leq \epsilon_0$ where $\epsilon_0$ is as in \eqref{eq:mini} and $u_0:=\sum_{i=1}^m\alpha_i^0\varphi_{a_i^0, \l_i^0}+\sum_{r=1}^{\bar m}\beta_k^0(v_r-\ov{(v_r)}_Q)+\bar w((\bar\alpha^0, A^0, \bar \l^0, \bar \beta^0))\in V_4(m, \epsilon, \eta)\cap  \tilde{V}(m, \epsilon, \eta) $ (where $\bar \alpha^0:=(\alpha^0_1, \cdots, \alpha^0_m)$, $A^0:=(a_1^0, \cdots, a^0_m)$,  $\bar \l:=(\l_1^0, \cdots, \l_m^0)$ and  $\bar \beta^0:=(a_1^0, \cdots, \beta_{\bar m}^0)$) with the concentration  points $a_i^0$,  the masses $\alpha_i^0$, the concentrating parameters $\l_i^0$ ($i=1, \cdots, m$) and the negativity parameters $\beta_r^0$ ($r=1, \cdots, \bar m$) satisfying \eqref{eq:afpara} and furthermore $A^0\in Crit(\mathcal{F}_K)$, then there exists an open neighborhood  $U$ of $(\bar\alpha^0, A^0, \bar \l^0, \bar \beta^0)$ such that for every $(\bar \alpha, A, \bar \l, \bar \beta)\in U$ with $\bar \alpha:=(\alpha_1, \cdots, \alpha_m)$, $A:=(a_1, \cdots, a_m)$,  $\bar \l:=(\l_1, \cdots, \l_m)$, $\bar \beta:=(\beta_1, \cdots, \beta_{\bar m})$, and the $a_i$,  the $\alpha_i$, the $\l_i$ ($i=1, \cdots, m$)  and the $\beta_r$ ($r=1, \cdots, \bar m$) satisfying \eqref{eq:afpara}, and $w$ satisfying \eqref{eq:afpara}  with $\sum_{i=1}^m\alpha_i\varphi_{a_i, \l_i}+\sum_{r=1}^{\bar m}\beta_r(v_r-\ov{(v_r)}_Q)+w\in V_4(m, \epsilon, \eta)\cap  \tilde{V}(m, \epsilon, \eta)$, we have the existence of a change of variable
\begin{equation}\label{eq:morsevinf}
\begin{split}
&\alpha_i\longrightarrow y_i ,i=1, \cdots, m,\\& A\longrightarrow \tilde A=(\tilde A_{-}, \tilde A_{+})\\&\l_1\longrightarrow x_1,\\&\tau_i\longrightarrow x_i, i=2, \cdots, m,\\&\beta_r\longrightarrow \tilde \beta_r\\ &V\longrightarrow \tilde V,
\end{split}
\end{equation}
such that
\begin{equation}
J(\sum_{i=1}^m\alpha_i\varphi_{a_i, \l_i}+\sum_{r=1}^{\bar m}\beta_r(v_r-\ov{(v_r)}_Q+w)=-|\tilde A_{-}|^2+|\tilde A_{+}|^2+\sum_{i=1}^{m}y_i^2-\sum_{r=1}^{\bar m}\tilde \beta_r^2+x_1^2-\sum_{i=2}^mx_i^2+||\tilde V||^2
\end{equation}
where $\tilde A=(\tilde A_{-}, \tilde A_{+})$ is the Morse variable of the map $E_J: M^m\setminus F(M^m)\longrightarrow \R$ which is defined by the right hand side of \eqref{eq:infinitycriticallevel}. Hence a critical point at infinity \;$A_{\infty}=(a_1,\cdots, a_m)_{\infty}$ of $J$ has Morse index at infinity  $M_{\infty}(A_{\infty})=i_{\infty}(A)+\bar m$.
\end{lem}

\section{Proof of the results}\label{eq:proofresults}
In this section, we present the proofs of Theorem \ref{eq:morsepoincare1}-Theorem \ref{eq:Cm}. We start with the ones of Theorem \ref{eq:morsepoincare1} and Corollary \ref{eq:existence1}. \\\\
 \begin{pfn}{ of Theorem \ref{eq:morsepoincare1} and Corollary \ref{eq:existence1}}\\
We argue by contraction, hence suppose that $J$ has no critical points and look for a contradiction. To continue, we will consider the following two cases: $\bar m=0$ and $\bar m\geq 1$.\\\\
{\bf Case 1}:\;$\bar m=0$.\\
In this case,  we first take $L>L^1$ where $L^1$ is given by Lemma \ref{eq:tophigh} and infer from the same Lemma that
\begin{equation}\label{eq:ccontractible1}
J^L\;\;\text{is contractible}.
\end{equation}
On the other hand, defining the polynomials $M(t)$ and the series $P(t)$ as follows
\begin{equation}\label{eq:morsepoly1}
M(t):=\sum_{i=0}^4m_i^1t^{i},
\end{equation}
\begin{equation}\label{eq:poincarepoly1}
P(t)=\sum_{i\geq 0}^{4}b_i(J^L)t^{i}
\end{equation}
which is clearly converging thanks to \eqref{eq:ccontractible1} and using the fact that $J$ has no critical point, Lemma \ref{eq:deformlemr}, Proposition \ref{eq:conspseudograd}, Corollary \ref{eq:loccritinf}, Lemma \ref{eq:morleminf}, and classical arguments in Morse theory (see for example the work of Bahri-Rabinowitz\cite{BR}), we get
%
\begin{equation}\label{eq:morsepoincareid1}
M(t)-P(t)=(1+t)R(t)
\end{equation}
for some series $R(t)=\sum_{i\geq 0}k_i t^{i}$ where $k_i\geq 0$ for all $i$. As already said, we have $P(t)$ is converging, but clearly \eqref{eq:ccontractible1} and \eqref{eq:poincarepoly1} implies
\begin{equation}\label{eq:pt1}
P(t)=1.
\end{equation}
Now, it is easy to see that \eqref{eq:morsepoly1} and \eqref{eq:pt1} implies that \eqref{eq:morsepoincareid1}  is equivalent to \eqref{eq:mp1} has a solution, hence reaching a contradiction in the case of Theorem \ref{eq:morsepoincare1}, thereby ending the proof of Theorem \ref{eq:morsepoincare1} in this case. On the other hand, as in classical Morse theory, equation \eqref{eq:morsepoincareid1} with $t=-1$ implies
\begin{equation}\label{eq:eulpoin}
\chi(J^L)=\sum_{A\in \mathcal{F}_{\infty}}(-1)^{i_{\infty}(A)}.
\end{equation}
Now, clearly \eqref{eq:ccontractible1}, \eqref{eq:eulpoin}, and the assumption $\sum_{A\in \mathcal{F}_{\infty}}(-1)^{i_{\infty}(A)}\neq 1$ lead to a contradiction, thereby ending the proof Theorem \ref{eq:existence1} in this case. Hence the proof of Theorem \ref{eq:morsepoincare1} and Corollary \ref{eq:existence1} in the case $\bar m=0$ is complete.\\\\
{\bf Case 2}: $\bar m\geq 1$\\
In this case, we first take $L>\max\{L^1, L_{1, \bar m}\}$ where $L^1$ is given by Lemma \ref{eq:tophigh}, $L_{1, \bar m}$ is given by Lemma \ref{eq:topnegative}, and infer from Lemma \ref{eq:tophigh} that
\begin{equation}\label{eq:ccontractible2}
J^L\;\;\text{is contractible}.
\end{equation}
Next, appealing to Lemma \ref{eq:topnegative}, we infer that
\begin{equation}\label{eq:topbary2}
J^{-L} \;\,\text{has the same homotopy type as}\;S^{\bar m-1}
\end{equation}
On the other hand, using Lemma \ref{eq:deformlem}, Corollary \ref{eq:energybddinf}, Corollary \ref{eq:energybddcrit}, the the fact that $L>L^1>2\max\{C_0^1, C_1^1\}$, we deduce that $J^{-L}$ is a deformation retract of some of its neighborhood in $J^L$, namely  $(J^L, J^{-L})$ is a good pair. Hence, using \eqref{eq:ccontractible2} and and exact sequence, we get
\begin{equation}\label{eq:homl1}
H_0(J^{L}, J^{-L})\simeq 0, \;\;H_1(J^{L}, J^{-L})\simeq 0,\;\,\text{and}\;\;H_p(J^L, J^{-L})\simeq  H_{p-1}(J^{-L})\;\;\;\forall p\in \N\setminus\{0, 1\} .
\end{equation}
Now, as above, we define the polynomial $M(t)$ and the series $P(t)$ as follows
\begin{equation}\label{eq:morsepoly2}
M(t):=\sum_{i=\bar m}^{4+\bar m}m_i^{1, \bar m}t^{i},
\end{equation}
\begin{equation}\label{eq:poincarepoly2}
P(t)=\sum_{i\geq 0}b_i(J^L, J^{-L})t^{i},
\end{equation}
where
\begin{equation}\label{eq:auxcarind}
m_i^{1, \bar m}:=card\{A\in Crit(\mathcal{F}_K): 4+\bar m-Morse(A, \mathcal{F}_K)=i\}, \;i=\bar m, \cdots, 4+\bar m,
\end{equation}
Clearly \eqref{eq:topbary2} and \eqref{eq:homl1} implies that $P(t)$ is converging. On the other hand, as above, using the fact that $J$ has no critical point, Lemma \ref{eq:deformlemr}, Proposition \ref{eq:conspseudograd}, Corollary \ref{eq:loccritinf}, Lemma \ref{eq:morleminf}, and classical arguments in Morse theory (see for example the work of Bahri-Rabinowitz\cite{BR}), we get
%
\begin{equation}\label{eq:morsepoincareid2}
M(t)-P(t)=(1+t)R(t),
\end{equation}
for some series $R(t)=\sum_{i\geq 0}k_i t^{i}$ where $k_i\geq 0$ for all $i$. Now, using \eqref{eq:topbary2}, \eqref{eq:homl1}, and \eqref{eq:poincarepoly2}, we obtain
\begin{equation}\label{eq:pt2}
P(t)=t^{\bar m}.
\end{equation}
On the other hand, it is easy to see that \eqref{eq:morsepoly2}, \eqref{eq:pt2} implies that \eqref{eq:morsepoincareid2} is equivalent to the following system
 has a solution
 \begin{equation}\label{eq:mp2}
\begin{cases}
0=k_0,\\
0=k_i+k_{i-1}, \;&i=1, \cdots, \bar m-1,\\
m^{1, \bar m}_{\bar m}=1+k_{\bar m}+k_{\bar m-1},\\
m_{i}^{1, \bar m}=k_{i}+k_{i-1}, \;\;i=\bar m+1, 4+\bar m,\\
0=k_{4+\bar m},\\
k_i\geq 0,\;\;& i=0, \cdots, 4+\bar m.
\end{cases}
\end{equation}
On the other hand, it is easy to that the following holds
\begin{equation}\label{eq:relationmibarm}
m^{1, \bar m}_i=m^1_{i-\bar m}, \;\;\;i=\bar m, \cdots, 4+\bar m.
\end{equation}
 Now, clearly \eqref{eq:relationmibarm} implies that the system \eqref{eq:mp1} and \eqref{eq:mp2} are equivalent, thus both are equivalent to \eqref{eq:morsepoincareid2}. Hence reaching a contradiction for Theorem \ref{eq:morsepoincare1} in this case. Next, as above, using equation \eqref{eq:morsepoincareid2} with $t=-1$ and the relation $\chi(J^{L}, J^{-L})=\chi(J^L)-\chi(J^{-L})$, we get
\begin{equation}\label{eq:eulpoin2}
\chi(J^L)=(-1)^{\bar m}\sum_{A\in \mathcal{F}_{\infty}}(-1)^{i_{\infty}(A)}+\chi(J^{-L}).
\end{equation}
Now, clearly \eqref{eq:ccontractible2}, \eqref{eq:topbary2}, \eqref{eq:eulpoin2}, the assumption $\sum_{A\in \mathcal{F}_{\infty}}(-1)^{i_{\infty}(A)}\neq 1$, and the well-known  relation \;$1-\chi(S^{\bar m-1})=(-1)^{\bar m}$ lead to a contradiction, thereby ending the proof of Theorem \ref{eq:existence1} in this case.
 Hence the proof of Theorem \ref{eq:morsepoincare1} and Corollary \ref{eq:existence1} is complete.
\end{pfn}
\vspace{6pt}

\noindent
\begin{pfn}{ of Theorem \ref{eq:morsepoincare2} and Corollary \ref{eq:existence2}}\\
As above, here also, we argue by contraction, hence we suppose that $J$ has no critical points and look for a contradiction. Still as above, to continue, we will consider the following two cases: $\bar m=0$ and  $\bar m\geq 1$.\\\\
{\bf Case 1}:\;$\bar m=0$.\\
In this case,  we first take $L>\max\{L^m, L_{m, 0}\}$ where $L^m$ is given by Lemma \ref{eq:tophigh}, $L_{m, 0}$ is given by Lemma \ref{eq:topnegative} and infer from the same Lemma that
\begin{equation}\label{eq:ccontractible3}
J^L\;\;\text{is contractible}.
\end{equation}
Next, appealing to Lemma \ref{eq:topnegative}, we infer that
\begin{equation}\label{eq:topbary3}
J^{-L} \;\,\text{has the same homotopy type as}\;B_{m-1}(M).
\end{equation}
On the other hand, using gain Lemma \ref{eq:deformlem}, Corollary \ref{eq:energybddinf}, Corollary \ref{eq:energybddcrit}, the the fact that $L>L^m>2\max\{C_0^m, C_1^m\}$, we deduce that $J^{-L}$ is a deformation retract of some of neighborhood of it in $J^L$, hence  $(J^L, J^{-L})$ is a good pair. Therefore,  as above, using \eqref{eq:ccontractible3} and exact sequence, we get
\begin{equation}\label{eq:homl3}
H_0(J^{L}, J^{-L})\simeq 0, \;\;\;H_1(J^{L}, J^{-L})\simeq 0, \;\;\;\;\text{and}\;\;\; H_p(J^L, J^{-L})\simeq H_{p-1}(J^{-L})\;\;\;\forall p\in \N\setminus\{0, 1\}.
\end{equation}
On the other hand, defining the polynomials $M(t)$ and the series $P(t)$ as follows
\begin{equation}\label{eq:morsepoly3}
M(t):=\sum_{i=0}^{5m-1}m_i^mt^{i},
\end{equation}
\begin{equation}\label{eq:poincarepoly3}
P(t)=\sum_{i\geq 0}b_i(J^L, J^{-L})t^{i},
\end{equation}
which is clearly converging thanks to  \eqref{eq:homl3} and \eqref{eq:topbary3} and and using the fact that $J$ has no critical point, Lemma \ref{eq:deformlem}, Proposition \ref{eq:conspseudograd}, Corollary \ref{eq:loccritinf}, Lemma \ref{eq:morleminf}, and classical arguments in Morse theory (see for example the work of Bahri-Rabinowitz\cite{BR}), we get
%
\begin{equation}\label{eq:morsepoincareid3}
M(t)-P(t)=(1+t)R(t),
\end{equation}
for some series $R(t)=\sum_{i\geq 0}k_i t^{i}$ where $k_i\geq 0$ for all $i$. As already said, we have $P(t)$ is converging, but clearly  \eqref{eq:cpm}, \eqref{eq:homl3}, \eqref{eq:topbary3} and \eqref{eq:poincarepoly3} implies
\begin{equation}\label{eq:pt3}
P(t)=\sum_{i=2}^{5m-6}c_{i-1}^{m-1}t^{i}.
\end{equation}
Now, it is easy to see that \eqref{eq:morsepoly3} and \eqref{eq:pt3} implies that \eqref{eq:morsepoincareid3}  is equivalent to \eqref{eq:mp3} has a solution, hence reaching a contradiction in the case of Theorem \ref{eq:morsepoincare2}, thereby ending the proof of Theorem \ref{eq:morsepoincare2} in this case. Now, as above, using equation \eqref{eq:morsepoincareid3} with $t=-1$  combined with the relation $\chi(J^L, J^{-L})=\chi(J^L)-\chi(J^{-L})$, we get
\begin{equation}\label{eq:eulpoin3}
\chi(J^L)=\sum_{A\in \mathcal{F}_{\infty}}(-1)^{i_{\infty}(A)}+\chi(J^{-L}).
\end{equation}
Now, clearly \eqref{eq:ccontractible3}, \eqref{eq:topbary3}, \eqref{eq:eulpoin3}, the assumption $\sum_{A\in \mathcal{F}_{\infty}}(-1)^{i_{\infty}(A)}\neq \frac{1}{(m-1)!}\Pi_{i=1}^{m-1}(i-\chi(M))$ and the well-known fact $1-\chi(B_{m-1}(M))=\frac{1}{(m-1)!}\Pi_{i=1}^{m-1}(i-\chi(M))$ (see for example \cite{maldeg}) lead to a contradiction, thereby ending the proof Theorem \ref{eq:existence2} in this case. Hence the proof of Theorem \ref{eq:morsepoincare2} and Theorem \ref{eq:existence2} in the case $\bar m=0$ is complete.\\\\
{\bf Case 2}: $\bar m\geq 1$\\
In this case, we first take $L>\max\{L^m, L_{m, \bar m}\}$ where $L^m$ is given by Lemma \ref{eq:tophigh}, $L_{m, \bar m}$ is given by Lemma \ref{eq:topnegative}, and infer from Lemma \ref{eq:tophigh} that
\begin{equation}\label{eq:ccontractible4}
J^L\;\;\text{is contractible}.
\end{equation}
Next, appealing to Lemma \ref{eq:topnegative}, we infer that
\begin{equation}\label{eq:topbary4}
J^{-L} \;\,\text{has the same homotopy type as}\;A_{m-1, \bar m}.
\end{equation}
On the other hand, using Lemma \eqref{eq:deformlem}, \eqref{eq:energybddinf}, \eqref{eq:energybddcrit}, the the fact that $L>L^m>2\max\{C_0^m, C_1^m\}$, we deduce that $J^{-L}$ is a deformation retract of some of neighborhood of it in $J^L$, hence  $(J^L, J^{-L})$ is a good pair. Thus, as above, using \eqref{eq:ccontractible4} and exact sequence, we derive that here also,
\begin{equation}\label{eq:homl4}
H_0(J^{L}, J^{-L})\simeq 0, \;\;\;H_1(J^{L}, J^{-L})\simeq 0\;\;\;\;\text{and}\;\;\; H_p(J^L, J^{-L})\simeq H_{p-1}(J^{-L})\;\;\;\forall p\in \N\setminus\{0, 1\}.
\end{equation}
Now, as above, we define the polynomial $M(t)$ and the series $P(t)$ as follows
\begin{equation}\label{eq:morsepoly4}
M(t):=\sum_{i=\bar m}^{5m-1+\bar m}m_i^{m, \bar m}t^{i},
\end{equation}
\begin{equation}\label{eq:poincarepoly4}
P(t)=\sum_{i\geq 0}b_i(J^L, J^{-L})t^{i},
\end{equation}
where
\begin{equation}\label{eq:auxcarind1}
m_i^{m, \bar m}:=card\{A\in Crit(\mathbb{•}hcal{F}_K): 5m-1+\bar m-Morse(A, \mathcal{F}_K)=i\}, i=\bar m, \cdots, 5m-1+\bar m,
\end{equation}
with clearly $m_i^{m, \bar m}=0$ for $\bar m\leq i\leq m-2+\bar m$. Now, we have  \eqref{eq:topbary4} and \eqref{eq:homl4} implies that $P(t)$ is converging. On the other hand, as above, using the fact that $J$ has no critical point, Lemma \ref{eq:deformlem}, Proposition \ref{eq:conspseudograd}, Corollary \ref{eq:loccritinf}, Lemma \ref{eq:morleminf}, and classical arguments in Morse theory (see for example the work of Bahri-Rabinowitz\cite{BR}), we get
\begin{equation}\label{eq:morsepoincareid4}
M(t)-P(t)=(1+t)R(t),
\end{equation}
for some series $R(t)=\sum_{i\geq 0}k_i t^{i}$ where $k_i\geq 0$ for all $i$. Now, using \eqref{eq:relationbary}, \eqref{eq:topbary4}, \eqref{eq:homl4}, \eqref{eq:poincarepoly4}, we obtain
\begin{equation}\label{eq:pt4}
P(t)=\sum_{i=\bar m+2}^{5m-5+\bar m}c_{i-1-\bar m}^{m-1}t^{i}.
\end{equation}
Now, it is easy to see that \eqref{eq:morsepoly4}, \eqref{eq:pt4} implies that \eqref{eq:morsepoincareid4} is equivalent to the following system
 has a solution
 \begin{equation}\label{eq:mp4}
\begin{cases}
0=k_i,\;&i=0, \cdots, \bar m\\
m^{m, \bar m}_{\bar m+1}==k_{1+\bar m}, \\
m_{i}^{m, \bar m}=c_{i-1-\bar m}^{m-1}+k_{i}+k_{i-1}, \;\;&i=\bar m+2, \cdots 5m-5+\bar m,\\
m^{m, \bar m}_i=k_i+k_{i-1},\;\;;& i=5m-4+\bar m,\cdots, 5m-1+\bar m,\\
0=k_{5m-1+\bar m},\\
k_i\geq 0,\;\;& i=0, \cdots, 5m-1+\bar m.
\end{cases}
\end{equation}
On the other hand, it is easy to that the following holds
\begin{equation}\label{eq:relationmibarm2}
m^{m, \bar m}_i=m^m_{i-\bar m}, \;\;\;i=\bar m, \cdots, 5m-1+\bar m.
\end{equation}
 Now, clearly \eqref{eq:relationmibarm2} implies that the system \eqref{eq:mp3} and \eqref{eq:mp4} are equivalent. Hence reaching a contradiction for Theorem \ref{eq:morsepoincare2} in this case. Next, as above, using equation \eqref{eq:morsepoincareid4} with $t=-1$ and the relation $\chi(J^{L}, J^{-L})=\chi(J^L)-\chi(J^{-L})$, we get
\begin{equation}\label{eq:eulpoin4}
\chi(J^L)=(-1)^{\bar m}\sum_{A\in \mathcal{F}_{\infty}}(-1)^{i_{\infty}(A)}+\chi(J^{-L}).
\end{equation}
Now, clearly \eqref{eq:ccontractible4}, \eqref{eq:topbary4}, \eqref{eq:eulpoin4}, the assumption $\sum_{A\in \mathcal{F}_{\infty}}(-1)^{i_{\infty}(A)}\neq 1$, and the well-known  relation \;$1-\chi(A_{m-1, \bar m})=(1-\chi(S^{\bar m-1}))(1-\chi(B_{m-1}(M)))=(-1)^{\bar m}(1-\chi(B_{m-1}(M)))=(-1)^{\bar m}\frac{1}{(m-1)!}\Pi_{i=1}^{m-1}(i-\chi(M))$ lead to a contradiction, thereby ending the proof of Corollary \ref{eq:existence2} in this case.
 Hence the proof of Theorem \ref{eq:morsepoincare2} and Corollary\ref{eq:existence2} is complete.
\end{pfn}

\begin{pfn}{ of Theorem \ref{t:C} and Theorem \ref{eq:Cm}}\\
 We first notice that, thanks to our Morse decomposition (Lemma \ref{eq:morleminf}), we can associate to each "critical point at infinity" $A_{\infty}$ a stable $W_s(A_{\infty})$, as well as an unstable manifold $W_u (A_{\infty})$. Moreover, we say that a critical point at infinity $A_{\infty}$ dominates another critical point at infinity $B_{\infty}$ if
 $$
 W_u(A_{\infty}) \cap W_s (B_{\infty})  \neq \emptyset.
 $$
Assuming that $A_{\infty}$ dominates $B_{\infty}$ and that the intersection is transverse, and recalling that $M_{\infty}(A_{\infty})=i_{\infty}(A)+\bar m$, we have   then that
$$
i_{\infty}(A) \, \geq i_{\infty}(B).
$$
We first prove Theorem \ref{t:C}. To do that, for $l \in \N$ as in Theorem \ref{t:C}, we define the following closed set
$$
X_{\infty} := \bigcup_{ A \in \mathcal{F}_{\infty}:\;  i_{\infty}(A) \leq l - 1} \ov{W_u(A_{\infty})}= \bigcup_{ A \in \mathcal{F}_{\infty}:\; M_{\infty}(A_{\infty}) \leq l - 1+\bar m} \ov{W_u(A_{\infty})}.
$$
By Proposition 4.21 of Bahri-Rabinowitz  \cite{BR}, the closure of $W_u(A_{\infty})$ adds to it the unstable manifold of all critical points or critical points at infinity dominated by $A_{\infty}$. Arguing by contradiction, we assume that $J$ has no critical point. Denoting by $A_0$ a global maximum of $\mathcal{F}_K$ and $(A_0)_{\infty}$ its associated critical point at infinity, we define the following set 
$$
\mathcal{C}(X_{\infty}) := \big \{  u + (1-t) (A_0)_{\infty}: \quad u \in X_{\infty},\; t\in [0, 1]  \big \}.
$$
It follows from the work of Bahri-Rabinowitz \cite{BR}(sections 7 and 8) that $\mathcal{C}(X_{\infty})$ is a stratified set of top dimension $l+\bar m$. Now, we use the pseudogradient in Proposition \ref{eq:conspseudograd} to deform the set $\mathcal{C}(X_{\infty})$. By transversality arguments, we may assume that the deformation avoids all "critical points at infinity" of Morse indices at infinity bigger than or equal to $l +1+\bar m$. Since, by assumption, there are no critical points at infinity of Morse index at infinity $l+\bar m$, we conclude that $\mathcal{C}(X_{\infty})$ is contractible. Hence, by the Poincar\'e-Hopf theorem, we derive that 
\begin{equation}\label{eq:con1}
1 \, = \, \chi(\mathcal{C}(X_{\infty})) \, = \, \sum_{A\in \mathcal{F}_{\infty}:\; i_{\infty}(A) \leq l -1} (-1)^{i_{\infty}(A)} 
\end{equation}
for $\bar m=0$ and 
\begin{equation}\label{eq:con2}
1 \, = \, \chi(\mathcal{C}(X_{\infty})) \, = \, \sum_{ A\in \mathcal{F}_{\infty}:\;i_{\infty}(A) \leq l -1} (-1)^{i_{\infty}(A)+\bar m} +\chi(J^{-L})\
\end{equation}
for $\bar m\geq 1$ with  $L$ as in the proof of Theorem \ref{eq:morsepoincare1} and Corollary \ref{eq:existence1}. Now clearly \eqref{eq:con1} contradicts the assumption of the theorem in the case $\bar m=0$.  Furthermore, in the case $\bar m\geq 1$ and for the same reason as in of Theorem \ref{eq:morsepoincare1} and Corollary \ref{eq:existence1}, we have  \eqref{eq:con2} contradicts the assumption of the theorem in the case $\bar m\geq 1$. Hence the proof of Theorem \ref{t:C} is complete. \\
Now, following the same arguments one deduces Theorem \ref{eq:Cm}. Indeed let  $X_{\infty}$  and $\mathcal{C}(X_{\infty})$ be as above where $A_0$ is replaced by $A^l$.  Just like in the above proof, we derive that
\begin{equation}\label{eq:con3}
1 \, = \, \chi(\mathcal{C}(X_{\infty})) \, = \, \sum_{A\in \mathcal{F}_{\infty}:\; i_{\infty}(A) \leq l -1}(-1)^{i_{\infty}(A)+\bar m} \, + \, \chi(J^{-L}),
\end{equation}
with  $L$ as in the proof of Theorem \ref{eq:morsepoincare2} and Corollary \ref{eq:existence2}. Thus, for the same reason as in of Theorem \ref{eq:morsepoincare2} and Corollary \ref{eq:existence2}, we have  \eqref{eq:con3} contradicts the assumption of the theorem. Hence the proof of Theorem \ref{eq:Cm} is complete. 
\end{pfn}

\section{Appendix}\label{eq:appendix}
In this section, we collect some technical Lemmas. First of all, using the conformal invariance of the Paneitz operator, the properties of the metric $g_{a}$ (see \eqref{eq:detga}-\eqref{eq:proua}), the PDE satisfied by the $\varphi_{a, \l}$'s (see \eqref{eq:projbubble}),  and the Green's representation formula \eqref{eq:G4integral}, we derive the following two Lemmas.
\begin{lem}\label{eq:intbubbleest}(Global Bubbles estimates)\\
Assuming that $\epsilon$ is positive and small,  $a\in M$ and $\l\geq \frac{1}{\epsilon}$, then\\
1)
$$
\varphi_{a, \l}(\cdot)=\hat \d_{a, \l}(\cdot)+\log\frac{\l}{2}+H(a, \cdot)+\frac{1}{4\l^2}\D_{g_{a}}H(a, \cdot)+O\left(\frac{1}{\l^3}\right),
$$
2)
$$
\l\frac{\partial \varphi_{a, \l}(\cdot)}{\partial \l}=\frac{2}{1+\l^2\chi_{\varrho}^2(d_{g_{a}}(a, \cdot))}-\frac{1}{2\l^2}\D_{g_{a}}H(a, \cdot)+O\left(\frac{1}{\l^3}\right),
$$
3)
$$
\frac{1}{\l}\frac{\partial \varphi_{a, \l}(\cdot)}{\partial a}=\frac{\chi_{\varrho}(d_{g_{a}}(a, \cdot))\chi_{\varrho}^{'}((d_{g_{a}}(a, \cdot))}{d_{g_{a}}(a, \cdot)}\frac{2\l exp_a^{-1}(\cdot)}{1+\l^2\chi_{\varrho}^2(d_{g_a}(a, \cdot))}+\frac{1}{\l}\frac{\partial  H(a, \cdot)}{\partial a}+O\left(\frac{1}{\l^3}\right),
$$
where \;$O(1)$ means $O_{a, \l, \epsilon}(1)$ and for it meaning see section \ref{eq:notpre}.
\end{lem}
\begin{lem}\label{eq:outbubbleest}(Global Bubbles estimates outside concentration points)\\
Assuming that  $\epsilon$ is small and d positive, $a\in M$, $\l\geq \frac{1}{\epsilon}$, and \;$0<2\eta<\varrho$ with $\varrho$ as in \eqref{eq:cutoff}, then  there holds
$$
\varphi_{a, \l}(\cdot)=G(a, \cdot)+\frac{1}{4\l^2}\D_{g_{a}}G(a, \cdot)+O\left(\frac{1}{\l^3}\right)\;\;\text{in} \;\;\;M\setminus B^{a}_{a}(\eta),
$$
$$
\l\frac{\partial \varphi_{a, \l}(\cdot)}{\partial \l}=-\frac{1}{2\l^2}\D_{g_{a}}G_(a, \cdot)+O\left(\frac{1}{\l^3}\right)\;\;\text{in} \;\;\;M\setminus B^{a}_{a}(\eta),
$$
and
$$
\frac{1}{\l}\frac{\partial \varphi_{a, \l}(\cdot)}{\partial a}=\frac{1}{\l}\frac{\partial G(a, \cdot)}{\partial a}+O\left(\frac{1}{\l^3}\right)\;\;\text{in} \;\;\;M\setminus B^{a}_{a}(\eta),
$$
where $O(1)$ means $O_{a, \l, \epsilon}(1)$ and for it meaning see section \ref{eq:notpre}.
\end{lem}
\vspace{4pt}

\noindent
Next, using the above two Lemmas, we obtain the following three Lemmas.
\begin{lem}\label{eq:selfintest}(Selfaction estimates)\\
Assuming that $\epsilon$ is small and positive, $a\in M$ and $\l\geq \frac{1}{\epsilon}$, then there holds
$$
<P_g\varphi_{a, \l}, \varphi_{a, \l}>=32\pi^2\log \l-\frac{40\pi^2}{3}+16\pi^2 H(a, a)+\frac{8\pi^2}{\l^2}\D_{g_a} H(a, a)+O\left(\frac{1}{\l^3}\right),
$$
$$
<P_g\varphi_{a, \l}, \l\frac{\varphi_{a, \l}}{\partial \l} >=16\pi^2-\frac{8\pi^2}{\l^2}\D_{g_a}H(a, a)+O\left(\frac{1}{\l^3}\right),
$$
$$
<P_g\varphi_{a, \l}, \frac{1}{\l}\frac{\varphi_{a, \l}}{\partial a} >=\frac{16\pi^2}{\l}\frac{\partial H(a, a)}{\partial a}+O\left(\frac{1}{\l^3}\right),
$$
where \;$O(1)$\; means \;$O_{a, \l, \epsilon}(1)$\; and for its meaning see section \ref{eq:notpre}.
\end{lem}

\begin{lem}\label{eq:interactest} (Interaction estimates)\\
Assuming that $\epsilon$ is small and positive $a_i, a_j\in M$, \;$d_g(a_i, a_j)\geq 4\ov C\eta$, $0<2\eta<\varrho$, $\frac{1}{\L}\leq \frac{\l_i}{\l_j}\leq \L$, and $\l_i, \l_j\geq \frac{1}{\epsilon}$, $\ov C$ as in \eqref{eq:proua}, and $\varrho$ as in \eqref{eq:cutoff}, then there holds
$$
<P_g\varphi_{a_i, \l_i}, \varphi_{a_j, \l_j}>=16\pi^2G (a_j, a_i)+\frac{8\pi^2}{\l_j^2}\D_{g_{a_j}} G(a_j, a_i)+O\left(\frac{1}{\l^3}\right),
$$
$$
<P_g\varphi_{a_i, \l_i}, \l_j\frac{\partial \varphi_{a_j, \l_j}}{\partial \l_j}>=-\frac{8\pi^2}{\l_j^2}\D_{g_{a_j}}G(a_j, a_i)+O\left(\frac{1}{\l^3}\right),
$$
and
$$
<P_g\varphi_{a_i, \l_i}, \frac{1}{\l_j}\frac{\partial \varphi_{a_j, \l_j}}{\partial a_j}>= \frac{16\pi^2}{\l_j}\frac{\partial G(a_j, a_i)}{\partial a_j}+O\left(\frac{1}{\l^3}\right),
$$
where \;$O(1)$\; means here \;$O_{A, \bar \l, \epsilon}(1)$\; with \;$A=(a_i, a_j)$\; and \;$\bar \l=(\l_i, \l_j)$\; and for the meaning of $O_{A, \bar \l, \epsilon}(1)$, see section \ref{eq:notpre}.
\end{lem}
\begin{lem}\label{eq:auxibubbleest1}
Assuming that $\epsilon>0$ is very small, we have that for $a\in M$, $\l\geq \frac{1}{\epsilon}$, there holds
\begin{equation}
||\l\frac{\partial \varphi_{a, \l}}{\partial \l}||=\tilde O(1),
\end{equation}
\begin{equation}
||\frac{1}{\l}\frac{\partial \varphi_{a, \l}}{\partial a}||=\tilde O(1),
\end{equation}
and
\begin{equation}
||\frac{1}{\sqrt{\log \l}}\varphi_{a, \l}||=\tilde O(1),
\end{equation}
where here $\tilde O(1)$ means bounded by positive constants form below and above independent of $\epsilon$, $a$, and $\l$.
\end{lem}
\vspace{4pt}

\noindent
Finally, using the relation between $(\R^4, g_{\R^4})$ and $(S^4, g_{S^4})$,  the spectral property of the Paneitz operator $P_{g_{S^4}}$, and the blow-up argument of Brendle \cite{bren2}, we obtain the following last two Lemmas.
\begin{lem}\label{eq:positive}
There exists \;$\Gamma_0$\; and \;$\tilde \L_0$ two large positive constant such that for every $a\in M$, $\l\geq \tilde \L_0$, and \;$w\in F_{a,  \l}:=\{w\in W^{2, 2}(M), <Q_g, w>=<\varphi_{a, \l}, w>_P=\}$, we have
\begin{equation}\label{eq:positivity}
 \int_Me^{4\hat\d_{a,\l}}w^2dV_{g_{a}}\leq \Gamma_0||w||^2.
\end{equation}
\end{lem}

\begin{lem}\label{eq:positiveg}
Assuming that\;$\eta$\; is a small positive real number with $0<2\eta<\varrho$ where $\varrho$ is as in \eqref{eq:cutoff}, then there exists a small positive constant $c_0:=c_0(\eta)$ and $\L_0:=\L_0(\eta)$ such that for every $a_i\in M$ concentrations points  with $d_g(a_i, a_j)\geq 4\ov C\eta$ where $\bar C$ is as in \eqref{eq:proua}, for every \;$\l_i>0$\; concentrations parameters satisfying $\l_i\geq \L_0$, with $i=1, \cdots, m$,  and for every \;$w\in E_{A, \bar \l}^*=\cap_{i=1}^m E^*_{a_i, \l_i}$ with $A:=(a_1, \cdots, a_m$), $\bar \l:= (\l_1, \cdots, \l_m)$ and $E^*_{a_i, \l_i}=\{w\in W^{2, 2}(M): \;\;<\varphi_{a_i, \l_i}, w>_P=<\frac{\partial\varphi_{a_i, \l_i}}{\partial \l_i}, w>_P=<\frac{\partial\varphi_{a_i, \l_i}}{\partial a_i}, w>_P=<w, Q_g>=0\}$, there holds
\begin{equation}\label{eq:positivity}
||w||^2-48 \sum_{i=1}^m\int_Me^{4\hat\d_{a_i,\l_i}}w^2dV_{g_{a_i}}\geq c_0||w||^2,
\end{equation}
\end{lem}

\end{document}